\documentclass{article}

\usepackage{microtype}
\usepackage{graphicx}
\usepackage{subfigure}
\usepackage{booktabs} 

\usepackage{hyperref}

 \usepackage[accepted]{icml2023}

\usepackage{amsmath}
\usepackage{amssymb}
\usepackage{mathtools}
\usepackage{amsthm}

\usepackage{hyperref}       
\usepackage{url}            
\usepackage{booktabs}       
\usepackage{amsfonts}       
\usepackage{nicefrac}       
\usepackage{microtype}      
\usepackage{xcolor}   

\usepackage{amsmath}
\usepackage{amssymb}
\usepackage{graphicx}
\usepackage{multirow}
\usepackage{enumerate}
\usepackage{arydshln}
\usepackage{multicol}
\usepackage[utf8]{inputenc} 
\usepackage[T1]{fontenc}    
\usepackage{hyperref}       

\usepackage{url}            
\usepackage{booktabs}       
\usepackage{amsfonts}       
\usepackage{nicefrac}       
\usepackage{microtype}      
\usepackage{xcolor}         
\usepackage{wrapfig}
\usepackage{multicol}

\usepackage[capitalize,noabbrev]{cleveref}

\theoremstyle{plain}
\newtheorem{theorem}{Theorem}[section]
\newtheorem{proposition}[theorem]{Proposition}
\newtheorem{lemma}[theorem]{Lemma}

\theoremstyle{definition}

\newtheorem{assumption}[theorem]{Assumption}
\theoremstyle{remark}
\newtheorem{remark}[theorem]{Remark}

\def\x{\mathbf{x}}
\def\y{\mathbf{y}}
\def\z{\mathbf{z}}

\def\v{\mathbf{v}}

\def\I{\mathbf{I}}

\def\S{\mathcal{S}}

\def\A{\mathbf{A}}

\usepackage[textsize=tiny]{todonotes}

\icmltitlerunning{
Averaged Method of Multipliers for Bi-Level Optimization without Lower-Level Strong Convexity
}

\begin{document}

\twocolumn[
\icmltitle{ Averaged Method of Multipliers for Bi-Level Optimization \\
	without Lower-Level Strong Convexity}



\icmlsetsymbol{equal}{*}

\begin{icmlauthorlist}
\icmlauthor{Risheng Liu}{DUT}
\icmlauthor{Yaohua Liu}{DUT2}
\icmlauthor{Wei Yao}{sustech}
\icmlauthor{Shangzhi Zeng}{uvic}
\icmlauthor{Jin Zhang}{sustech,sustech2}
%
\end{icmlauthorlist}

\icmlaffiliation{DUT}{International School of Information Science \& Engineering, Dalian University of Technology, Dalian, China}

\icmlaffiliation{DUT2}{School of Software Technology, Dalian University of Technology}

\icmlaffiliation{sustech}{Department of Mathematics and National Center for Applied Mathematics Shenzhen, Southern University of Science and Technology, Shenzhen, China}

\icmlaffiliation{uvic}{Department of Mathematics and Statistics, University of Victoria, Victoria, British Columbia, Canada}

\icmlaffiliation{sustech2}{SUSTech International Center for Mathematics, Southern University of Science and Technology, Shenzhen, China}

\icmlcorrespondingauthor{Jin Zhang}{zhangj9@sustech.edu.cn}
%
%

\icmlkeywords{Machine Learning, ICML}

\vskip 0.3in
]



\printAffiliationsAndNotice{}  

\begin{abstract}
Gradient methods have become mainstream techniques for Bi-Level Optimization (BLO) in learning fields. The validity of existing works heavily rely on either a restrictive Lower- Level Strong Convexity (LLSC) condition or on solving a series of approximation subproblems with high accuracy or both. In this work, by averaging the upper and lower level objectives, we propose a single loop Bi-level Averaged Method of Multipliers (sl-BAMM) for BLO that is simple yet efficient for large-scale BLO and gets rid of the limited LLSC restriction. We further provide non-asymptotic convergence analysis of sl-BAMM towards KKT stationary points, and the comparative advantage of our analysis lies in the absence of strong gradient boundedness assumption, which is always required by others. Thus our theory safely captures a wider variety of applications in deep learning, especially where the upper-level objective is quadratic w.r.t. the lower-level variable. Experimental results demonstrate the superiority of our method.
\end{abstract}

\section{Introduction}
\label{intro}

Bi-Level Optimization (BLO) has recently attracted growing interests due to its wide applications in the field of deep learning, especially for hyperparameter optimization \cite{pedregosa2016hyperparameter,franceschi2017forward,franceschi2018bilevel,okuno2018hyperparameter,mackay2018self}, meta learning~\cite{franceschi2018bilevel,zugner2018adversarial,rajeswaran2019meta,KaiyiJi2020ConvergenceOM}, neural architecture search~\cite{liu2018darts,liang2019darts+,chen2019progressive,elsken2020meta}, adversarial learning~\cite{goodfellow2014generative,pfau2016connecting}, and reinforcement learning~\cite{yang2019provably,hong2020two}, etc. 
BLO tackles nested optimization structures appearing in applications. In the last decade, BLO has emerged as a prevailing optimization technique for modern machine learning tasks with an underlying hierarchy. 

BLO is notoriously challenging due to its nested nature. A variety of applications in deep learning have greatly promoted the development of BLO, with a broad collection of methods and algorithms. Among them, due to the effectiveness and simplicity, gradient-based algorithms have become mainstream techniques for BLO in learning and vision fields \cite{liu2021investigating}, which can be generally categorized into the iterative differentiation (ITD) based approach \cite{maclaurin2015gradient,franceschi2017forward,shaban2019truncated,grazzi2020iteration,liu2020generic,li2020improved,liu2021value,ji2021bilevel,ji2022will,liu2022optimization,liu2023hierarchical,liu2023augmenting} and the approximate implicit differentiation (AID) based approach \cite{pedregosa2016hyperparameter,ghadimi2018approximation,rajeswaran2019meta,lorraine2020optimizing,hong2020two,huang2021enhanced,ji2021bilevel,ji2022will,arbel2022amortized}. 
Despite the large literature, most of existing methods for BLO are slow and unsatisfactory in various ways. For example, most of these studies (including both algorithmic design and theoretical investigation) heavily rely on a restrictive Lower-Level Strong Convexity (LLSC) condition, 
except a few attempts recently (which will be discussed later on). 
However, it is observed that these conventional gradient-based algorithms may lead to incorrect solutions in the absence of LLSC; see, e.g., Example 1 in \cite{liu2020generic} and the experimental results in \cite{liu2020generic,NeurIPS2022-Liu}. 
Hence, a practical issue lingers: {\it how to design a fast BLO algorithm, which enjoys both wider applicability and provable convergence guarantee, especially without LLSC?}

The studies beyond LLSC are rather limited. 
Several recent works partially addressed the above issue by ITD-based approach using a bi-level sequential averaging method \cite{sabach2017first,liu2020generic,li2020improved} or using a bi-level value function based penalty and barrier methods \cite{liu2021value,liu2021value2}. Without using ITD-based approach, recent studies \cite{NeurIPS2022-Liu,sow2022constrained} proposed two kinds of first-order BLO algorithms that depend only on first-order gradient information, based on the value function approach for BLO \cite{outrata1990numerical,ye1995optimality}. However, all of these approaches take a double-loop optimization structure (involving numerical Lower-Level (LL) optimization loops), which heavily relies on solving a series of approximation subproblems with high accuracy. While double-loop algorithms have many applications across deep learning, they suffer from time and memory complexity, and could be difficult to implement in practice as the number of required inner loop iterations are difficult to adjust. 

\subsection{Main Contributions}\label{Our Motivations and Contributions}

The focus of this paper is to develop a provably faster single-loop BLO algorithm beyond LLSC. 
Such extension of the existing single-loop BLO methods with LLSC is not straightforward, requiring introducing the stationarity measure as stopping criterion in an appropriate way. We start from reformulating BLO to an equality-constrained optimization using the first-order stationarity condition of LL problem. Clearly, when LL problem is unconstrained and convex, this reformulation keeps the equivalence. Our {\bf first insight} is that KKT condition of the resulting equality-constrained optimization can be regarded as an appropriate stationarity measure for BLO in the absence of LLSC. 
Observe that the equality constraint in the resulting problem comes from a convex or even strongly convex unconstrained LL problem. 
Our {\bf second insight} is to well exploit the bi-level structure of the resulting equality-constrained problem, keeping the update of LL variable by gradient descent of LL objective and its variants. Notice that this feature is distinctively different from MPEC approach \cite{kim2020mpec} and the penalty method for BLO in \cite{mehra2021penalty} where the bi-level structure disappeared after reformulation. 
Our {\bf third insight} is to introduce the dual multiplier associated to the resulting equality-constrained optimization as a new variable in the update of BLO algorithm, which will enable us to develop a provably faster single-loop BLO algorithm in the absence of LLSC.
The main contributions of this paper can be summarized as follows:
\begin{itemize}
	
	\item We first propose a single loop Bi-level Averaged Method of Multipliers (sl-BAMM) beyond LLSC. The striking feature of si-BAMM is simple yet highly efficient for solving large-scale BLO, and the comparative advantage lies in theoretically guaranteed convergence. 
	
	\item We use KKT stationarity as the trackable  measure of solution quality returned by sl-BAMM. Further, a comprehensive convergence rate analysis for sl-BAMM characterized in terms of KKT stationarity is conducted. Compared to exiting works which are either lack of convergence towards stationarity or convergence rate, our result is evidently stronger and general enough to cover BLO with multiple LL minimizers.
	
	\item Finally, we are the first to conduct theoretical analysis for general BLO in the absence of the strong boundedness assumption of gradients
	of UL objective w.r.t. LL variable,  which was prevailing among existing literature. Hence, our result captures a wider range of
	applications in deep learning, especially where UL objective is quadratic w.r.t. LL variable.

\end{itemize}

\subsection{Related Work}
\label{relatedwork}

In this subsection we give a brief review of some recent works that are directly related to ours and refer to Appendix~\ref{appendix1} for a detailed comparison of these methods.

{\bf BLO Beyond Lower-Level Strong Convexity (LLSC).} 
The studies beyond LLSC are rather limited. To our knowledge, there are three main types of study used to get rid of LLSC. The first one is the sequential averaging method (SAM) used in \cite{sabach2017first,liu2020generic,li2020improved}. In the setting of BLO, for fixed UL variable, SAM iteratively generates a sequence LL points by averaging the gradient descents of UL and LL objectives. The second one is based on the value function approach for BLO, such as \cite{liu2021value,NeurIPS2022-Liu,sow2022constrained} where they all need numerical LL optimization loops to approximate the value function accurately. The last one is replacing LL problem by its first-order stationarity condition that was used in \cite{mehra2021penalty,liu2021towards} and be called MPEC approach \cite{kim2020mpec} when LL problem is constrained. However, all of these approaches take a double-loop optimization structure (involving numerical LL optimization loops), which could be difficult to implement in practice as the number of required inner loop iterations are difficult to adjust. In contrast, our algorithm sl-BAMM  has a single loop structure, which admits a much simpler implementation.  Moreover, the updates in sl-BAMM evolve simultaneously and then can be performed in parallel instead of sequentially.

{\bf Stationarity Measure for BLO Algorithm.} 
The most common stationarity measure for BLO is the gradient norm of the total UL objective (i.e., the norm of hypergradient) \cite{pedregosa2016hyperparameter,ghadimi2018approximation,chen2021closing,ji2021bilevel,ji2022will,li2022fully,dagreou2022framework,arbel2022amortized}. 
Clearly, the norm of hypergradient heavily rely on LLSC. For general BLO, the stationarity measure has been much less studied, from an algorithmic point of view, especially in the context of machine learning. Recently, stationarity measures involving value function or its regularized variant have been proposed in \cite{NeurIPS2022-Liu} and \cite{sow2022constrained} respectively. 




{\bf Discussion on Gradient Boundedness Assumption.} 
The most of existing analyses in BLO literature reply on a strong assumption on the boundedness of gradients of UL objective w.r.t. LL variable (a.k.a. gradient boundedness assumption, GBA), to guarantee the smoothness of the hypergradient and then to get tight convergence rates \cite{ghadimi2018approximation,hong2020two,ji2021bilevel,chen2021closing,li2022fully,NeurIPS2022-Liu}. However, GBA is fairly restrictive for applications. For example, it does not hold if UL objective is quadratic w.r.t. LL variable, which is representative in the machine learning context; see e.g., kernel ridge regression, biased regularization and hyper-representation in \cite{grazzi2020iteration}. Recently, when LL objective is strongly convex, a weaker boundedness assumption was made in \cite{ji2022will,dagreou2022framework}. Further, for the case where the total UL objective is strongly convex or convex, by using an induction analysis that all iterates are bounded, \cite{ji2021lower} provided the first-known analyses of the complexity bound without GBA. However, the argument in \cite{ji2021lower} heavily rely on the (strong) convexity assumption on the total UL objective. In comparison, there is no assumption on the total UL objective and also no LLSC assumption in this work.

\section{Proposed Algorithm} \label{sec2}

We consider BLO problem in the following form:
\begin{equation}\label{blo problem}
	\mathop{\min}\limits_{\x, \y} F(\x,\y)\quad \mathrm{s.t.} \quad  \y \in \S(\x):=\arg\min_{\y \in \mathbb{R}^m} f(\x,\y),
\end{equation}
where $\x\in\mathbb{R}^n$, $\y\in\mathbb{R}^m$ are UL and LL variables respectively. Here UL and LL objectives $F,f$ are smooth functions with Lipschitz continuous first and second order derivatives. 

In this section, we propose a single loop Bi-level Averaged Method of Multipliers (sl-BAMM) for BLO  without relying on LL strong convexity (LLSC). 

\subsection{Stationarity Measure}

Most of previous works assume that $f(\x,\cdot)$ is strongly convex for any $\x$. Thus $\S(\x)$ is always a singleton, denoted by $\y^*(\x)$. Then the total UL objective $\Phi(\x):=F(\x,\y^*(\x))$  is continuously differentiable. Following the standard analysis in nonconvex optimization, since $\Phi(\x)$ is nonconvex, algorithms are expected to find a near stationary point of $\Phi(\x)$ as measured by the norm of its gradient (i.e., hypergradient). However, without LLSC, this goal is not suitable since $\mathcal{S}(\x)$ may be a set-valued mapping and $\Phi(\x)$'s variants (in optimistic and pessimistic) may be nonsmooth (even be discontinuous), which makes traditional definitions of stationary point inapplicable. Therefore, one fundamental question arises:

 {\it Question: What is a good stationarity measure for BLO algorithm without LLSC?}

To answer the above question, we start from formulating
BLO to an equality-constrained optimization:
\begin{equation}\label{ecp}
	\min_{\x,\y}\  F(\x,\y)\quad
	\mathrm{s.t.}\quad \nabla_{\y} f(\x,\y)=0.
\end{equation} 
Clearly, Problems~\eqref{blo problem} and~\eqref{ecp} are equivalent when LL objective $f(\x,\y)$ is convex w.r.t. $\y$ for any fixed $\x$. This motivates us to use KKT condition of Problem~\eqref{ecp} as a measure of stationarity of the solution returned by BLO algorithms. 
Let $\mathcal{L}(\x,\y,\v):=F(\x,\y)-\v^T\nabla_{\y} f(\x,\y)$ be the Lagrangian function of Problem~\eqref{ecp}, where $\v$ is the dual multiplier. We say that $(\x^*,\y^*)$ is a KKT point if there exists a dual multiplier $\v^*$ such that $\nabla\mathcal{L}(\x^*,\y^*,\v^*)
=0$. That is, the following KKT condition holds: 
\begin{align}\label{KKT}
	\left(
	\begin{array}{ccc}
		\nabla_{\x} F(\x^*,\y^*)-\nabla_{\x\y}^2 f(\x^*,\y^*)\v^*\\
		\nabla_{\y} F(\x^*,\y^*)-\nabla_{\y\y}^2 f(\x^*,\y^*)\v^*\\
		-\nabla_{\y} f(\x^*,\y^*)
	\end{array}
	\right)=0.
\end{align}
To appropriately characterize the convergence rate of the solution returned by a BLO algorithm, we define KKT residual for $(\x,\y,\v)$ as following
\begin{equation*}
	\mathrm{KKT}(\x, \y,\v):= \|\nabla \mathcal{L}(\x,\y, \v)\|^2.
\end{equation*}
Then, $\mathrm{KKT}(\x, \y,\v) = 0$ if and only if $(\x,\y)$ is a KKT point and $\v$ is a multiplier. 
The next proposition provides a strong connection between KKT condition and the hypergradient using in the existing BLO methods. 

\begin{proposition}\label{prop-equiv}
	Assume that for all $\x$, $f(\x,\cdot)$ is strongly convex. Then $\nabla\Phi(\x)=0$ if and only if $\mathrm{KKT}(\x, \y,\v)=0$ for some $\y,\v\in\mathbb{R}^m$. In fact, $\y=\y^*(\x)$ and 
	\begin{equation*}
		\v=\left[\nabla_{\y\y}^2 f(\x,\y^*(\x))\right]^{-1} \nabla_{\y} F\left(\x,\y^*(\x)\right).
	\end{equation*}
\end{proposition}
Since the nonsingularity of $\nabla_{\y\y}^2 f$ is not necessary for KKT condition, KKT residual is more applicable to general BLOs. 
Further, unlike the computation of hypergradient, which would require knowledge of the unique LL solution to be usable as a stopping criterion, KKT residual is readily available to the algorithm as a stopping criterion. In a novel framework of \cite{dagreou2022framework}, assuming LLSC, the newly involved variable $\v$ is seen as the solution of a linear system: $\left[\nabla_{\y\y}^2 f(\x,\y^*(\x))\right] \v=\nabla_{\y} F\left(\x,\y^*(\x)\right)$. We offer a new perspective on $\v$ that help us to develop a novel single-loop BLO algorithm beyond LLSC.

\subsection{Bi-level Averaged Method of Multipliers}

As shown in Algorithm~\ref{alg:BMA}, by introducing dual multiplier from the viewpoint of KKT condition and adopting the sequencing averaging method in \cite{sabach2017first,liu2020generic,li2020improved}, we propose a single loop Bi-level Averaged Method of Multipliers (sl-BAMM) in which UL variable and the dual multiplier evolve simultaneously with LL variable, following the directions given by 
\begin{align}
	d_{\y}^k=& \nabla_{\y} \psi_{\mu_{k}}(\x_k, \y_k), \label{direction y}\\
	d_{\v}^k=& \nabla_{\y} F(\x_k, \y_k) - \nabla_{\y\y}^2 \psi_{\mu_{k}}(\x_k, \y_k) \v_k, \label{direction v}\\
	d_{\x}^k=& \nabla_{\x} F(\x_k, \y_k) - \nabla_{\x\y}^2 \psi_{\mu_{k}}(\x_k, \y_k) \v_k,
	\label{direction x}
\end{align}
where the aggregation function 
\begin{equation}
	\psi_{\mu}(\x,\y):=\mu  F(\x,\y)+(1-\mu)f(\x,\y),
\end{equation}
averaging UL and LL objectives, is inspired by BDA \cite{liu2020generic}. Here $\mu$ is the averaging (aggregation) parameter that will iteratively goes to zero. 

\begin{algorithm}[h]
	\caption{single loop Bi-level Averaged Method of Multipliers (sl-BAMM)}\label{alg:BMA}
	{\bf Input:} 
	initial points $\x_0,\y_0,\v_0$, stepsizes $\alpha_k,\beta_k,\eta_k$, aggregation parameter $\mu_k$;
	\begin{algorithmic}[1]
		\FOR{$k=0,1,\dots,K-1$} 
		\STATE update $\y_{k+1}=\y_k-\beta_k d_{\y}^k$ according to Eq.~\eqref{direction y};
		\STATE update 
		$\v_{k+1}=\v_k+\eta_k d_{\v}^k$ according to Eq.~\eqref{direction v};
		\STATE update 
		$\x_{k+1}=\x_k-\alpha_k d_{\x}^k$ according to Eq.~\eqref{direction x}.
		\ENDFOR
	\end{algorithmic}
\end{algorithm}
\begin{remark}
	The first and most important motivation of sl-BAMM is to leverage KKT conditions to update all variables, but keeping the update of LL variable by gradient descent of LL objective and its variants. This will enable sl-BAMM to be a provably single loop hence faster algorithm, differently from MPEC approach \cite{kim2020mpec} and the penalty method for BLO in \cite{mehra2021penalty}. 
\end{remark}
\begin{remark}
	Secondly, it is worth noting that all the directions in Eq.~\eqref{direction y}-\eqref{direction x} are linear in $F$ and $f$, differently from the existing BLO algorithms without using an extra variable. This kind of linear structure will motivate more stochastic and global variance reduction algorithms as done in \cite{dagreou2022framework} with LL strong convexity.  
\end{remark}

Note that the dual multiplier in BLO algorithm will result in practical procedures for correcting DARTS \cite{liu2018darts} and other recent algorithms attempt to simplify multi-step iteration with the single gradient,  and removing an unavoidable non-vanishing convergence error in ITD-based approach \cite{ji2022will}.

\section{Theoretical Investigations}\label{Convergence Analysis}

In this section, we provide convergence rates of sl-BAMM with or without LL strong convexity. Note that, unlike existing analyses, {\it we do not require a strong assumption on the boundedness of $\nabla_{\y} F(\x,\y)$}. Hence, our theoretical analysis can be applied to a wider variety of applications in deep learning. The proofs are deferred in Appendix.

\subsection{General Assumptions}

We start by stating some standard assumptions on UL and LL objectives.

\begin{assumption}\label{Assump0}
	The objectives $F$ and $f$ satisfy 
	\begin{enumerate}[(a)]
		\item UL objective $F$ is twice differentiable. The first order derivatives $\nabla_{\x} F(\cdot, \y)$, $\nabla_{\x} F(\x, \cdot)$, $\nabla_{\y} F(\cdot,\y)$, $\nabla_{\y} F(\x,\cdot)$ are Lipschitz continuous with respective Lipschitz constants $L_{F_{\x1}},L_{F_{\x2}}, L_{F_{\y1}}, L_{F_{\y2}}$.
		\item LL objective $f$ is twice differentiable. The derivatives $\nabla_{\y} f$ and $\nabla_{\x\y}^2 f,\nabla_{\y\y}^2 f$ are Lipschitz continuous in $(\x,\y)$ with respective Lipschitz constants $L_{f_{\y1}}, L_{f_{\y2}}$ and $L_{f_{\x\y1}}, L_{f_{\x\y2}}$, $L_{f_{\y\y1}}, L_{f_{\y\y2}}$.
	\end{enumerate}
\end{assumption}

The smoothness assumption above is common in BLO literature, see, e.g., \cite{ghadimi2018approximation,ji2021bilevel,khanduri2021near,ji2021lower,chen2022single,dagreou2022framework,ji2022will}. 
In this paper, we study two classes of BLO: (1) LL objective $f(\x, \cdot)$ is merely convex and UL objective $F(\x, \cdot)$ is strongly convex; (2) LL objective $f(\x, \cdot)$ is strongly convex. 
This can be unified by the following assumption. 
\begin{assumption}\label{asump-strongconvex}
	For any $\x$, the aggregation function $\psi_{\mu}(\x,\cdot)=\mu  F(\x,\cdot)+(1-\mu)f(\x,\cdot)$ is strongly convex for $\mu=0$ or $\mu>0$.
\end{assumption}

\subsection{General Results for LL Merely Convex Case}

To handle BLO with multiple LL minimizers, we take the following standard assumption as in \cite{liu2020generic}.
\begin{assumption}\label{assump_convex0}
	\begin{enumerate}[(a)]
		\item For any fixed $\x$, LL objective $f(\x,\cdot)$ is convex, and UL objective $F(\x,\cdot)$ is $\sigma_F$-strongly convex and has a uniform lower bound denoted by $F_0$.
		\item The derivatives $\nabla_{\x\y}^2 F,\nabla_{\y\y}^2 F$ are Lipschitz continuous in $(\x,\y)$ with respective Lipschitz constants $L_{F_{\x\y1}}, L_{F_{\x\y2}}$, $L_{F_{\y\y1}}, L_{F_{\y\y2}}$.
	\end{enumerate}
\end{assumption}

Under the above assumption, the aggregation function $\psi_{\mu}(\x,\cdot)$ is $\sigma_{\psi_{\mu}}$-strongly convex with $\sigma_{\psi_{\mu}}=\mu \sigma_{F}$. Hence $\psi_{\mu}(\x,\cdot)$ has a unique minimizer, denoted by $\y^*_{\mu}(\x)$. We define the approximate overall UL objective by $\Phi_{\mu_{k}}(\x_k):=F(\x_k,\y^*_{\mu_k}(\x_k))$, and also denote the correct dual multiplier by $\v^*_{\mu}(\x):=\left[\nabla_{\y\y}^2 \psi_{\mu}(\x,\y^*_{\mu}(\x))\right]^{-1}\nabla_{\y} F\left(\x,\y^*_{\mu}(\x)\right)$.
Let $(\x_k,\y_k,\v_k)$ be a sequence, we define 
\begin{equation*}
	\begin{aligned}
		\Pi(\x_k,\y_k,\v_k)=&\| \nabla \Phi_{\mu_{k}}(\x_k) \|^2
		+ \|\x_{k+1}-\x_k\|^2	\\
		&+  \|\y_k-\y_{\mu_k}^*(\x_k)\|^2 
		+ \|\v_k-\v_{\mu_k}^*(\x_k)\|^2.
	\end{aligned}
\end{equation*}

Now we provide a series of convergence rate analysis for sl-BAMM towards KKT stationary points, using three different simple and handy strategies for choosing the step sizes. Specially, the first one does {\bf NOT} require any assumption on the boundedness of $\nabla_{\y} F(\x,\y)$.  

\begin{theorem}\label{thmconvex1}
	Suppose Assumptions \ref{Assump0} and \ref{assump_convex0} hold. Choose $\mu_k= \bar{\mu} (k+1)^{-p}$ with $0<p<1/10$, and stepsizes $\beta_{k}\in \left[ \bar{\beta}(k+1)^{-\tau/2}, 1/ \left(L_{F_{\y2}} + L_{f_{\y2}}\right) \right]$ and 
	\begin{equation}\label{stepsizes1}
		\begin{aligned}
			\eta_k= (k+1)^{-\tau/2} \beta_{k} \mu_{k}^2 , 
			\quad
			\alpha_k= (k+1)^{-3\tau/2} \beta_{k} \mu_{k}^7,
		\end{aligned}
	\end{equation}
with $0<\tau<1/30$. Let $(\x_k,\y_k,\v_k)$ be the sequence generated by sl-BAMM. If $(\x_k, \y_k)$ is bounded, then  
	\begin{align*}
		&\min_{0\leq k\leq K}
		\big\{
		\Pi(\x_k,\y_k,\v_k)
		\big\}
		=O\left(\frac{1}{K^{1-9p-3\tau}}\right).
	\end{align*}
	If $\v_k$ is also bounded, then 
	\begin{equation*}
		\min_{0 \leq k \leq K} \left\{ \mathrm{KKT}(\x_k, \y_k,\v_k) \right\}
		= O\left(\frac{1}{K^{1-9p-3\tau}} + \frac{1}{K^{2p}}\right).
	\end{equation*}
	Furthermore, any limit point of $\left(\x_k,\y_k\right)$ is a KKT point of Problem~\eqref{ecp}, and the convergence rate of $\x$ is given by $\|\x_{k+1}-\x_{k}\|=O(k^{-\frac{5p}{2}-\frac{\tau}{4}})$. 
\end{theorem}

Next we establish tighter convergence rates when a weaker gradient boundedness assumption holds.
\begin{theorem}\label{thmconvex2}
	Suppose Assumptions \ref{Assump0} and \ref{assump_convex0} hold. Under the same setting of Theorem \ref{thmconvex1} except $0<p<1/6$, $0<\tau<1/18$ and replacing stepsizes \eqref{stepsizes1} by 
	\begin{equation}\label{stepsizes2}
		\begin{aligned}
			\eta_k= (k+1)^{-\tau/2} \beta_{k} \mu_{k} ,
			\quad 
			\alpha_k= (k+1)^{-3\tau/2} \beta_{k} \mu_{k}^5 .
		\end{aligned}
	\end{equation}
	Let $(\x_k,\y_k,\v_k)$ be the sequence generated by sl-BAMM. If $\nabla_{\y} F(\x_k, \y^*_{\mu_{k}}(\x_{k}))$ is bounded, then  
	\begin{align*}
		&\min_{0\leq k\leq K}
		\big\{
		\Pi(\x_k,\y_k,\v_k)
		\big\}
		=O\left(\frac{1}{K^{1-5p-3\tau}}\right).
	\end{align*}
	If $\v_k$ is also bounded, then 
	\begin{equation*}
		\min_{0 \leq k \leq K} \left\{ \mathrm{KKT}(\x_k, \y_k,\v_k) \right\}= O\left(\frac{1}{K^{1-5p-3\tau}} + \frac{1}{K^{2p}}\right).
	\end{equation*}
	Furthermore, any limit point of $\left(\x_k,\y_k\right)$ is a KKT point of Problem~\eqref{ecp}, and the convergence rate  of $\x$ is given by $\|\x_{k+1}-\x_{k}\|=O(k^{-\frac{5p}{2}-\frac{\tau}{4}})$. 
\end{theorem}
Note that $\nabla_{\y} F(\x_k, \y^*_{\mu_{k}}(\x_{k}))$ is bounded if $\x_k$ is bounded, and for any $B>0$ there are positive constants $C$ and $\mu_0$ such that $\| \nabla_{\y} F(\x,\y^*_{\mu}(\x)) \|\leq C$ for all $\| \x \| \leq B$ and $0<\mu\leq\mu_0$. It is satisfied for a broad collection of applications, including BLOs with multiple LL minimizers in our experiment. Specially, assuming LLSC, $\nabla_{\y} F(\x_k, \y^*_{\mu_{k}}(\x_{k}))$ is bounded if there exists $C_F>0$ such that $\| \nabla_{\y} F(\x,\y^*(\x)) \|\leq C_{F}$ for all $\x$, a weaker boundedness assumption in \cite{ji2022will,dagreou2022framework}.

Finally, if $\| \v^*_{\mu_k}(\x_k) \|$  is bounded, the above results can be improved further.

\begin{theorem}\label{thmconvex3}
	Suppose Assumptions \ref{Assump0} and \ref{assump_convex0} hold. Under the same setting of Theorem \ref{thmconvex1} except $0<p<1/4$, $0<\tau<1/12$ and replacing stepsizes \eqref{stepsizes1} by 
	\begin{equation}\label{stepsizes3}
		\begin{aligned}
			\eta_k= (k+1)^{-\tau/2} \beta_{k} , 
			\quad
			\alpha_k= (k+1)^{-3\tau/2} \beta_{k} \mu_{k}^3 .
		\end{aligned}
	\end{equation}
	Let $(\x_k,\y_k,\v_k)$ be the sequence generated by sl-BAMM. If $\|\v^*_{\mu_k}(\x_k)\| $ is bounded, then we have  
	\begin{align*}
		&\min_{0\leq k\leq K}
		\big\{
		\Pi(\x_k,\y_k,\v_k)
		\big\}
		=O\left(\frac{1}{K^{1-3p-3\tau}}\right).
	\end{align*}
	If $\v_k$ is also bounded, then 
	\begin{equation*}
		\min_{0 \leq k \leq K} \left\{ \mathrm{KKT}(\x_k, \y_k,\v_k) \right\}= O\left(\frac{1}{K^{1-3p-3\tau}} + \frac{1}{K^{2p}}\right).
	\end{equation*}
	Furthermore, any limit point of $\left(\x_k,\y_k\right)$ is a KKT point of Problem~\eqref{ecp}, and the convergence rate of $\x$ is given by $\|\x_{k+1}-\x_{k}\|=O(k^{-\frac{3p}{2}-\frac{\tau}{4}})$. 
\end{theorem}

\subsection{Special Discussion for LL Strongly Convex Case}

For the special case where LL objective $f(\x,\cdot)$ is strongly convex, Assumption \ref{asump-strongconvex} holds where $\mu=0$. Hence Assumption \ref{assump_convex0} can be removed. 

\begin{assumption}\label{assump_convex20}
	For fixed $\x$, $f(\x,\cdot)$ is $\sigma_{f}$-strongly convex.
\end{assumption}

Under Assumption \ref{assump_convex20}, LL solution $\y^*(\x)$ is well defined and differentiable. Hence, one can use the gradient of $\Phi(\x)=F(\x,\y^*(\x))$ to characterize the optimality of the generated sequence by the algorithms. Actually, as shown in Proposition \ref{prop-equiv}, it is equivalent to KKT residual. 

Relinquishing any assumption on the boundedness of $\nabla_{\y} F(\x,\y)$, our next result provides a new convergence rate guarantee for LL strongly convex case.
\begin{theorem}\label{thm_Strong_convex}
	Suppose Assumption \ref{assump_convex20} holds. Consider sl-BAMM with $\mu_k=0$. We choose $\tau>0$ and stepsizes
	\begin{equation}
		\begin{aligned}
			&\bar{\beta}(k+1)^{-\tau/2}\leq \beta_k
			\leq
			\frac{1}{ L_{F_{\y2}} + L_{f_{\y2}}} , \\
			&\eta_k= \bar{\eta} (k+1)^{-\tau/2} \beta_{k} , 
			\ 
			\alpha_k= \bar{\alpha} (k+1)^{-\tau} \beta_{k}  ,
		\end{aligned} 
	\end{equation}
	where $\bar{\beta}, \bar{\eta}, \bar{\alpha}$ are positive constants satisfying certain conditions. Let $(\x_k,\y_k,\v_k)$ be the sequence generated by sl-BAMM. If $(\x_k, \y_k)$ is bounded, then $\x_k$ satisfies
	\begin{equation}
		\min_{0 \le k \le K} \{\|\nabla \Phi(\x_k)\|^2 \}= O\left(\frac{1}{K^{1-\frac{3\tau}{2}}}\right).
	\end{equation}
\end{theorem}
\begin{remark}
	If GBA holds, that is, there is a constant $C$ such that $\| \nabla_{\y} F(\x,\y) \|\leq C$ for all $\x, \y$, then $\tau$ in Theorem \ref{thm_Strong_convex} can be taken to be $0$ and then $\min_{0 \le k \le K} \{\|\nabla \Phi(\x_k)\|^2 \}= O\left(1/K\right)$, which has been proved recently in \cite{li2022fully} (in the stochastic setting) and in  \cite{ji2022will} (in the deterministic setting). Since $\tau$ can go to $0$, Theorem \ref{thm_Strong_convex} shows that sl-BAMM with or without GBA have  almost the same convergence rate. It is worth noting that GBA does not hold if UL objective is quadratic in LL variable, which is representative in the machine learning context; referring to  kernel ridge regression, biased regularization and hyper-representation in \cite{grazzi2020iteration}.
\end{remark}

To analyze the convergence of sl-BAMM without LLSC and also relinquishing any assumption on the boundedness of $\nabla_{\y} F(\x,\y)$, according to different step size strategies, we will utilize different intrinsic Lyapunov functions with flexible coefficients given by 
\begin{equation*}
	\begin{aligned}
		V_k
		:=&a_k \left[ F(\x_k,\y^*_{\mu_k}(\x_k)) - F_0\right]\\
		&+b_k \|\y_k-\y^*_{\mu_k}(\x_k)\|^2
		+ c_k \|\v_k-\v^*_{\mu_k}(\x_k)\|^2.
	\end{aligned}
\end{equation*}
Here $\{ a_k, b_k, c_k \}_{k=1}^\infty$ is a sequence of positive constants depending on the step size strategy. In particular, when $f(\x, \cdot)$ is strongly convex, we take $\mu_{k}=0$ for all $k$ and then $F(\x, \y^*(\x))$ is the total UL objective. 
Generally, the first term $\Phi_{\mu}(\x)=F(\x,\y^*_{\mu}(\x))$ quantifies the (approximate) overall UL objective functions, the second term $\|\y_k-\y^*_{\mu_k}(\x_k)\|^2$ characterizes LL solution errors, and the third term $\|\v_k-\v^*_{\mu_k}(\x_k)\|^2$ delineates the multiplier errors. 

We will first analyze the descent of the approximate overall UL objective in the next lemma.
\begin{lemma}\label{funda_lemma}
	Suppose Assumptions \ref{Assump0}, and either Assumption \ref{assump_convex0} or Assumption \ref{assump_convex20} holds. Let $\mu_{k+1}\leq\mu_{k}\leq\frac{1}{2}$ for all $k$. Then the sequence of $\x_k,\y_k,\v_k$ generated by sl-BAMM satisfies
	\begin{align*}
		&\Phi_{\mu_{k+1}}(\x_{k+1})-\Phi_{\mu_{k}}(\x_k)  \nonumber\\
		\leq & -\frac{\alpha_k}{2}\|\nabla \Phi_{\mu_{k}}(\x_k)\|^2\\
		&-\frac{1}{2}\left(\frac{1}{\alpha_k}
		-\frac{C_{\Phi1}\|\v^*_{\mu_k}(\x_k)\|+C_{\Phi2}}{\sigma_{\psi_{\mu_k}}^2}\right)\|\x_{k+1}-\x_k\|^2 
		\nonumber\\  
		&+\alpha_k \left( L_{\psi_{\x\y2}} \left\| \v^*_{\mu_{k}}(\x_k)\right\| + L_{F_{\x2}} \right)^2 
		\|\y_{k}-\y^*_{\mu_k}(\x_k)\|^2\\
		&+\alpha_k L_{\psi_{\y1}}^2\|\v_{k}-\v^*_{\mu_{k}}(\x_k)\|^2\\
		&+\frac{2\big\|\nabla_{\y} F(\x_{k+1},\y^*_{\mu_{k+1}}(\x_{k+1}))\big\|^2}{\sigma_{F}} \left(\frac{\mu_k-\mu_{k+1}}{\mu_k}\right),
	\end{align*}
	where $C_{\Phi1}, C_{\Phi2}$ are some explicit positive constants given in Lemma~\ref{lem12}.
\end{lemma}
\begin{remark}
	Compared to the previous results using both LLSC and GBA, the coefficients of $\| \x_{k+1} - \x_{k} \|^2$ and $\|\y_{k}-\y^*_{\mu_k}(\x_k)\|^2$ depend on $\| \v^*_{\mu_{k}}(\x_k) \| $ since Lipschitz continuity of the hypergradient and its surrogate could not be guaranteed without any assumption on the boundedness  of $\nabla_{\y} F(\x,\y)$ and its variant. To address this issue, we characterize a weaker smoothness of $\Phi_{\mu}(\cdot)$ in Lemma \ref{lem12}. Note also that there is an extra term in the inequality of Lemma \ref{funda_lemma}  involving the descent of the averaging parameter. Here we adopt the convention that $0/0=0$. 
\end{remark}

For the convenience of the reader, we give a unified proof sketch of Theorems in Section \ref{Convergence Analysis} in Appendix \ref{proofsketch}.

\section{Experimental Results}\label{experiments}

In this section, we conduct experiments to study basic properties of sl-BAMM such as correctness, scalability and practicality. We also compare the performances of sl-BAMM with competitive methods on different tasks. Specially, we first conduct experiments on toy problems to verify the convergence results established in Section~\ref{Convergence Analysis}. Then we test the performance of sl-BAMM on two real-world BLO applications, including data hyper-cleaning and few-shot classification, compared with representative BLO methods. Note that the detailed hyper-parameter settings for numerical examples, data hyper-cleaning and few-shot classification tasks are provided in Appendix E.1, E.2 and E.3, respectively.  The code is available under~\url{https://github.com/vis-opt-group/sl-BAMM/}. %

\subsection{Numerical Verification}\label{4.1}
\begin{figure*}[htt]
	\centering
	\setlength{\tabcolsep}{0.3mm}{
		\begin{tabular}{c@{\extracolsep{0.1em}}c@{\extracolsep{0.8em}}c@{\extracolsep{0.1em}}c@{\extracolsep{0.8em}}c@{\extracolsep{0.1em}}c@{\extracolsep{0.8em}}c@{\extracolsep{0.1em}}c@{\extracolsep{0.1em}}}
			\includegraphics[height=0.111\textheight,width=0.16\linewidth]{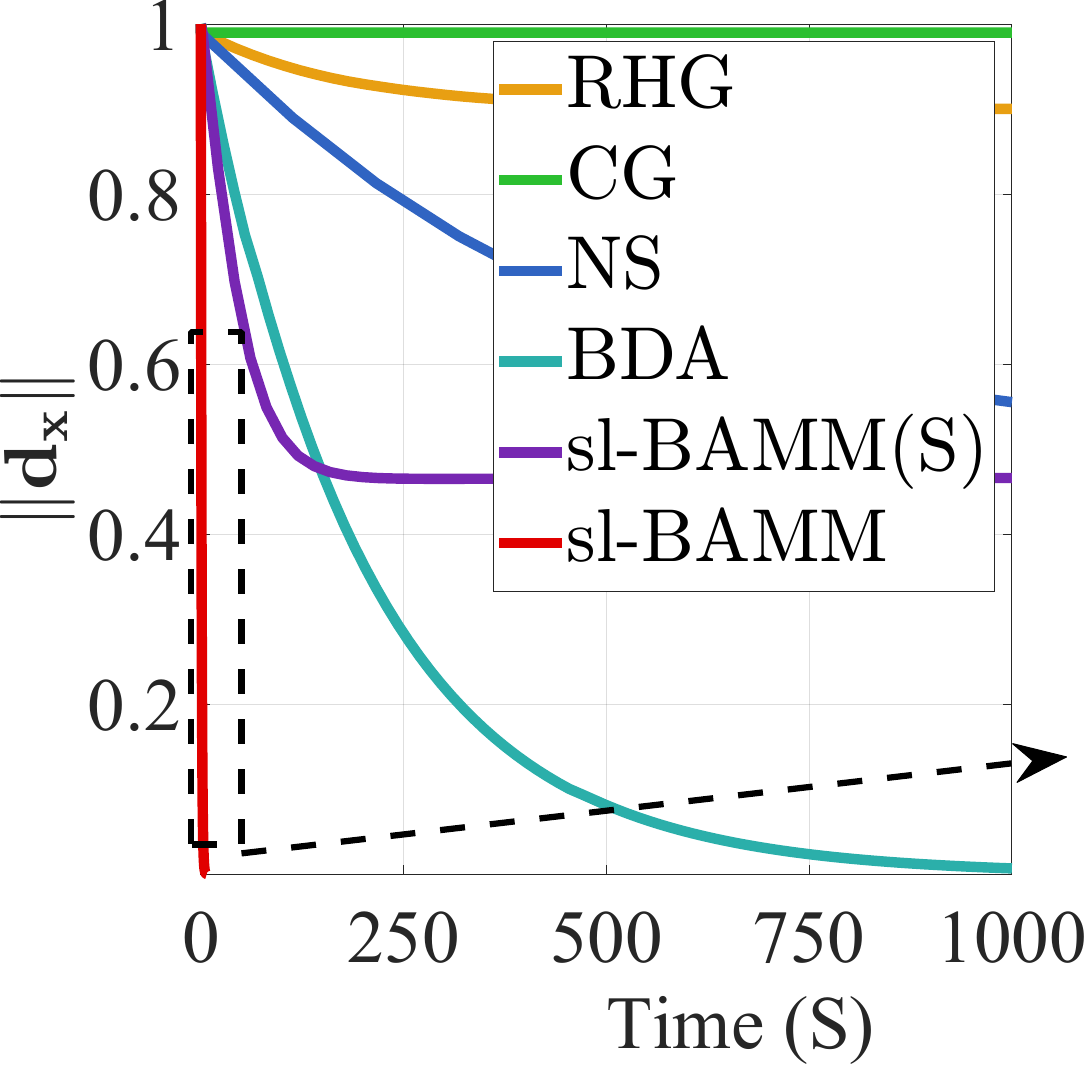}&
			\includegraphics[height=0.111\textheight,width=0.05\linewidth,trim=0 -55 0 0,clip]{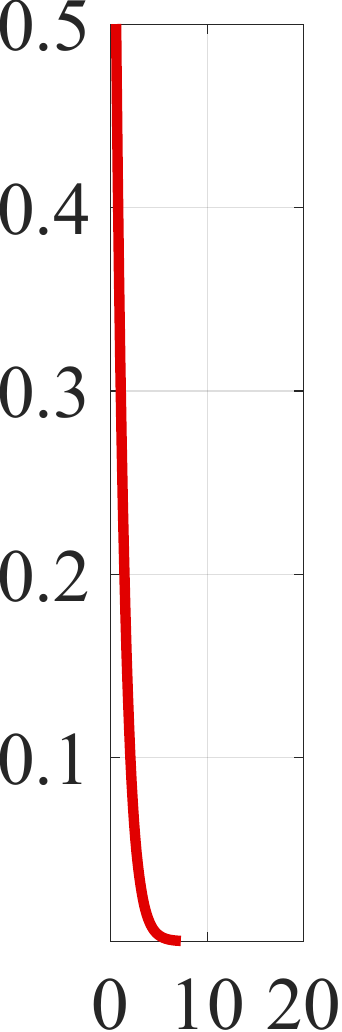}&
			\includegraphics[height=0.11\textheight,width=0.16\linewidth]{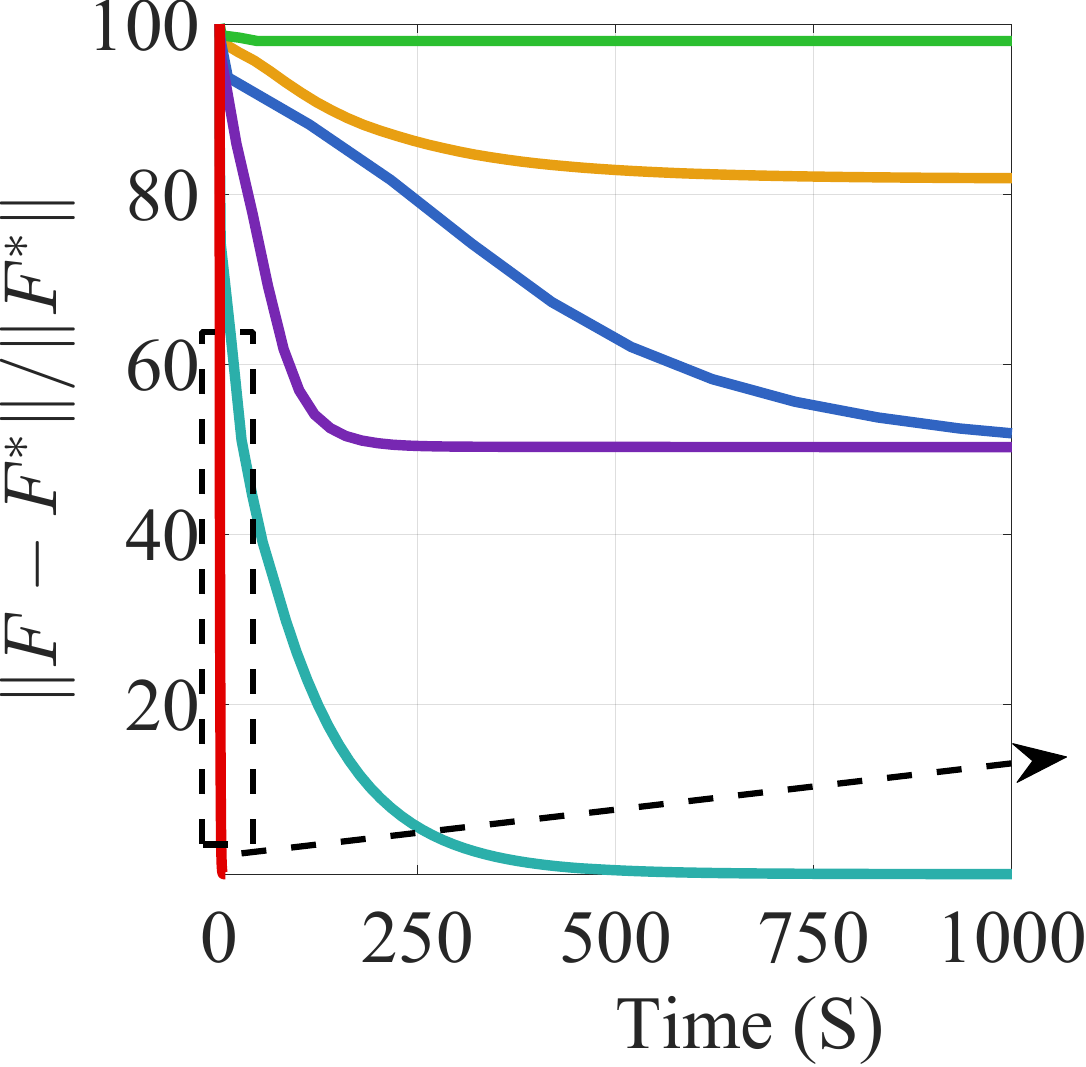}&
			\includegraphics[height=0.11\textheight,width=0.05\linewidth,trim=0 -55 0 0,clip]{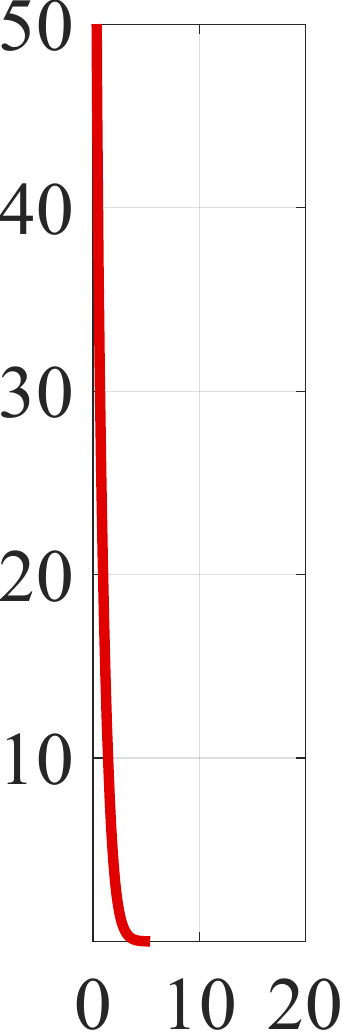}&
			\includegraphics[height=0.109\textheight,width=0.16\linewidth]{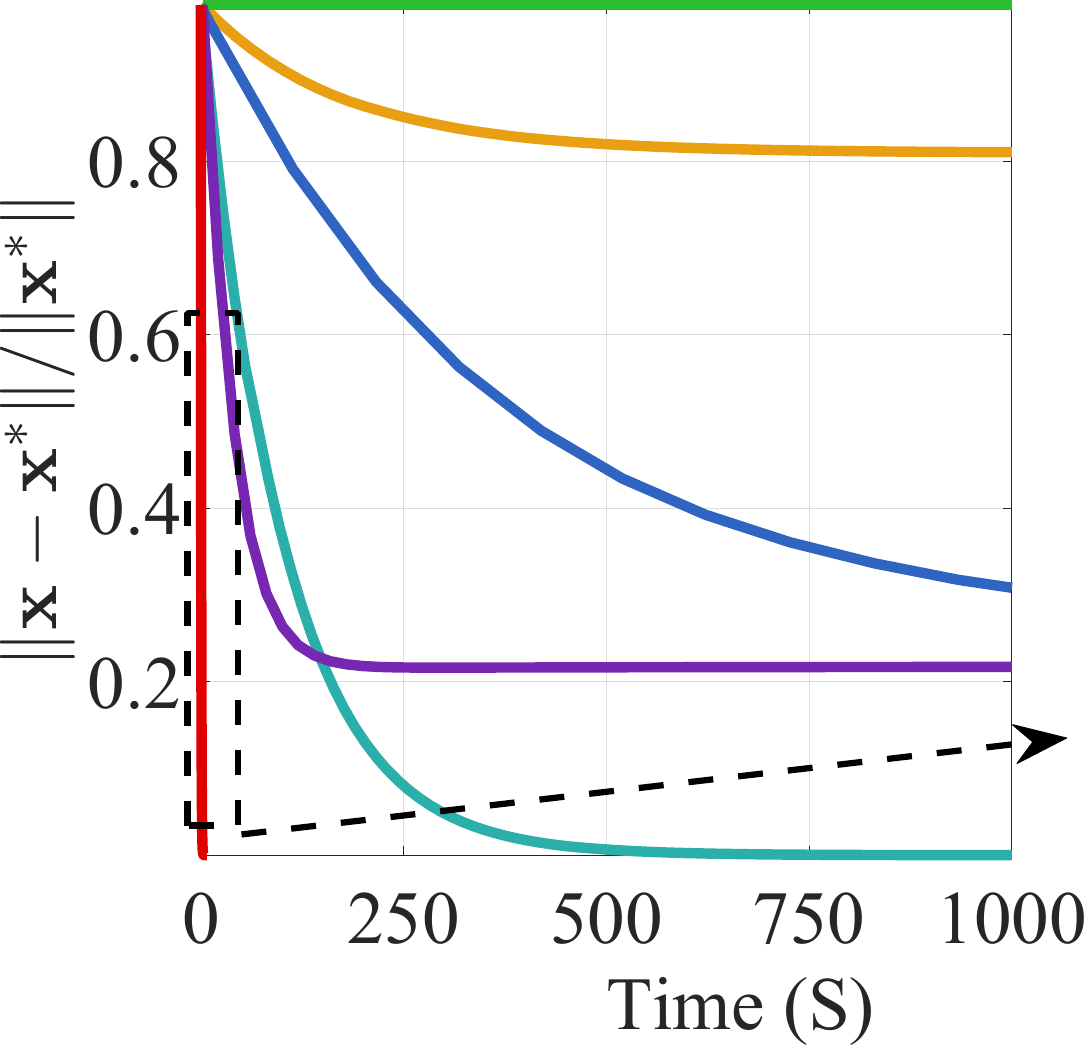}&
			\includegraphics[height=0.11\textheight,width=0.05\linewidth,trim=0 -55 0 0,clip]{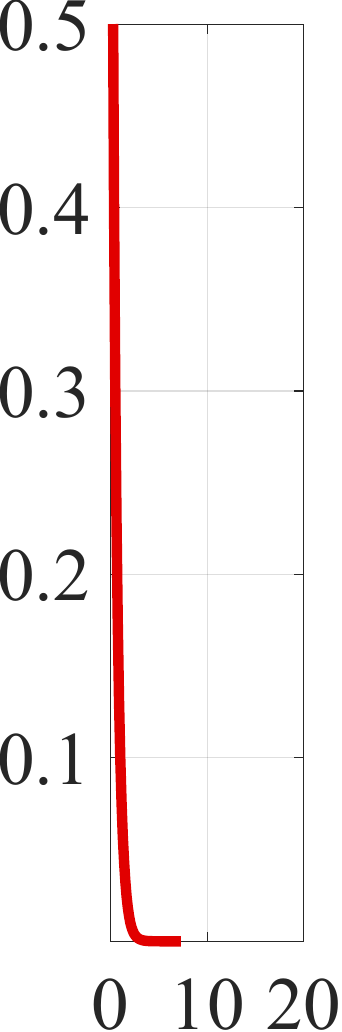}&
			\includegraphics[height=0.11\textheight,width=0.16\linewidth]{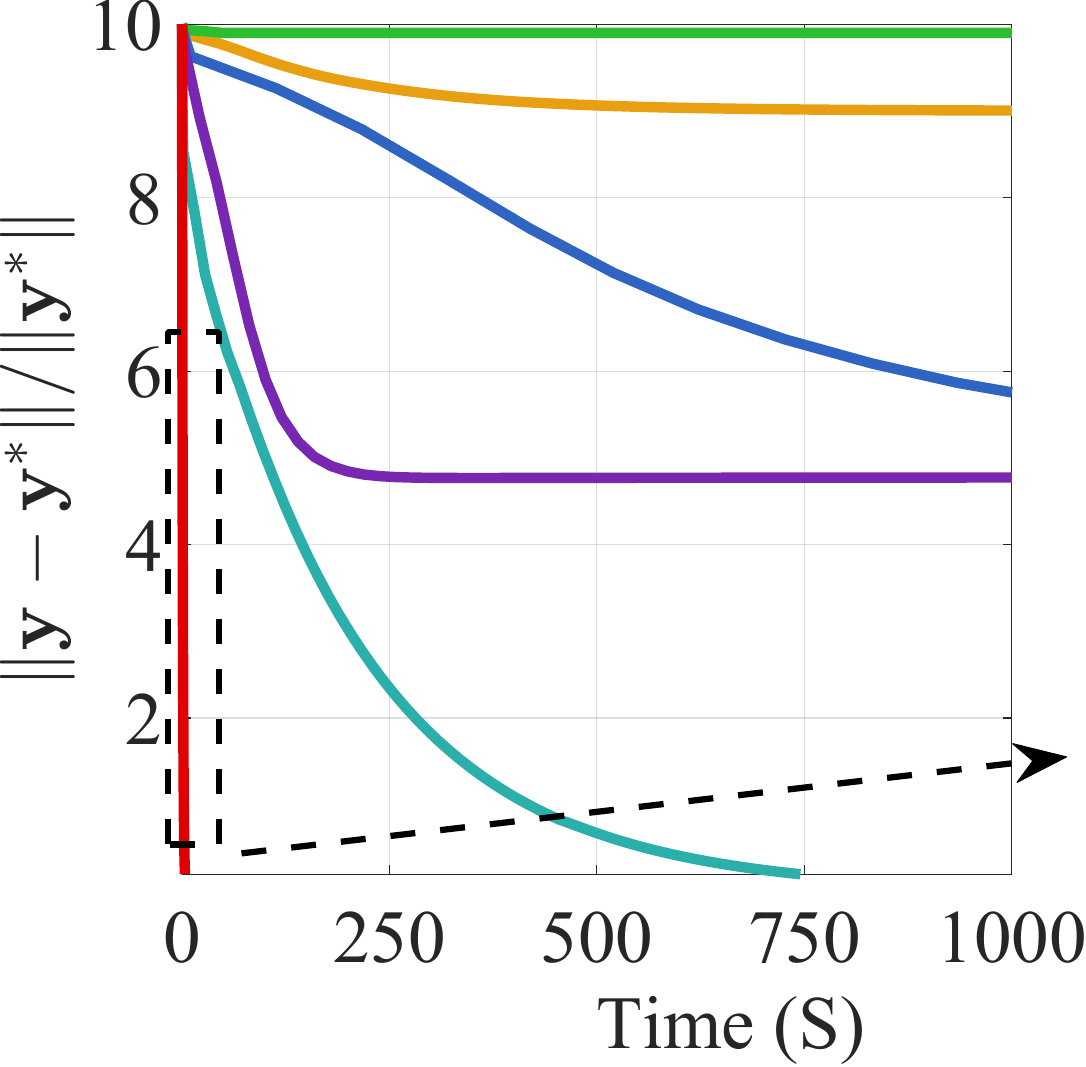}&
			\includegraphics[height=0.11\textheight,width=0.045\linewidth,trim=0 -55 0 0,clip]{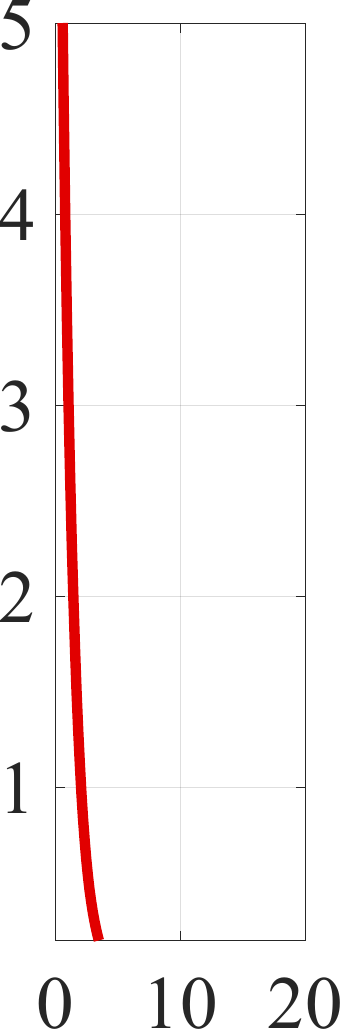}
			\\
		\end{tabular}
	}	 
	\caption{Illustrating the convergence curves of RHG, CG, NS, BDA and sl-BAMM based on different indicators, including $\mathbf{d}_{\x}$ ($\mathbf{d}_{\x}$ denotes the direction of descent of $\x$), $F$, $\x$, and $\y$. We also zoom in on the right side of each subfigure for better view of the  curve at the very beginning.  Note that sl-BAMM (S) represents the version of sl-BAMM with $\mu_k=0$.} 
	\label{fig:LLC_1}
\end{figure*}

\textbf{LL Merely Convex Case}. In the following, we first verify the convergence results under LL merely convex cases based on the toy example problem used in BDA~\cite{liu2020generic}:
\begin{eqnarray}\label{eq:LLC}
	\begin{aligned}
		&\min _{\x \in \mathbb{R}^n} \frac{1}{2}\|\x-\y_2\|^2+\frac{1}{2}\|\y_1-\mathbf{e}\|^2 \quad \\
		&\mathrm { s.t. }\quad\y=(\y_1,\y_2)\in\arg \min _{(\y_1,\y_2) \in \mathbb{R}^{2n}} \frac{1}{2}\|\y_1\|^2-\x^\top \y_1,\\
	\end{aligned}
\end{eqnarray} 
where $\mathbf{e}$ denotes the vector of which the elements are all equal to $1$. By simple calculation, we can calculate the optimal solution as $(\mathbf{e},\mathbf{e},\mathbf{e})$.  Note that we use $\mathbf{d}_{\x}$ to denote the update direction at $\x$ calculated by different tested BLO methods, including sl-BAMM, RHG, Conjugate Gradient (CG)~\cite{pedregosa2016hyperparameter}, Neumann series (NS)~\cite{rajeswaran2019meta} and BDA. 
In the following, though the updates of $\x$, $\y$ and $\v$ of  sl-BAMM are implemented sequentially, to show the potential efficiency of the parallelism, we use the maximum time cost for updating $\x$, $\y$ and $\v$ respectively as the runtime of sl-BAMM at each iteration.   
\begin{figure*}[h!]
	\setlength{\tabcolsep}{0.3mm}{
		\begin{tabular}{c@{\extracolsep{0.03em}}c@{\extracolsep{0.03em}}c@{\extracolsep{0.03em}}c@{\extracolsep{0.03em}}c@{\extracolsep{0.03em}}c@{\extracolsep{1em}}}  			
			
			\includegraphics[height=0.11\textheight,width=0.163\linewidth]{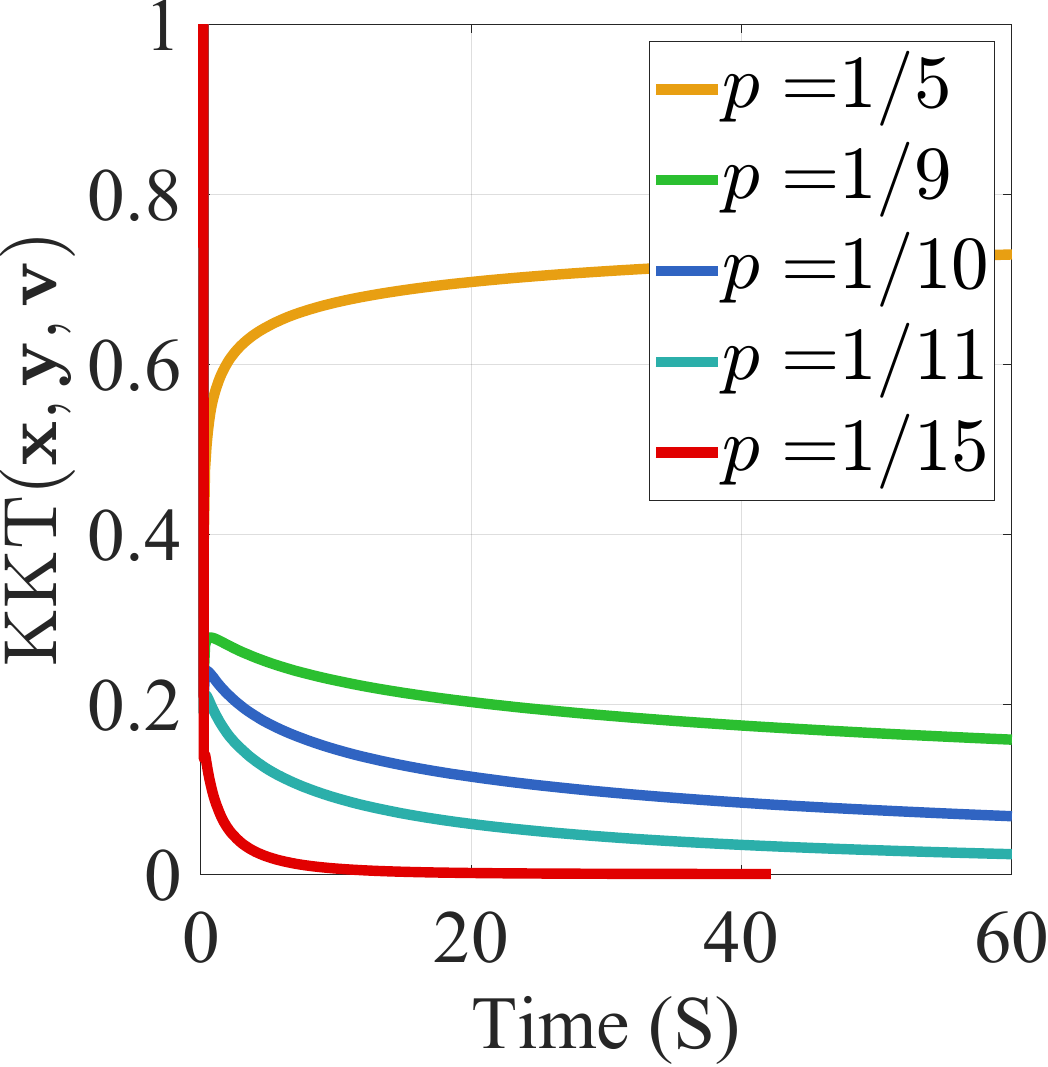}&
			\includegraphics[height=0.109\textheight,width=0.163\linewidth]{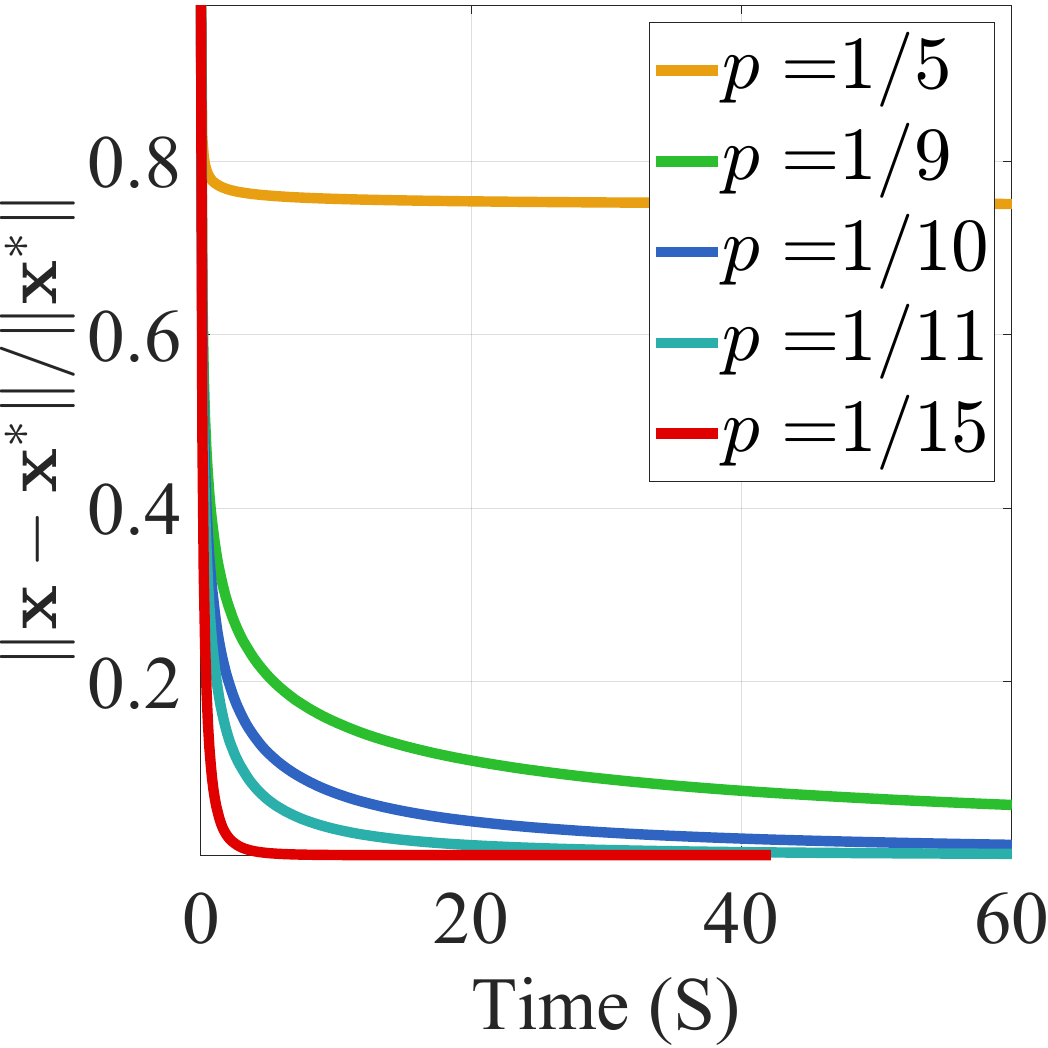}&
			\includegraphics[height=0.11\textheight,width=0.163\linewidth]{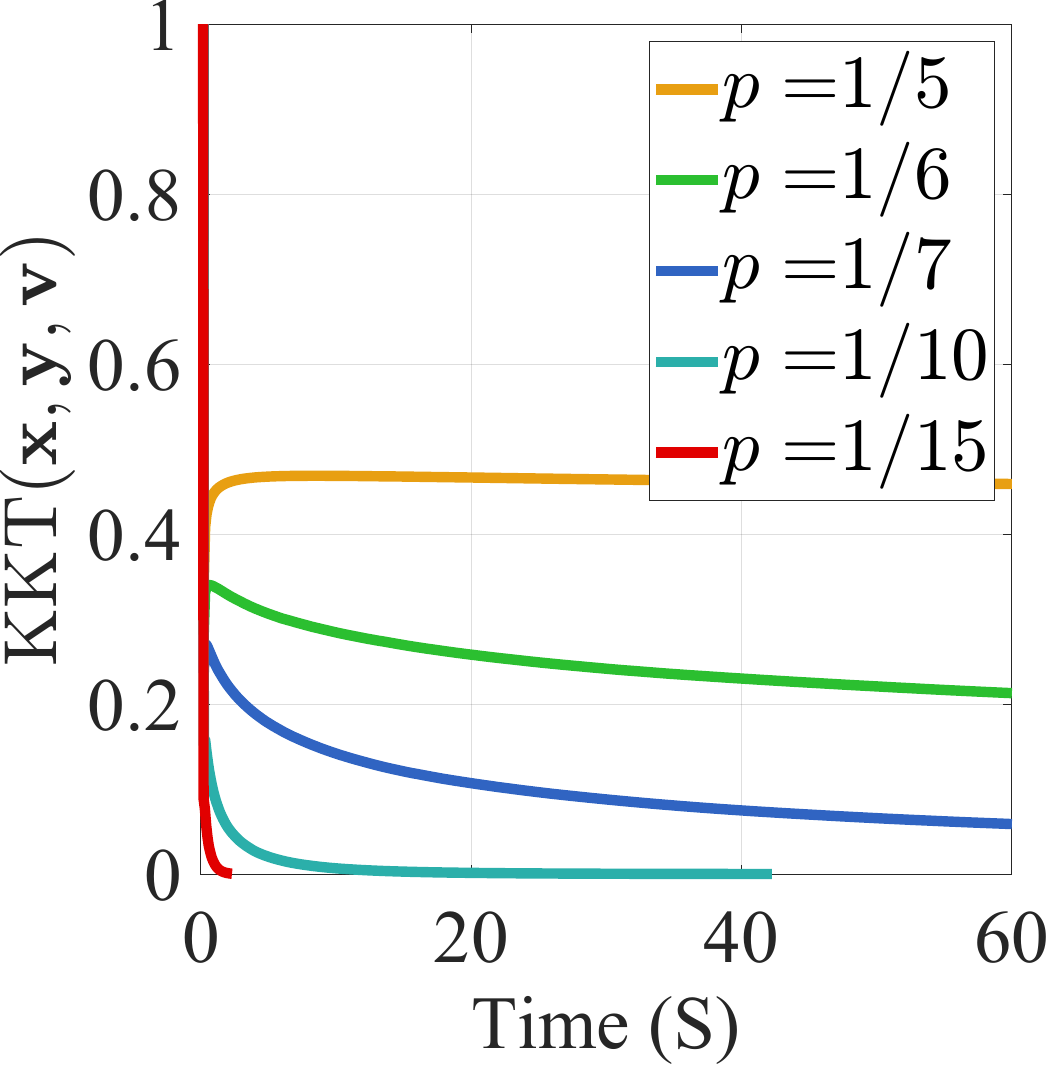}&
			\includegraphics[height=0.109\textheight,width=0.163\linewidth]{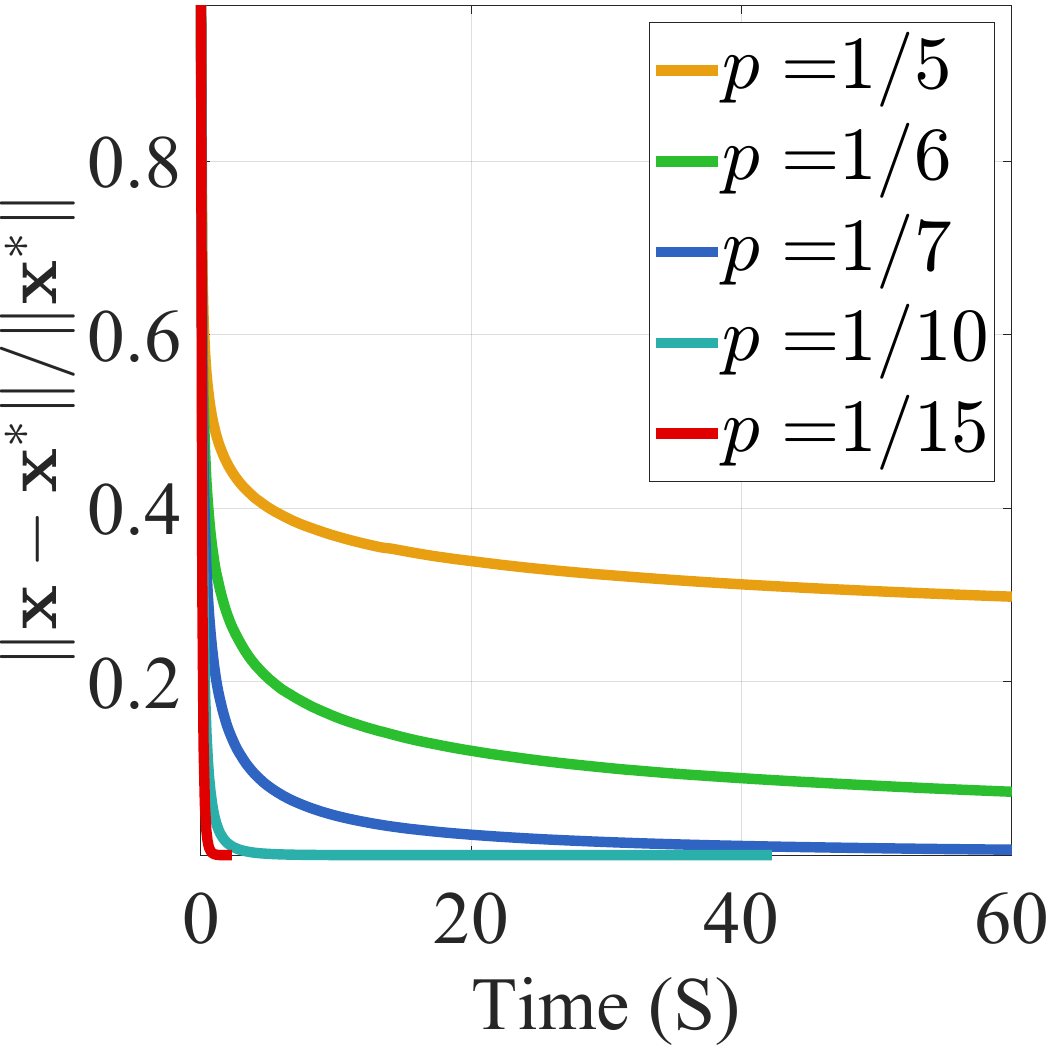}&
			\includegraphics[height=0.11\textheight,width=0.163\linewidth]{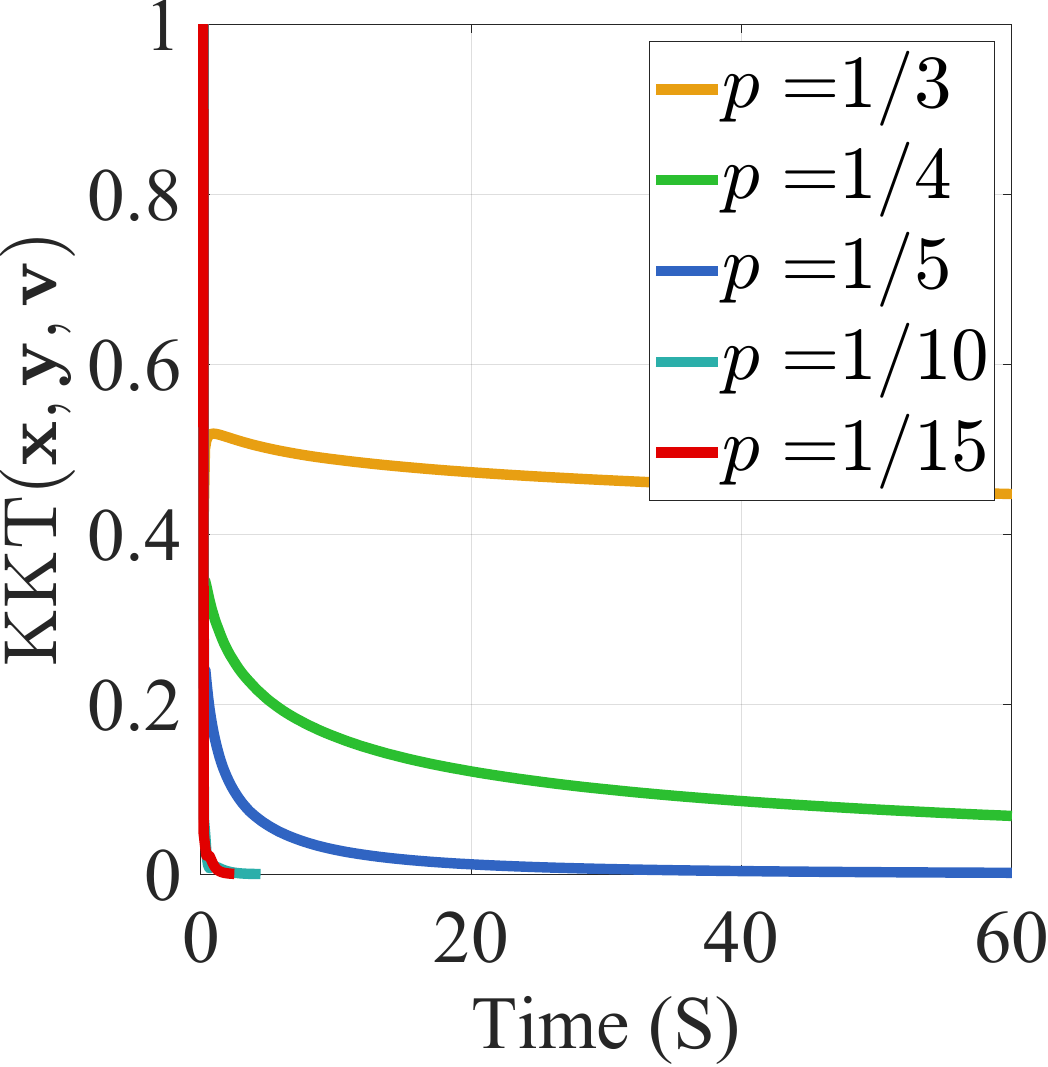}&
			\includegraphics[height=0.109\textheight,width=0.163\linewidth]{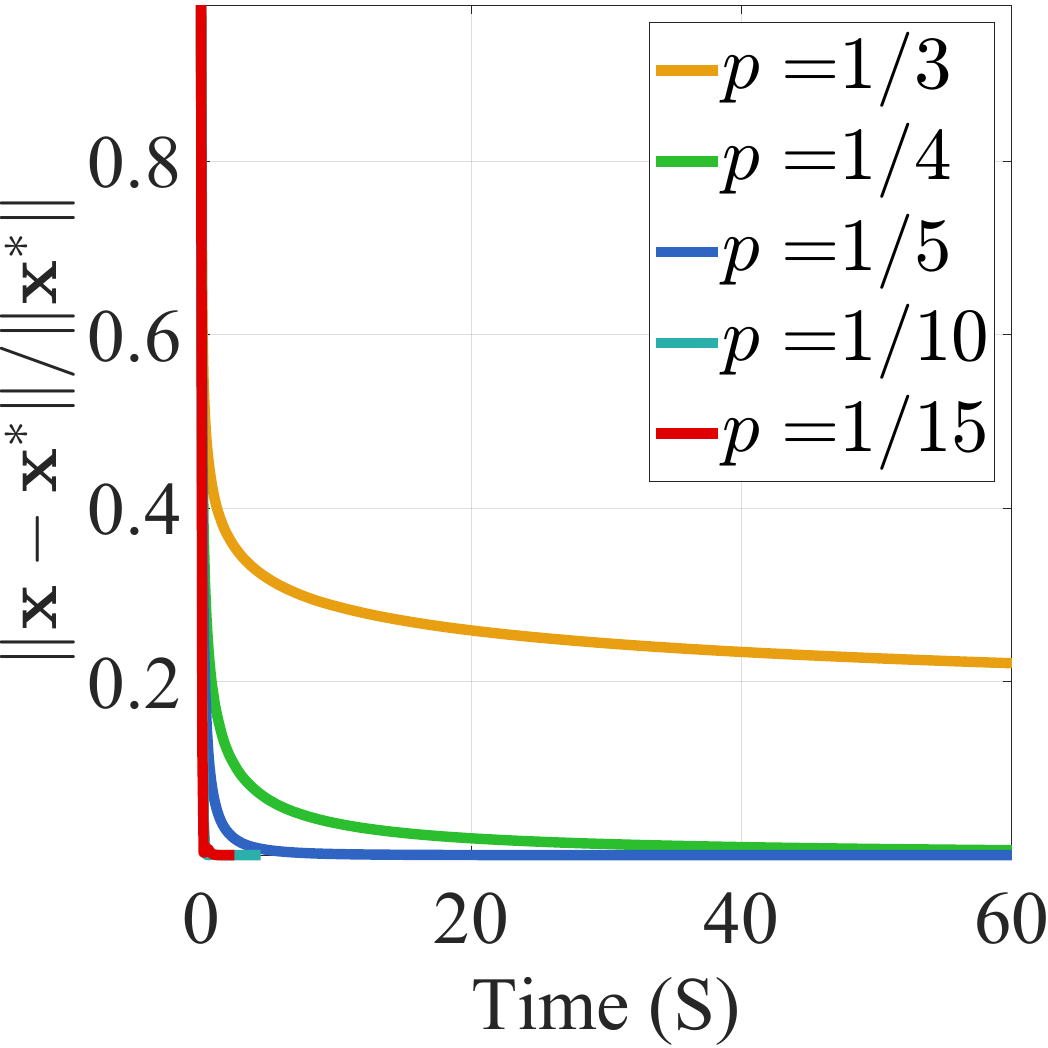}
			\\
			\multicolumn{2}{c}{\footnotesize (a) Strategy S1}&\multicolumn{2}{c}{\footnotesize (b) Strategy S2}&\multicolumn{2}{c}{\footnotesize (c) Strategy S3}\\
		\end{tabular}
	}	 	
	\vspace{-0.2cm}
	\caption{The convergence curves of sl-BAMM with three different strategies to update $\mu_{k}$, $\alpha_{k}$, $\beta_{k}$ and $\eta_k$. In each subfigure, we illustrate the convergence behavior of $\x$ and $\mathrm{KKT}(\x, \y,\v)$ by choosing different $p$ to update $\mu_{k}$.}
	
	\label{fig:LLC_2}
\end{figure*}

Since other methods don't involve the new introduced variable, for fairness and uniformity, we first compare all of the methods based on several widely accepted indicators in Fig.~\ref{fig:LLC_1}. As it is shown, since only BDA and sl-BAMM are capable to handle BLOs without LLSC assumption, they could converge to the truly global optimal solution, while other methods fail to  converge to the correct solution. Moreover, sl-BAMM shows significant advantage over BDA in terms of the convergence speed due to a single-loop update manner. Besides, it can be observed that if $\mu_{k}$ is set as $0$ in sl-BAMM  (denoted as sl-BAMM (S)), it no longer achieves the convergence property, which also validates that the aggregation step in sl-BAMM is indispensable.

Here, we use S1, S2 and S3 to denote three different parameter chosen strategies for sl-BAMM proposed in Theorem~\ref{thmconvex1},~\ref{thmconvex2} and~\ref{thmconvex3}, respectively. In Fig.~\ref{fig:LLC_2}, we test these three strategies with fixed $\tau$ and different values of $p$. It can be seen that when the value of $p$ violates the bound required in the convergence result in Section 3.2 much, sl-BAMM fails to converge in terms of both UL variable and KKT residual. Besides, in this toy example, which satisfies all the assumptions required by these three strategies, strategy S3 always has the fastest convergence speed.

\begin{figure}[htt]
	\vspace{0.4cm}
	\begin{center}
		\renewcommand\arraystretch{0.1}
		\begin{tabular}{c@{\extracolsep{1.2em}}c@{\extracolsep{1.2em}}c@{\extracolsep{0.8em}}}
			&\multicolumn{2}{r}{\includegraphics[width=0.6\linewidth,trim=70 0 0 350,clip]{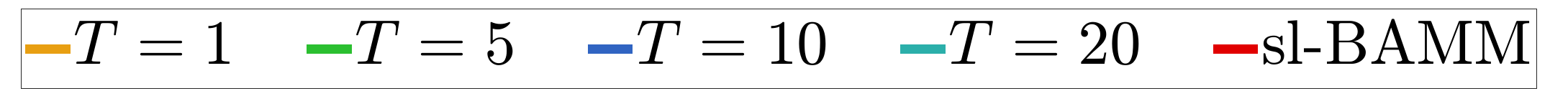}}\\
			\multirow{1}{*}{\rotatebox{90}{RHG~~~~~~~~}}&&\\
			&\includegraphics[height=0.1\textheight,width=0.333\linewidth]{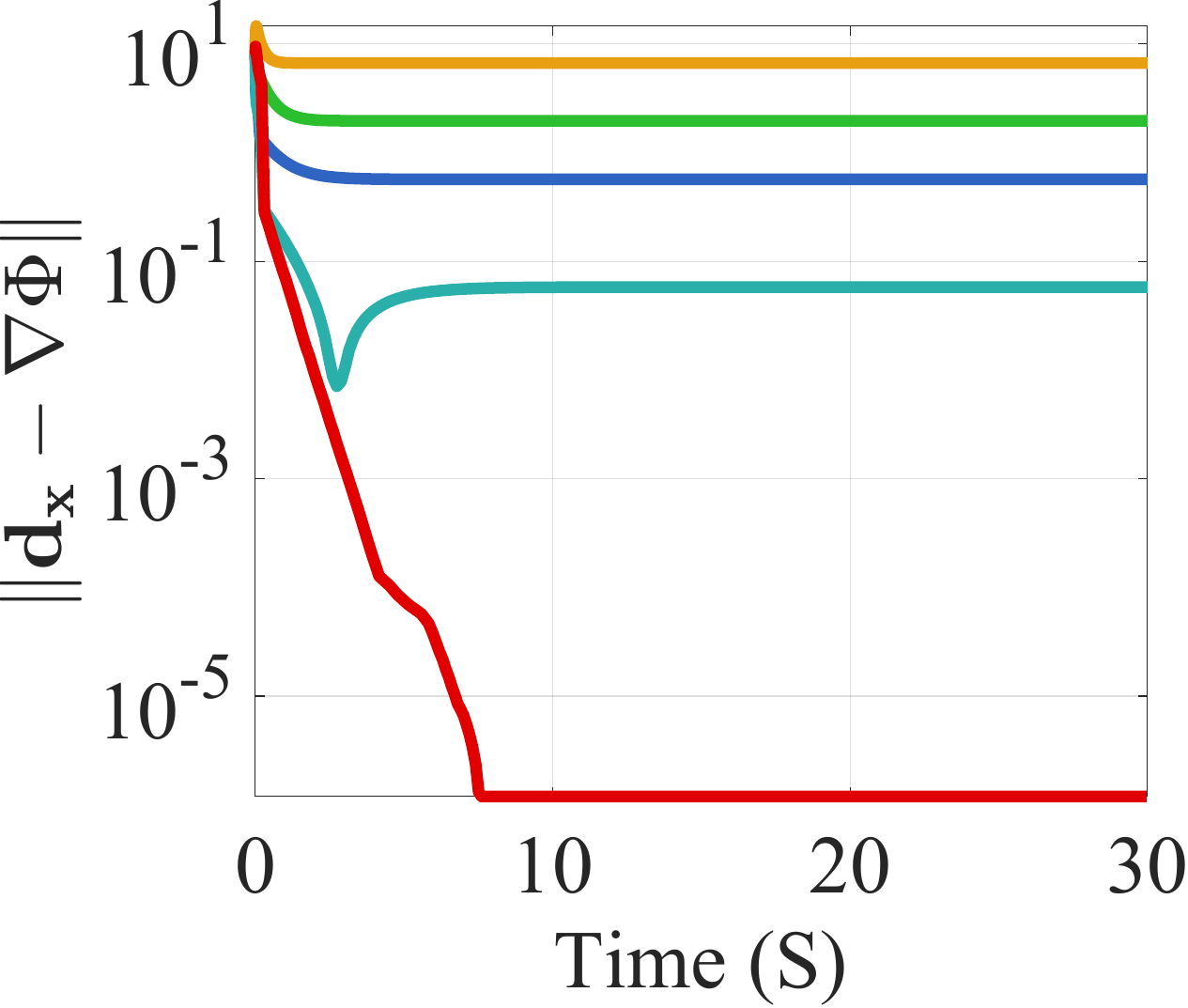}&
			\includegraphics[height=0.1\textheight,width=0.333\linewidth]{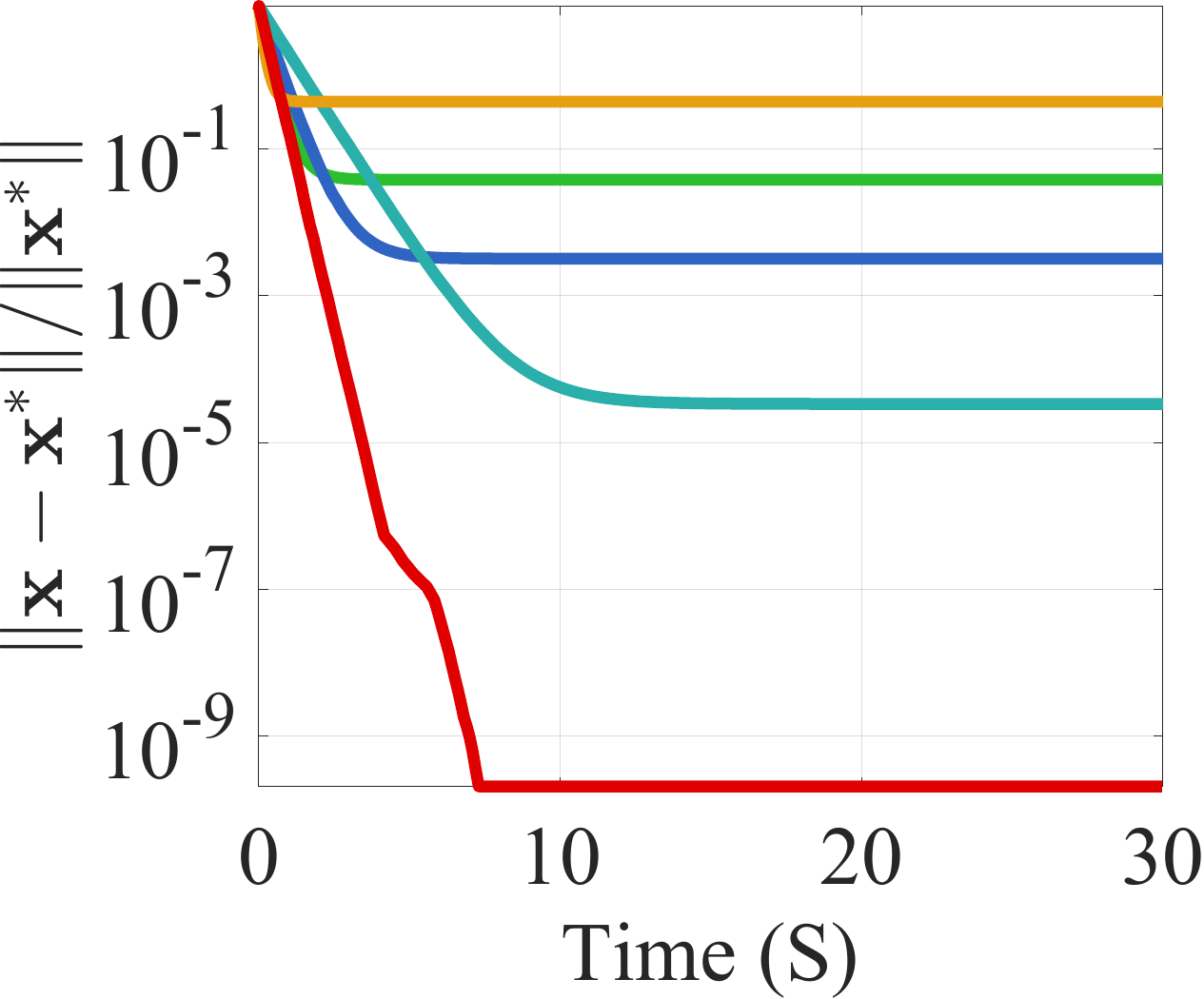}\\
			&\multicolumn{2}{r}{\includegraphics[width=0.6\linewidth,trim=10 0 0 350,clip]{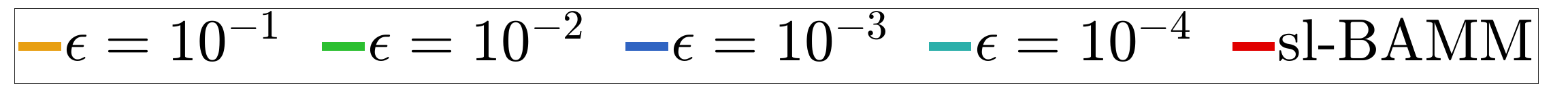}}\\
			\multirow{1}{*}{\rotatebox{90}{CG~~~~~~~~~}}&&\\
			&\includegraphics[height=0.1\textheight,width=0.333\linewidth]{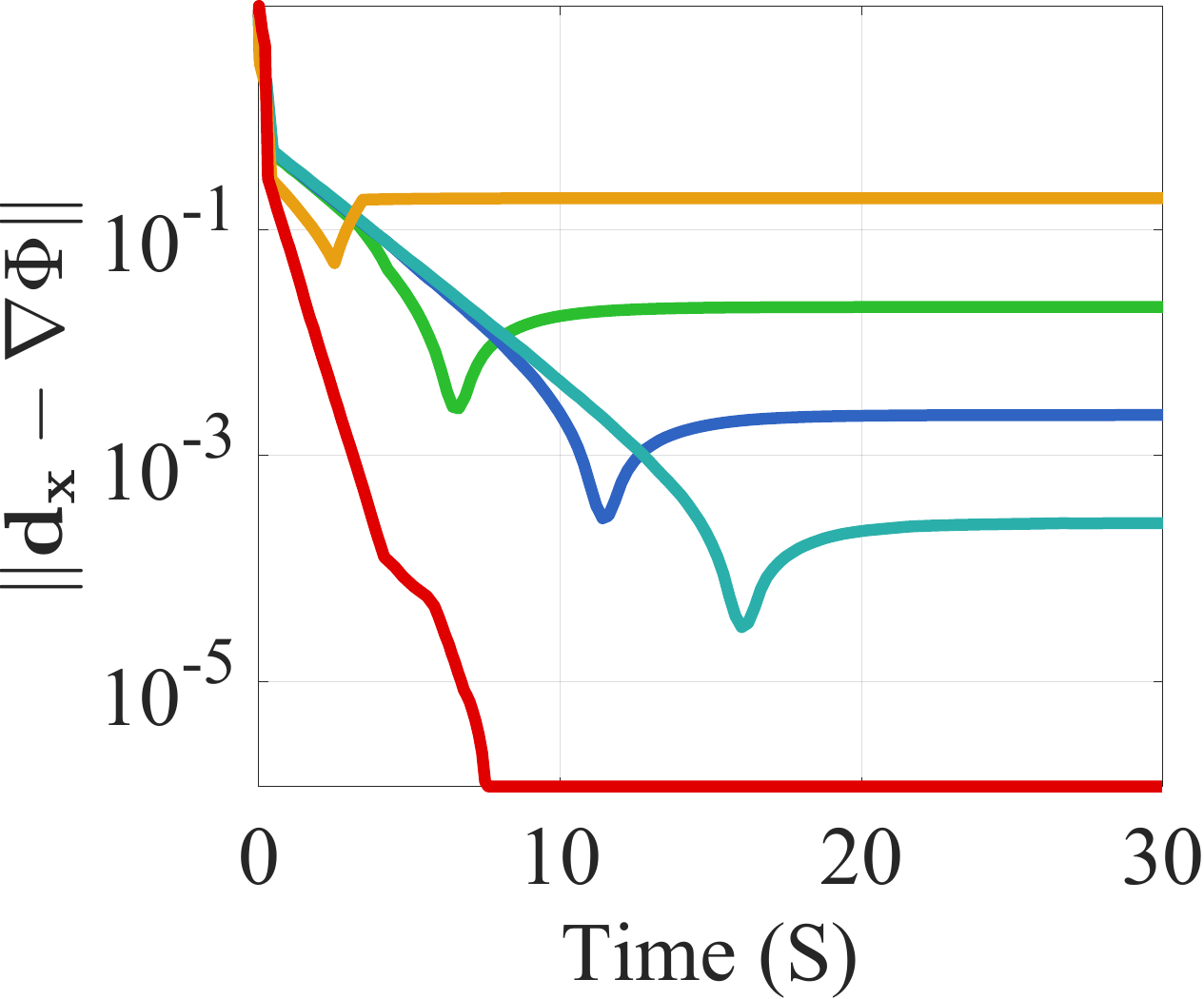}&
			\includegraphics[height=0.1\textheight,width=0.333\linewidth]{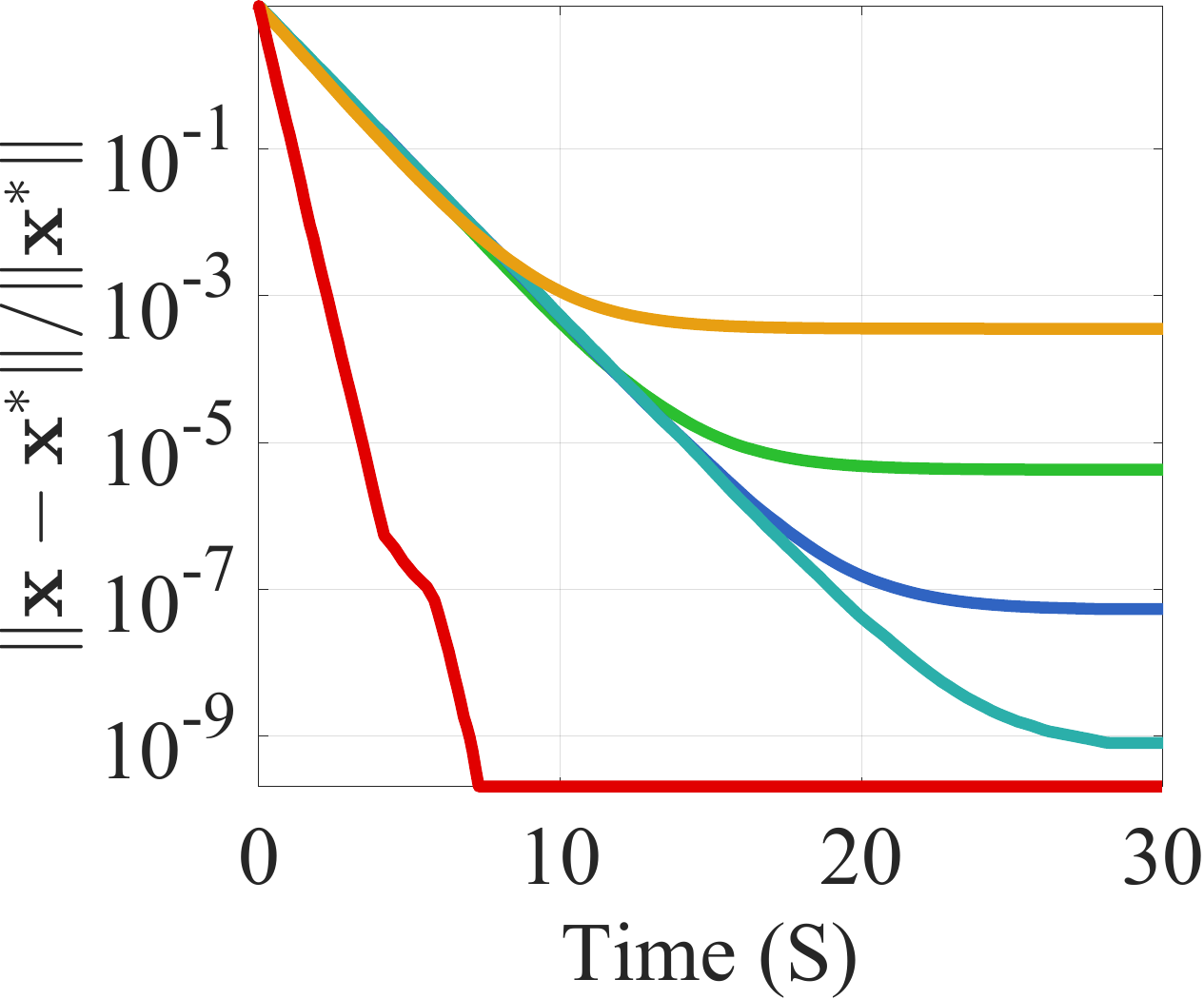}\\
			&\multicolumn{2}{r}{\includegraphics[width=0.6\linewidth,trim=20 0 0 350,clip]{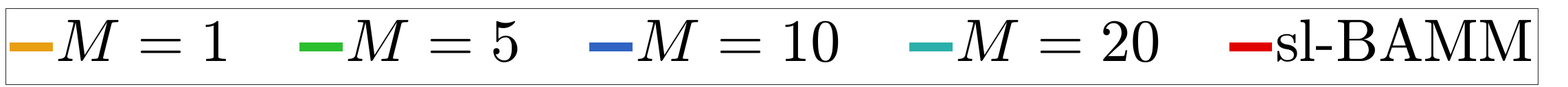}}\\
			\multirow{2}{*}{\rotatebox{90}{NS~~~~~~~~~}}&&\\
			
			&\includegraphics[height=0.1\textheight,width=0.333\linewidth]{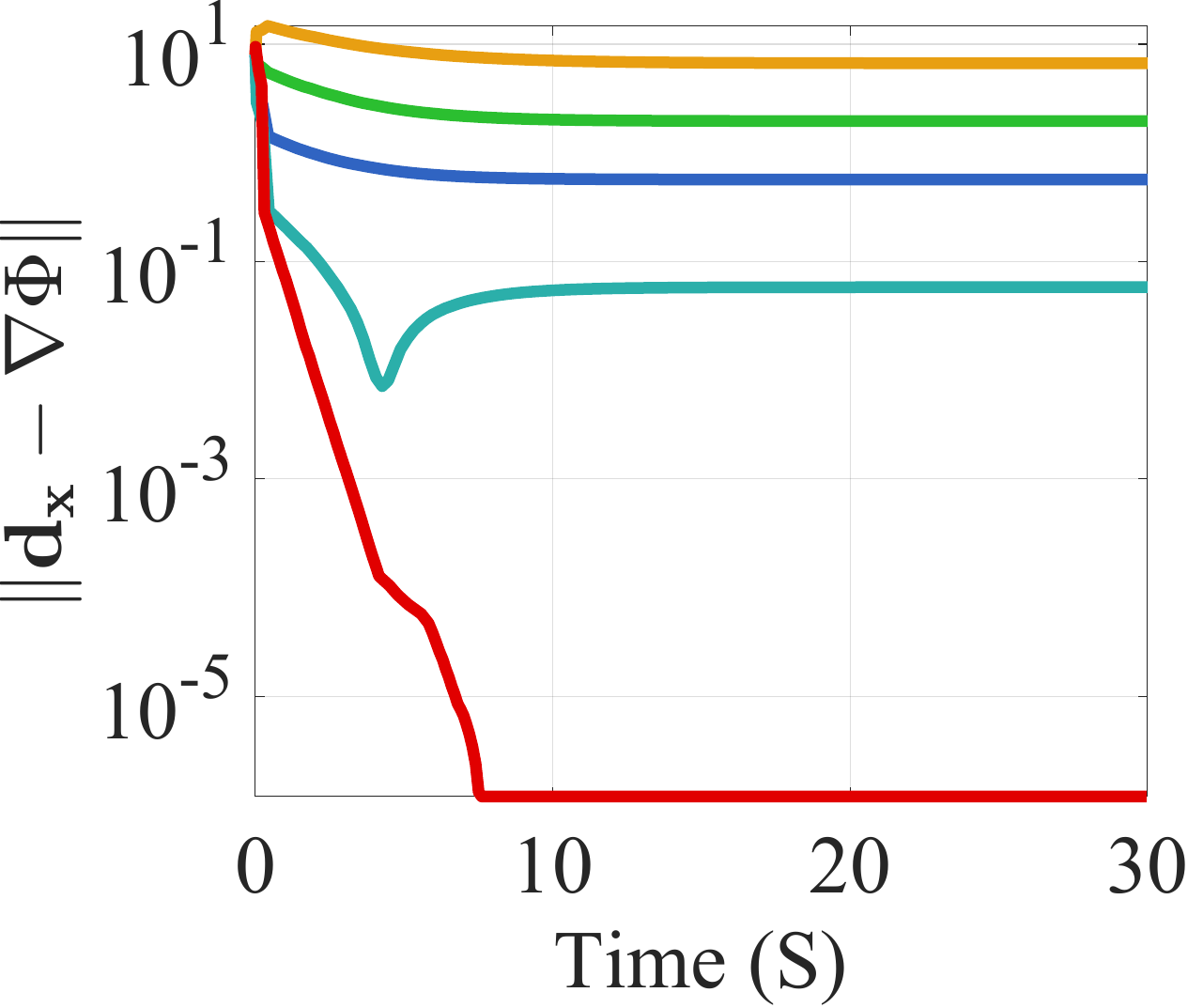}&
			\includegraphics[height=0.1\textheight,width=0.333\linewidth]{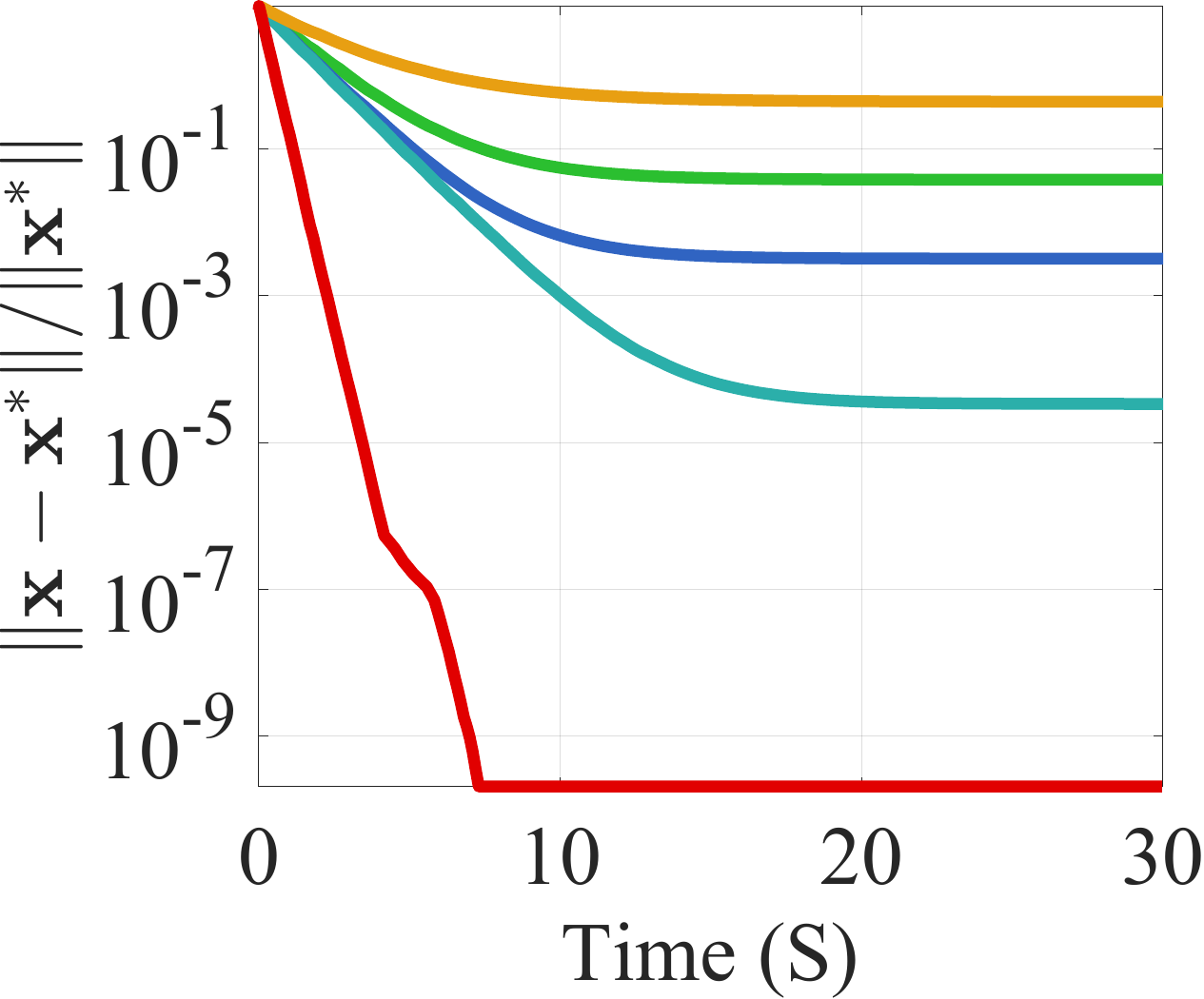}	\\	
		\end{tabular}
	\end{center}
	\vspace{-0.2cm}
	\caption{Convergence curves of sl-BAMM,  RHG with different LL iteration number $T$, CG with different error $\epsilon$ and NS with different sequence length $M$. }\label{fig:LLS_1}
	\vspace{-0.6cm}
\end{figure}

%

\textbf{LL Strongly Convex Case}. Next, we verify the performance of sl-BAMM to handle high-dimensional BLOs compared with existing methods based on the following toy example with LLSC assumption: 
\begin{eqnarray}\label{eq:LLS}
	\begin{aligned}
		&\min _{\x \in \mathbb{R}^n} \frac{1}{2}\|\x-\z_0\|^2+\frac{1}{2}\y^*(\x)^\top\A \y^*(\x) \quad \\&\mathrm { s.t. } \quad\y^*(\x) = \arg \min _{\y \in \mathbb{R}^n} f(\x, \y)=  \frac{1}{2}\y^\top\A\y-\x^\top \y,\\
	\end{aligned}
\end{eqnarray} 
where $\x\in \mathbb{R}^n$, $\y\in \mathbb{R}^n$, $\A\in\mathbb{S}^{n\times n}$ is a positive-definite symmetric matrix and $\z_0\neq0$ is a given point in $\mathbb{R}^n$. It can be easily verified that this example satisfies the assumptions required by RHG, CG and NS. Here, we first set $\A=\I$, $\z_0=\mathbf{e}$, then the unique solution of this problem is $\x^{*}=\y^{*}=\mathbf{e}/2$. Under LLSC assumption, function $\Phi(\x)=F(\x,\y^*(\x))$ is differentiable and we use $\nabla \Phi$ to denote its gradient at $\x$. 
\begin{figure*}[h!]
	\centering
	\begin{tabular}{c@{\extracolsep{1.2em}}c@{\extracolsep{1.2em}}c@{\extracolsep{1.2em}}c@{\extracolsep{1.2em}}}
		\includegraphics[height=0.13\textheight,width=0.18\linewidth]{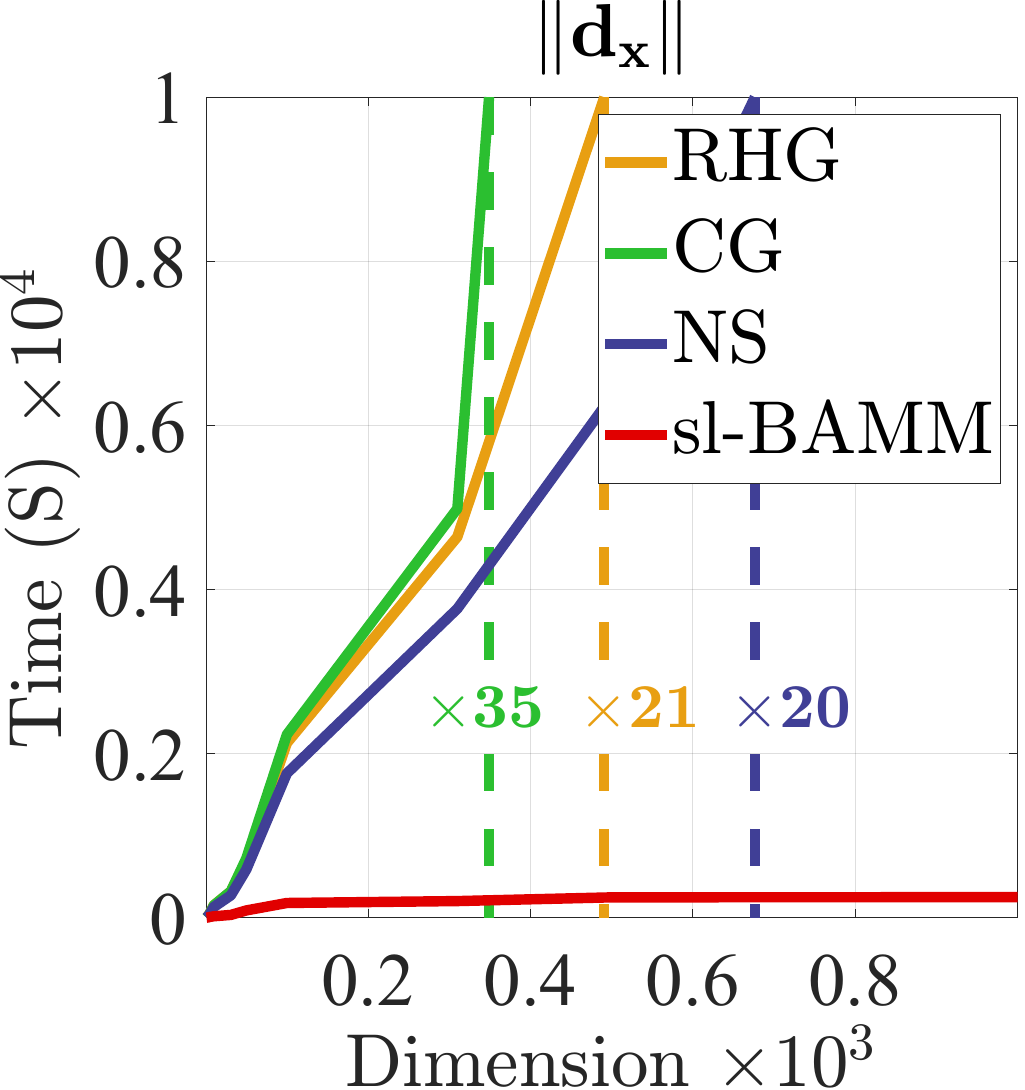}&
		\includegraphics[height=0.13\textheight,width=0.18\linewidth]{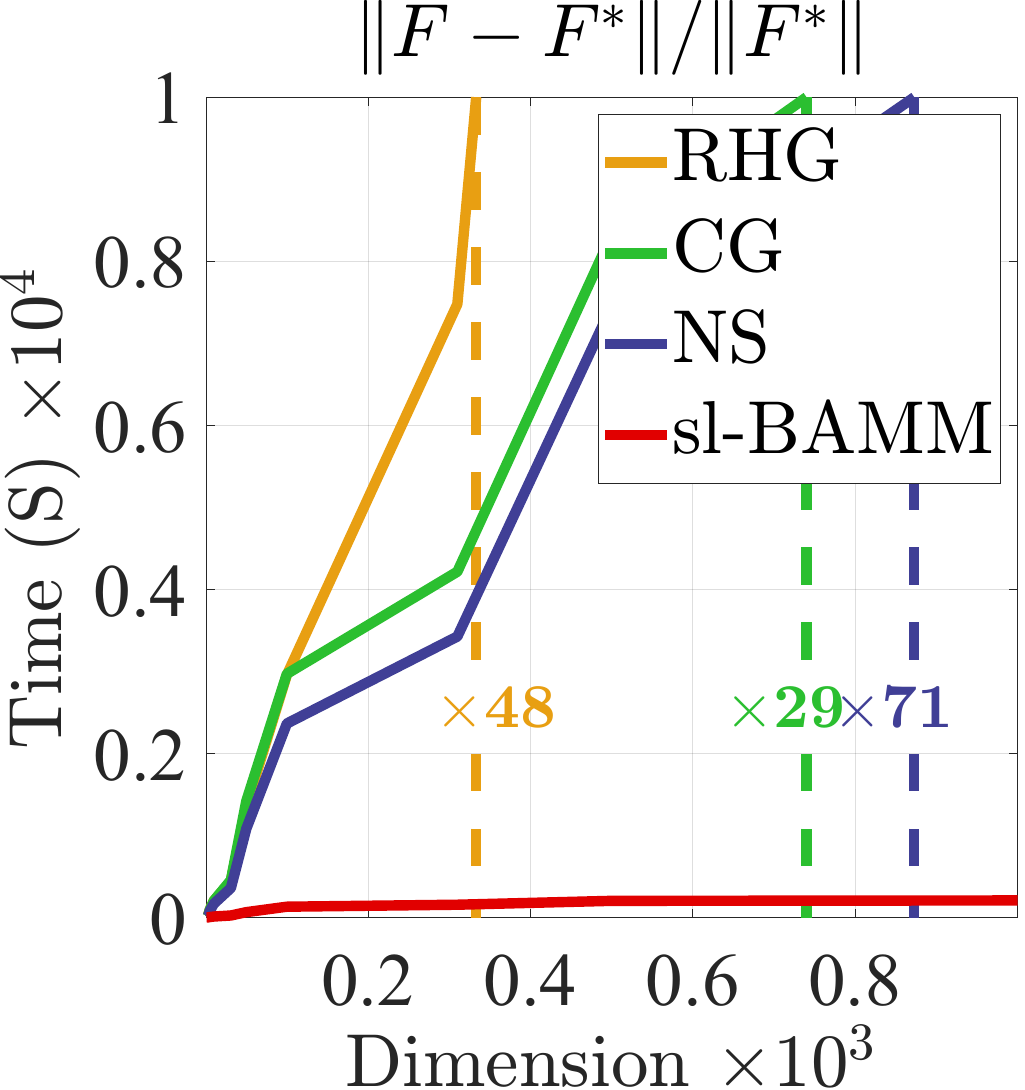}&
		\includegraphics[height=0.13\textheight,width=0.18\linewidth]{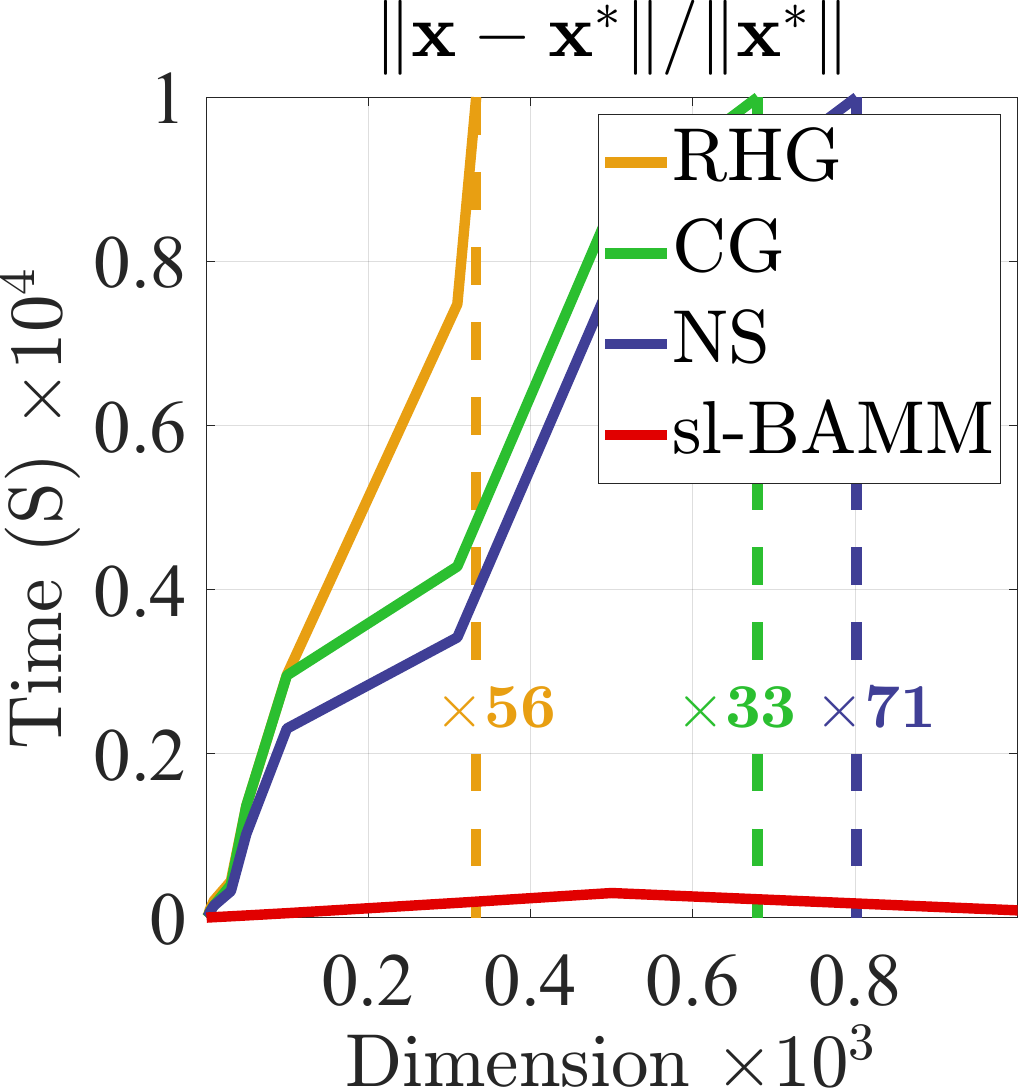}&
		\includegraphics[height=0.13\textheight,width=0.30\linewidth,trim=0 -34 0 0,clip]{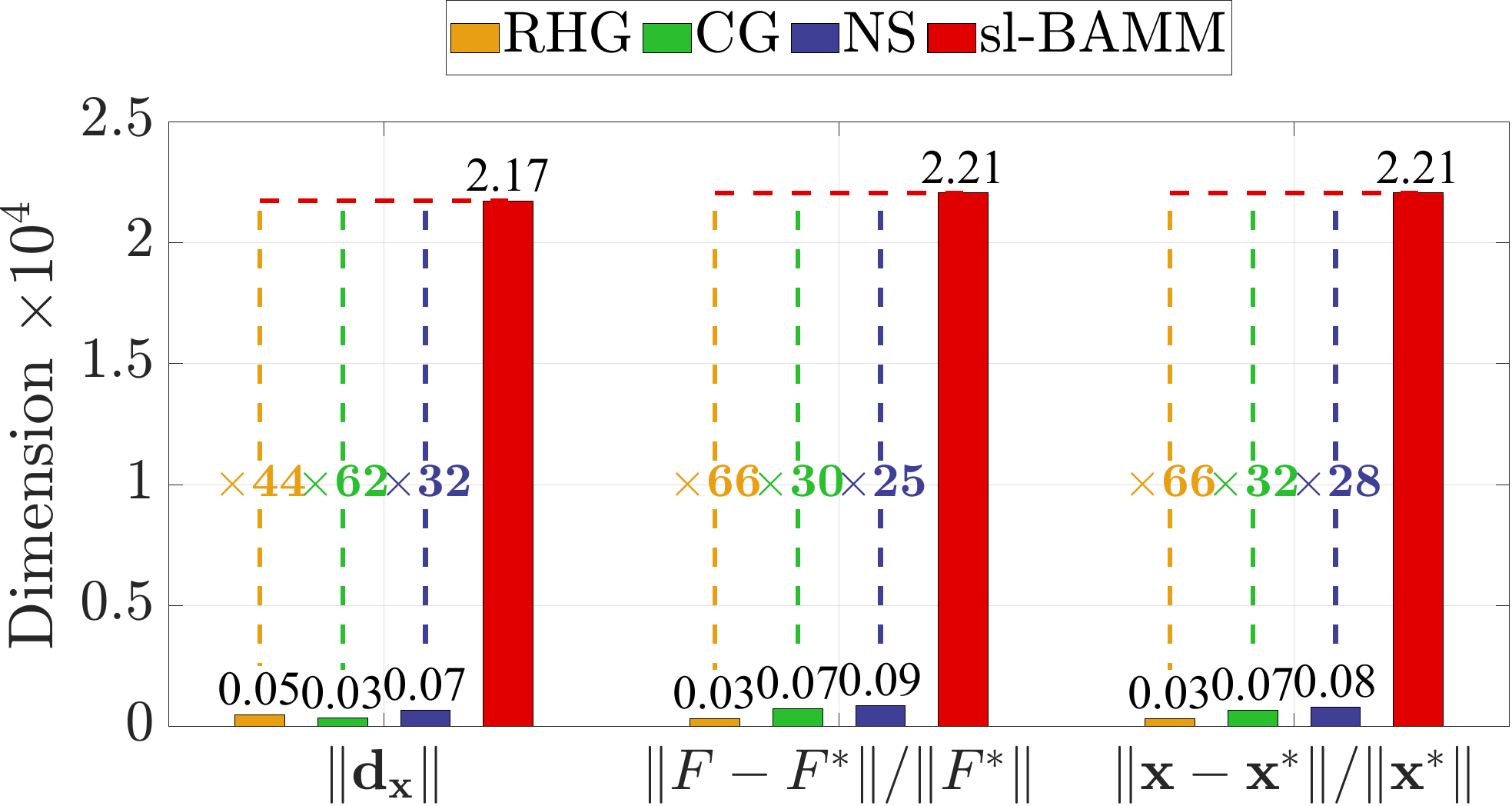}\\
		\multicolumn{3}{c}{\footnotesize (a) Convergence time}&\multicolumn{1}{c}{\footnotesize\quad (b) Convergence dimension}\\
	\end{tabular}
	\caption{Computational efficiency comparison of existing BLOs (RHG, CG and NS) and sl-BAMM. (a) The time required for different convergence metrics as the dimension increases. (b) The problem dimension when the runtime reached the upper limit ($10^{4} $s). }
	\label{fig:LLS_2}
\vspace{-0.2cm}
\end{figure*}
\begin{table*}[h!]
	\centering 
	\makeatletter\def\@captype{table}\makeatother \caption{Comparison of the results for hyper-cleaning and few-shot learning. We compared the time and F1 score for hyper-cleaning when achieving almost the same accuracy ($81\%$) and  the time for few-shot learning ($95\%$ for 5-way and $90\%$ for 20-way).} 
	\vspace{0.2cm}
	\setlength{\tabcolsep}{1.4mm}{
		\begin{tabular}{cccccccccc}
			\toprule
			\multirow{2}{*}{Method}&\multicolumn{3}{c}{Hyper-cleaning}&&&\multicolumn{2}{c}{ 5-way} & \multicolumn{2}{c}{20-way}\\
			&\multirow{1}{*}{Acc. (\%)}&\multirow{1}{*}{F1 score}&\multirow{1}{*}{Time (S)}&&&Acc. (\%) &Time (S)  &Acc. (\%)&Time (S) \\
			\midrule
			RHG&81.01&86.12&95.18&&&  95.01  &  536.17  & 90.27   & 185.28\\
			BDA&81.12&86.37&140.99&&&  95.27  &  297.35 & 90.08  & 425.72\\
			CG&81.08&89.30&35.96&&&  95.01  &  536.17  &  90.07  &  222.72\\
			NS&81.00&87.37&54.78&&&  95.06  &  318.58  & 90.78   &  176.78\\
			sl-BAMM&{81.04}&89.75&\textbf{{3.81}}&&&  95.11  &  \textbf{178.85 }  &  90.13&\textbf{154.06} \\ 
			\bottomrule		
		\end{tabular}
	}\label{tab:hypercleaning}%
\end{table*}
%

In Fig.\ref{fig:LLS_1}, we compare the convergence behavior of our proposed sl-BAMM with existing methods with different hyperparameter settings, including the LL iteration number $T$ for RHG, error $\epsilon$ for CG and sequence length $M$ for NS. 
And we illustrate the error of hypergradient $\Vert\mathbf{d}_{\x}-\nabla \Phi\Vert$ and UL variable $\Vert\x-\x^*\Vert/\Vert\x^*\Vert$ as the calculation time increases. As shown in Fig.\ref{fig:LLS_1}, accelerating RHG, CG and NS by directly simplifying the LL iterations will deteriorate the quality of the iterates generated. Whereas, though sl-BAMM only has one step update for LL problem at each iteration, the update of the introduced multiplier variable $v$ in sl-BAMM  serves as a correction helping the generated iterates converge to the solution.

Next, we test the toy example problem on high-dimensional case to illustrate the computational efficiency by increasing the dimension of $\A$. Here, we set the time limit as $10^4$ s. Fig.~\ref{fig:LLS_2} (a) shows the convergence time of different indicator $\Vert\mathbf{d}_{\x}\Vert\leq 10^{-3}$, $\Vert F-F^*\Vert/\Vert F^*\Vert\leq 10^{-3}$ and $\Vert\x-\x^*\Vert/\Vert\x^*\Vert\leq 10^{-4}$. It can be  observed that sl-BAMM saves much more time compared to RHG, CG and NS, which need multiple optimization iteration step on LL problem. Fig.~\ref{fig:LLS_2} (b) illustrates the upper limit of problem dimension that existing methods can solve when reaching the time limit. Since the implementation of sl-BAMM does not contain repeated nested matrix-vector products, it is capable to handle BLO problems with higher dimension.

\subsection{Real-world Applications}



In this section, we verify the performance of sl-BAMM based on two widely used real-world BLO applications, including data hyper-cleaning and few-shot learning.

\begin{figure}[h!]
	\vspace{-0.1cm}
	\centering
	\begin{tabular}{@{\extracolsep{1.2em}}c@{\extracolsep{1.2em}}c@{\extracolsep{1.2em}}}
		\includegraphics[height=0.135\textheight,width=0.38\linewidth]{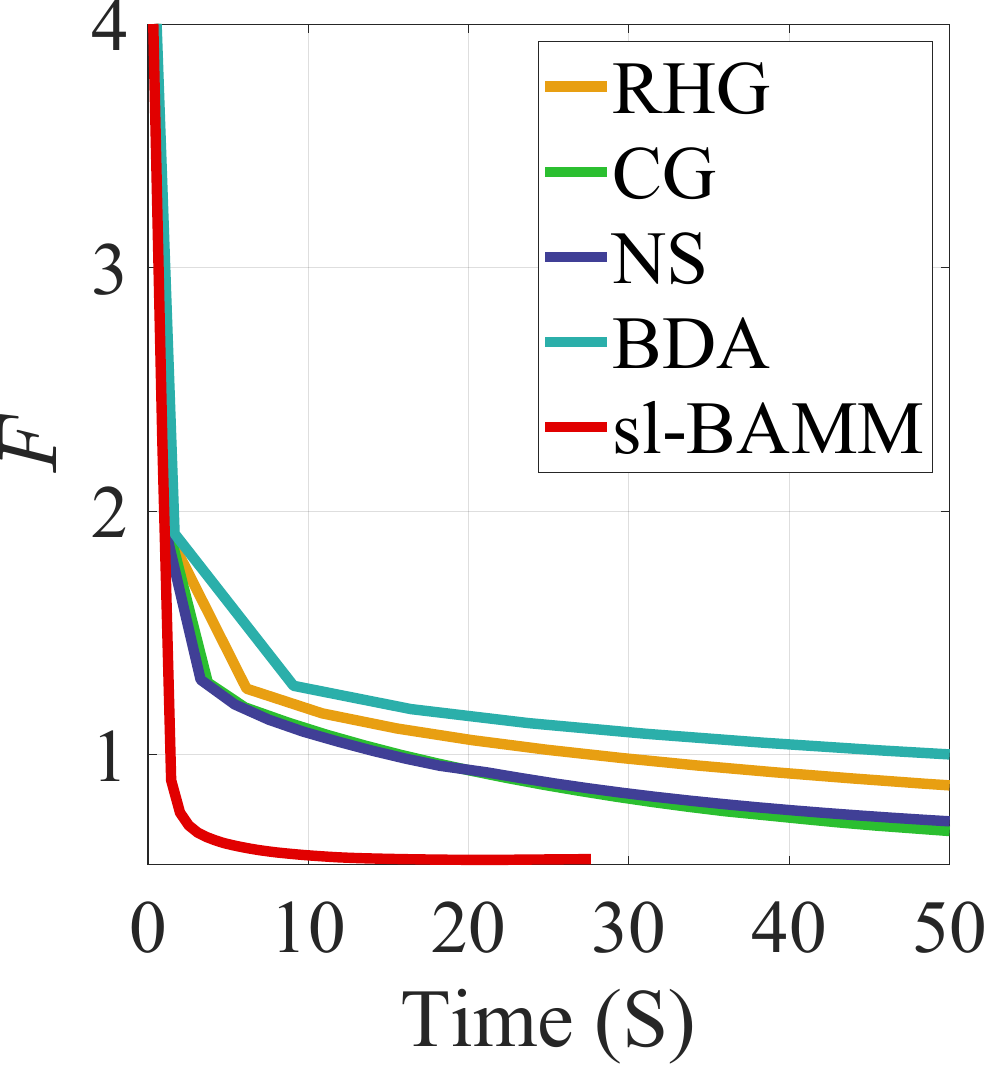}&
		\includegraphics[height=0.132\textheight,width=0.38\linewidth]{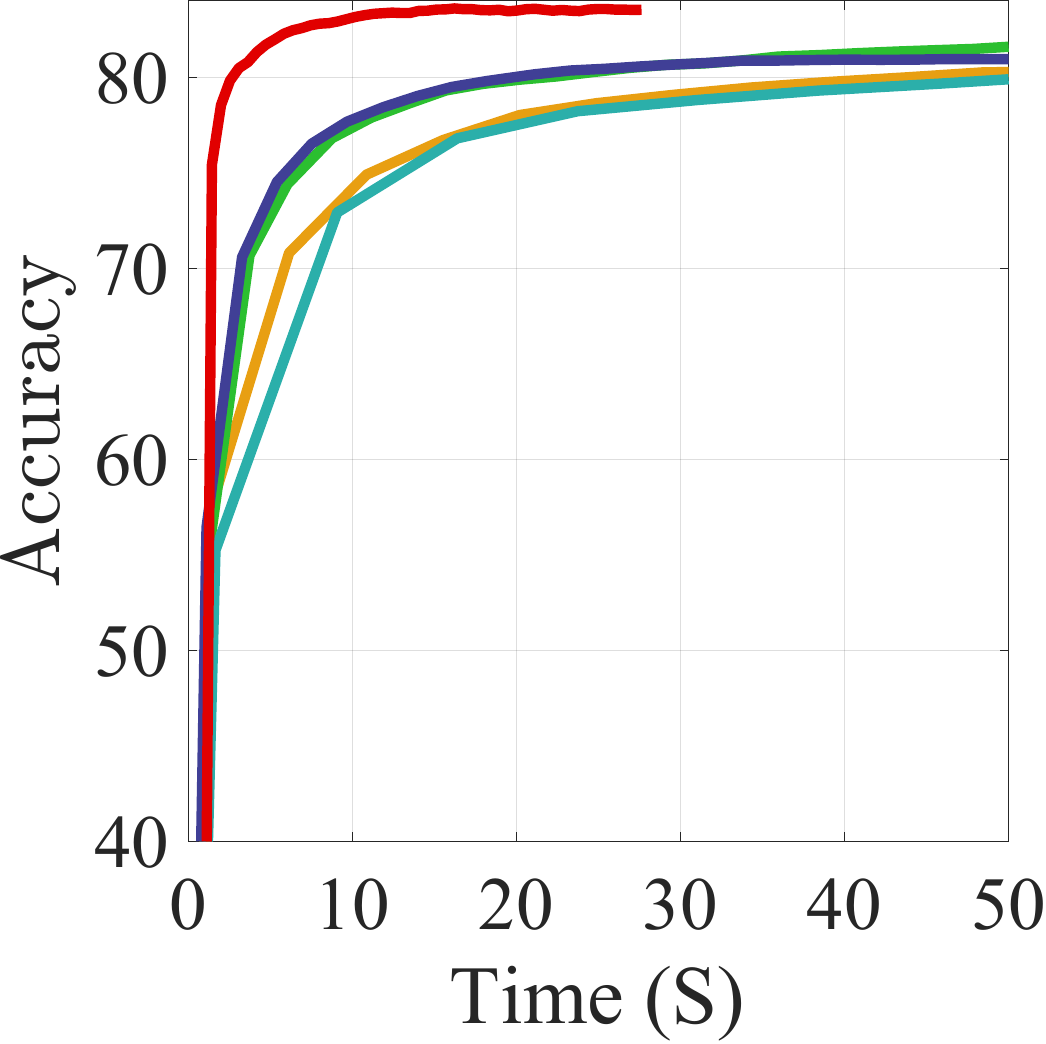}\\
	\end{tabular}
	\caption{Comparison of the validation loss $F$ and accuracy for hyper-cleaning on FashionMNIST dataset. }
	\label{fig:valloss}
	\vspace{-0.6cm}
\end{figure}

\textbf{Data Hyper-cleaning}. Assuming that some of the labels in the dataset are contaminated, data hyper-cleaning aims to reduce the impact of incorrect samples by adding hyper-parameters to label the corrupted data. In the experiment, we follow the general parameter settings in BDA \cite{liu2020generic} and conduct the experiments on FashionMNIST datasets. 
To demonstrate the significant improvement of sl-BAMM in terms of computational efficiency, in Tab.~\ref{tab:hypercleaning}, we report the F1 scores and the time required for different methods to achieve the similar accuracy. It can be seen that sl-BAMM significantly reduces the time needed to achieve desired solution.
Fig.~\ref{fig:valloss}  presents the validation loss and test accuracy for different methods. It can be seen that our method has the fastest convergence rate and keeps being the fastest at higher accuracies even after quickly achieving  $81\%$ accuracy rate.

%

\textbf{Few-shot classification}. The N-way M-shot classification aims to improve the fast-adaptability of hyper model such that the new task can be solved with comparable performance with quick updates based on few examples from each class. In the experiment, we conduct the 5-way 1-shot and 20-way 1-shot experiments on Omniglot dataset~\cite{finn2017model}.  In the right part of Tab.~\ref{tab:hypercleaning} , we report the runtime of existing methods to achieve the same accuracy (95\% for 5-way and 90\% for-20 way). It can be observed that sl-BAMM achieves similar accuracy while  significantly reducing the calculation time, which demonstrates the applicability of sl-BAMM for high dimension complicated real-world deep-learning applications.

\section{Conclusion}

In this paper, a simple yet general single-loop gradient method, named sl-BAMM, is proposed for BLO, without replying on Lower-Level Strong Convexity condition. Our results also provide the theoretical guarantee for various BLOs in the machine learning context where a strong gradient boundedness assumption does not hold. 
We anticipate that the convergence analysis that we develop will be useful for analyzing other BLOs, e.g., stochastic BLOs, and the proposed algorithm will be useful for other applications such as neural architecture search and reinforcement learning. 

\section*{Acknowledgements}
Authors listed in alphabetical order. 
This work is partially supported by the National Natural Science Foundation of China (Nos. U22B2052, 12222106), the Shenzhen Science and Technology Program (No. RCYX20200714114700072), the Guangdong Basic and Applied Basic Research Foundation (No. 2022B1515020082), the Pacific Institute for the Mathematical Sciences (PIMS), the LiaoNing Revitalization Talents Program(No. 2022RG04) and the Fundamental Research Funds for the Central Universities.


\nocite{langley00}

\bibliography{bilevel2023_bib}

\begin{thebibliography}{55}
\providecommand{\natexlab}[1]{#1}
\providecommand{\url}[1]{\texttt{#1}}
\expandafter\ifx\csname urlstyle\endcsname\relax
  \providecommand{\doi}[1]{doi: #1}\else
  \providecommand{\doi}{doi: \begingroup \urlstyle{rm}\Url}\fi

\bibitem[Arbel \& Mairal(2022{\natexlab{a}})Arbel and
  Mairal]{arbel2022amortized}
Arbel, M. and Mairal, J.
\newblock Amortized implicit differentiation for stochastic bilevel
  optimization.
\newblock In \emph{ICLR}, 2022{\natexlab{a}}.

\bibitem[Arbel \& Mairal(2022{\natexlab{b}})Arbel and Mairal]{arbel2022non}
Arbel, M. and Mairal, J.
\newblock Non-convex bilevel games with critical point selection maps.
\newblock In \emph{NeurIPS}, 2022{\natexlab{b}}.

\bibitem[Beck(2017)]{beck2017first}
Beck, A.
\newblock \emph{First-order methods in optimization}.
\newblock SIAM, 2017.

\bibitem[Chen et~al.(2023)Chen, Xu, and Zhang]{chen2023bilevel}
Chen, L., Xu, J., and Zhang, J.
\newblock On bilevel optimization without lower-level strong convexity.
\newblock \emph{arXiv preprint arXiv:2301.00712}, 2023.

\bibitem[Chen et~al.(2021)Chen, Sun, and Yin]{chen2021closing}
Chen, T., Sun, Y., and Yin, W.
\newblock Closing the gap: Tighter analysis of alternating stochastic gradient
  methods for bilevel problems.
\newblock \emph{NeurIPS}, 34:\penalty0 25294--25307, 2021.

\bibitem[Chen et~al.(2022)Chen, Sun, Xiao, and Yin]{chen2022single}
Chen, T., Sun, Y., Xiao, Q., and Yin, W.
\newblock A single-timescale method for stochastic bilevel optimization.
\newblock In \emph{AISTATS}, 2022.

\bibitem[Chen et~al.(2019)Chen, Xie, Wu, and Tian]{chen2019progressive}
Chen, X., Xie, L., Wu, J., and Tian, Q.
\newblock Progressive differentiable architecture search: Bridging the depth
  gap between search and evaluation.
\newblock In \emph{ICCV}, 2019.

\bibitem[Dagr{\'e}ou et~al.(2022)Dagr{\'e}ou, Ablin, Vaiter, and
  Moreau]{dagreou2022framework}
Dagr{\'e}ou, M., Ablin, P., Vaiter, S., and Moreau, T.
\newblock A framework for bilevel optimization that enables stochastic and
  global variance reduction algorithms.
\newblock In \emph{NeurIPS}, 2022.

\bibitem[Elsken et~al.(2020)Elsken, Staffler, Metzen, and
  Hutter]{elsken2020meta}
Elsken, T., Staffler, B., Metzen, J.~H., and Hutter, F.
\newblock Meta-learning of neural architectures for few-shot learning.
\newblock In \emph{CVPR}, 2020.

\bibitem[Finn et~al.(2017)Finn, Abbeel, and Levine]{finn2017model}
Finn, C., Abbeel, P., and Levine, S.
\newblock Model-agnostic meta-learning for fast adaptation of deep networks.
\newblock In \emph{ICML}, 2017.

\bibitem[Franceschi et~al.(2017)Franceschi, Donini, Frasconi, and
  Pontil]{franceschi2017forward}
Franceschi, L., Donini, M., Frasconi, P., and Pontil, M.
\newblock Forward and reverse gradient-based hyperparameter optimization.
\newblock In \emph{ICML}, 2017.

\bibitem[Franceschi et~al.(2018)Franceschi, Frasconi, Salzo, Grazzi, and
  Pontil]{franceschi2018bilevel}
Franceschi, L., Frasconi, P., Salzo, S., Grazzi, R., and Pontil, M.
\newblock Bilevel programming for hyperparameter optimization and
  meta-learning.
\newblock In \emph{ICML}, 2018.

\bibitem[Ghadimi \& Wang(2018)Ghadimi and Wang]{ghadimi2018approximation}
Ghadimi, S. and Wang, M.
\newblock Approximation methods for bilevel programming.
\newblock \emph{arXiv preprint arXiv:1802.02246}, 2018.

\bibitem[Gong et~al.(2021)Gong, Liu, and Liu]{gong2021automatic}
Gong, C., Liu, X., and Liu, Q.
\newblock Automatic and harmless regularization with constrained and
  lexicographic optimization: A dynamic barrier approach.
\newblock \emph{NeurIPS}, 34:\penalty0 29630--29642, 2021.

\bibitem[Goodfellow et~al.(2014)Goodfellow, Pouget-Abadie, Mirza, Xu,
  Warde-Farley, Ozair, Courville, and Bengio]{goodfellow2014generative}
Goodfellow, I., Pouget-Abadie, J., Mirza, M., Xu, B., Warde-Farley, D., Ozair,
  S., Courville, A., and Bengio, Y.
\newblock Generative adversarial nets.
\newblock \emph{NeurIPS}, 2014.

\bibitem[Grazzi et~al.(2020)Grazzi, Franceschi, Pontil, and
  Salzo]{grazzi2020iteration}
Grazzi, R., Franceschi, L., Pontil, M., and Salzo, S.
\newblock On the iteration complexity of hypergradient computation.
\newblock In \emph{ICML}, 2020.

\bibitem[Hong et~al.(2023)Hong, Wai, Wang, and Yang]{hong2020two}
Hong, M., Wai, H.-T., Wang, Z., and Yang, Z.
\newblock A two-timescale framework for bilevel optimization: Complexity
  analysis and application to actor-critic.
\newblock \emph{SIAM Journal on Optimization}, 33\penalty0 (1):\penalty0
  147--180, 2023.

\bibitem[Huang et~al.(2022)Huang, Li, Gao, and Huang]{huang2021enhanced}
Huang, F., Li, J., Gao, S., and Huang, H.
\newblock Enhanced bilevel optimization via bregman distance.
\newblock \emph{NeurIPS}, 35:\penalty0 28928--28939, 2022.

\bibitem[Ji \& Liang(2022)Ji and Liang]{ji2021lower}
Ji, K. and Liang, Y.
\newblock Lower bounds and accelerated algorithms for bilevel optimization.
\newblock \emph{Journal of Machine Learning Research}, 23:\penalty0 1--56,
  2022.

\bibitem[Ji et~al.(2020)Ji, Lee, Liang, and Poor]{KaiyiJi2020ConvergenceOM}
Ji, K., Lee, J.~D., Liang, Y., and Poor, H.~V.
\newblock Convergence of meta-learning with task-specific adaptation over
  partial parameters.
\newblock \emph{NeurIPS}, 33:\penalty0 11490--11500, 2020.

\bibitem[Ji et~al.(2021)Ji, Yang, and Liang]{ji2021bilevel}
Ji, K., Yang, J., and Liang, Y.
\newblock Bilevel optimization: Convergence analysis and enhanced design.
\newblock In \emph{ICML}, 2021.

\bibitem[Ji et~al.(2022)Ji, Liu, Liang, and Ying]{ji2022will}
Ji, K., Liu, M., Liang, Y., and Ying, L.
\newblock Will bilevel optimizers benefit from loops.
\newblock In \emph{NeurIPS}, 2022.

\bibitem[Khanduri et~al.(2021)Khanduri, Zeng, Hong, Wai, Wang, and
  Yang]{khanduri2021near}
Khanduri, P., Zeng, S., Hong, M., Wai, H.-T., Wang, Z., and Yang, Z.
\newblock A near-optimal algorithm for stochastic bilevel optimization via
  double-momentum.
\newblock \emph{NeurIPS}, 2021.

\bibitem[Kim et~al.(2020)Kim, Leyffer, and Munson]{kim2020mpec}
Kim, Y., Leyffer, S., and Munson, T.
\newblock Mpec methods for bilevel optimization problems.
\newblock In \emph{Bilevel Optimization}, pp.\  335--360. Springer, 2020.

\bibitem[Li et~al.(2020)Li, Gu, and Huang]{li2020improved}
Li, J., Gu, B., and Huang, H.
\newblock Improved bilevel model: Fast and optimal algorithm with theoretical
  guarantee.
\newblock \emph{arXiv preprint arXiv:2009.00690}, 2020.

\bibitem[Li et~al.(2022)Li, Gu, and Huang]{li2022fully}
Li, J., Gu, B., and Huang, H.
\newblock A fully single loop algorithm for bilevel optimization without
  hessian inverse.
\newblock In \emph{AAAI}, volume~36, pp.\  7426--7434, 2022.

\bibitem[Liang et~al.(2019)Liang, Zhang, Sun, He, Huang, Zhuang, and
  Li]{liang2019darts+}
Liang, H., Zhang, S., Sun, J., He, X., Huang, W., Zhuang, K., and Li, Z.
\newblock Darts+: Improved differentiable architecture search with early
  stopping.
\newblock \emph{arXiv preprint arXiv:1909.06035}, 2019.

\bibitem[Liu et~al.(2018)Liu, Simonyan, and Yang]{liu2018darts}
Liu, H., Simonyan, K., and Yang, Y.
\newblock Darts: Differentiable architecture search.
\newblock In \emph{ICLR}, 2018.

\bibitem[Liu et~al.(2020)Liu, Mu, Yuan, Zeng, and Zhang]{liu2020generic}
Liu, R., Mu, P., Yuan, X., Zeng, S., and Zhang, J.
\newblock A generic first-order algorithmic framework for bi-level programming
  beyond lower-level singleton.
\newblock In \emph{ICML}, 2020.

\bibitem[Liu et~al.(2021{\natexlab{a}})Liu, Gao, Zhang, Meng, and
  Lin]{liu2021investigating}
Liu, R., Gao, J., Zhang, J., Meng, D., and Lin, Z.
\newblock Investigating bi-level optimization for learning and vision from a
  unified perspective: A survey and beyond.
\newblock \emph{IEEE TPAMI}, 2021{\natexlab{a}}.

\bibitem[Liu et~al.(2021{\natexlab{b}})Liu, Liu, Yuan, Zeng, and
  Zhang]{liu2021value}
Liu, R., Liu, X., Yuan, X., Zeng, S., and Zhang, J.
\newblock A value-function-based interior-point method for non-convex bi-level
  optimization.
\newblock In \emph{ICML}, pp.\  6882--6892. PMLR, 2021{\natexlab{b}}.

\bibitem[Liu et~al.(2021{\natexlab{c}})Liu, Liu, Zeng, Zhang, and
  Zhang]{liu2021value2}
Liu, R., Liu, X., Zeng, S., Zhang, J., and Zhang, Y.
\newblock Value-function-based sequential minimization for bi-level
  optimization.
\newblock \emph{arXiv preprint arXiv:2110.04974}, 2021{\natexlab{c}}.

\bibitem[Liu et~al.(2021{\natexlab{d}})Liu, Liu, Zeng, and
  Zhang]{liu2021towards}
Liu, R., Liu, Y., Zeng, S., and Zhang, J.
\newblock Towards gradient-based bilevel optimization with non-convex followers
  and beyond.
\newblock \emph{NeurIPS}, 34:\penalty0 8662--8675, 2021{\natexlab{d}}.

\bibitem[Liu et~al.(2022{\natexlab{a}})Liu, Liu, Zeng, Zhang, and
  Zhang]{liu2022optimization}
Liu, R., Liu, X., Zeng, S., Zhang, J., and Zhang, Y.
\newblock Optimization-derived learning with essential convergence analysis of
  training and hyper-training.
\newblock In \emph{ICML}, pp.\  13825--13856. PMLR, 2022{\natexlab{a}}.

\bibitem[Liu et~al.(2022{\natexlab{b}})Liu, Mu, Yuan, Zeng, and
  Zhang]{liu2022general}
Liu, R., Mu, P., Yuan, X., Zeng, S., and Zhang, J.
\newblock A general descent aggregation framework for gradient-based bi-level
  optimization.
\newblock \emph{IEEE TPAMI}, 2022{\natexlab{b}}.

\bibitem[Liu et~al.(2023{\natexlab{a}})Liu, Liu, Zeng, Zhang, and
  Zhang]{liu2023hierarchical}
Liu, R., Liu, X., Zeng, S., Zhang, J., and Zhang, Y.
\newblock Hierarchical optimization-derived learning.
\newblock \emph{arXiv preprint arXiv:2302.05587}, 2023{\natexlab{a}}.

\bibitem[Liu et~al.(2023{\natexlab{b}})Liu, Liu, Zeng, and
  Zhang]{liu2023augmenting}
Liu, R., Liu, Y., Zeng, S., and Zhang, J.
\newblock Augmenting iterative trajectory for bilevel optimization:
  Methodology, analysis and extensions.
\newblock \emph{arXiv preprint arXiv:2303.16397}, 2023{\natexlab{b}}.

\bibitem[Lorraine et~al.(2020)Lorraine, Vicol, and
  Duvenaud]{lorraine2020optimizing}
Lorraine, J., Vicol, P., and Duvenaud, D.
\newblock Optimizing millions of hyperparameters by implicit differentiation.
\newblock In \emph{AISTATS}, 2020.

\bibitem[Mackay et~al.(2019)Mackay, Vicol, Lorraine, Duvenaud, and
  Grosse]{mackay2018self}
Mackay, M., Vicol, P., Lorraine, J., Duvenaud, D., and Grosse, R.
\newblock Self-tuning networks: Bilevel optimization of hyperparameters using
  structured best-response functions.
\newblock In \emph{ICLR}, 2019.

\bibitem[Maclaurin et~al.(2015)Maclaurin, Duvenaud, and
  Adams]{maclaurin2015gradient}
Maclaurin, D., Duvenaud, D., and Adams, R.
\newblock Gradient-based hyperparameter optimization through reversible
  learning.
\newblock In \emph{ICML}, 2015.

\bibitem[Mehra \& Hamm(2021)Mehra and Hamm]{mehra2021penalty}
Mehra, A. and Hamm, J.
\newblock Penalty method for inversion-free deep bilevel optimization.
\newblock In \emph{Asian Conference on Machine Learning}, pp.\  347--362. PMLR,
  2021.

\bibitem[Nesterov(2018)]{nesterov2018lectures}
Nesterov, Y.
\newblock \emph{Lectures on convex optimization}, volume 137.
\newblock Springer, 2018.

\bibitem[Okuno et~al.(2018)Okuno, Takeda, and Kawana]{okuno2018hyperparameter}
Okuno, T., Takeda, A., and Kawana, A.
\newblock Hyperparameter learning via bilevel nonsmooth optimization.
\newblock \emph{arXiv preprint arXiv:1806.01520}, 2018.

\bibitem[Outrata(1990)]{outrata1990numerical}
Outrata, J.~V.
\newblock On the numerical solution of a class of stackelberg problems.
\newblock \emph{Zeitschrift f{\"u}r Operations Research}, 34\penalty0
  (4):\penalty0 255--277, 1990.

\bibitem[Pedregosa(2016)]{pedregosa2016hyperparameter}
Pedregosa, F.
\newblock Hyperparameter optimization with approximate gradient.
\newblock In \emph{ICML}, 2016.

\bibitem[Pfau \& Vinyals(2016)Pfau and Vinyals]{pfau2016connecting}
Pfau, D. and Vinyals, O.
\newblock Connecting generative adversarial networks and actor-critic methods.
\newblock \emph{arXiv preprint arXiv:1610.01945}, 2016.

\bibitem[Rajeswaran et~al.(2019)Rajeswaran, Finn, Kakade, and
  Levine]{rajeswaran2019meta}
Rajeswaran, A., Finn, C., Kakade, S.~M., and Levine, S.
\newblock Meta-learning with implicit gradients.
\newblock In \emph{NeurIPS}, 2019.

\bibitem[Sabach \& Shtern(2017)Sabach and Shtern]{sabach2017first}
Sabach, S. and Shtern, S.
\newblock A first order method for solving convex bilevel optimization
  problems.
\newblock \emph{SIAM Journal on Optimization}, 27\penalty0 (2):\penalty0
  640--660, 2017.

\bibitem[Shaban et~al.(2019)Shaban, Cheng, Hatch, and
  Boots]{shaban2019truncated}
Shaban, A., Cheng, C.-A., Hatch, N., and Boots, B.
\newblock Truncated back-propagation for bilevel optimization.
\newblock In \emph{AISTATS}, 2019.

\bibitem[Sow et~al.(2022)Sow, Ji, Guan, and Liang]{sow2022constrained}
Sow, D., Ji, K., Guan, Z., and Liang, Y.
\newblock A constrained optimization approach to bilevel optimization with
  multiple inner minima.
\newblock \emph{arXiv preprint arXiv:2203.01123}, 2022.

\bibitem[Yang et~al.(2019)Yang, Chen, Hong, and Wang]{yang2019provably}
Yang, Z., Chen, Y., Hong, M., and Wang, Z.
\newblock Provably global convergence of actor-critic: A case for linear
  quadratic regulator with ergodic cost.
\newblock In \emph{NeurIPS}, 2019.

\bibitem[Ye \& Zhu(1995)Ye and Zhu]{ye1995optimality}
Ye, J. and Zhu, D.
\newblock Optimality conditions for bilevel programming problems.
\newblock \emph{Optimization}, 33\penalty0 (1):\penalty0 9--27, 1995.

\bibitem[Ye et~al.(2022)Ye, Liu, Wright, Stone, and Liu]{NeurIPS2022-Liu}
Ye, M., Liu, B., Wright, S., Stone, P., and Liu, Q.
\newblock Bome! bilevel optimization made easy: A simple first-order approach.
\newblock In \emph{NeurIPS}, 2022.

\bibitem[Zhang et~al.(2020)Zhang, Lin, Jegelka, Sra, and
  Jadbabaie]{zhang2020complexity}
Zhang, J., Lin, H., Jegelka, S., Sra, S., and Jadbabaie, A.
\newblock Complexity of finding stationary points of nonsmooth nonconvex
  functions.
\newblock In \emph{ICML}, pp.\  11173--11182, 2020.

\bibitem[Z{\"u}gner \& G{\"u}nnemann(2018)Z{\"u}gner and
  G{\"u}nnemann]{zugner2018adversarial}
Z{\"u}gner, D. and G{\"u}nnemann, S.
\newblock Adversarial attacks on graph neural networks via meta learning.
\newblock In \emph{ICLR}, 2018.

\end{thebibliography}
\bibliographystyle{icml2023}

 
\newpage
\appendix
\onecolumn
\section{Expanded Related Work}\label{appendix1}

We provide here a detailed review of recent studies that are closely related to ours. 

{\bf BLO Beyond Lower-Level Strong Convexity (LLSC).} 
The studies beyond LLSC are rather limited. A line of work (including this work) relaxed LLS by considering BLO problem with a merely convex LL objective. For example,
by formulating bi-level models from the optimistic viewpoint and aggregating UL and LL objectives (called sequential averaging method in \cite{sabach2017first}), 
the authors in \cite{liu2020generic,liu2022general} introduced Bilevel Descent Aggregation (BDA) framework with only the asymptotic convergence guarantee; see also a concurrent work in \cite{li2020improved}. Another study \cite{mehra2021penalty} presented a penalty method for BLO, by replacing LL problem by its first-order stationarity condition, which avoids computing Hessian inverse. Recently, based on the value function approach, by using the gradient descent and ascent method via smooth approximation, Sow et al.   \cite{sow2022constrained} proposed a primal-dual bilevel optimization (PDBO) algorithm and its proximal version (Proximal-PDBO) without involving second-order Hessian and Jacobian computations. Further, they provided the convergence rate analysis for PDBO and Proximal-PDBO, under certain convexity-type conditions and compactness of the domain. 

Another line of work considers BLO problem with a nonconvex LL objective. For example, the authors in \cite{liu2021value,liu2021value2} proposed a bi-level value function based penalty and barrier method. 
Another study \cite{liu2021towards} proposed an initialization auxiliary gradient-based algorithm to solve BLO problems with nonconvex LL objectives. Recently, Arbel and Mairal \cite{arbel2022non} introduced a bilevel game that resolves the ambiguity in BLO with nonconvex objectives using the notion of selection maps. Further, they extended implicit differentiation to degenerate critical points, provided LL objective satisfies a generalization of the Morse-Bott property. However, only asymptotic result is shown in \cite{liu2021value,liu2021value2,liu2021towards,arbel2022non}. Recently, based on the value function approach and a dynamic barrier gradient descent algorithm in \cite{gong2021automatic}, the authors in \cite{NeurIPS2022-Liu} proposed a simple first-order BLO algorithm that depends only on first-order gradient information and provided non-asymptotic convergence analysis for non-convex objectives under the global or local Polyak-Łojasiewicz conditions. 

However, all of these approaches take a double-loop optimization structure (involving numerical LL optimization loops), which could be difficult to implement in practice as the number of required inner loop iterations are difficult to adjust. Differently from the above studies, our sl-BAMM algorithm has a single loop structure, which admits a much simpler implementation.  Moreover, the updates in sl-BAMM evolve simultaneously and then can be performed in parallel instead of sequentially.

{\bf Stationarity Measure for BLO Algorithm.} 
The most common stationarity measure for BLO is the gradient norm of the total UL objective (i.e., the norm of hypergradient) \cite{pedregosa2016hyperparameter,ghadimi2018approximation,chen2021closing,ji2021bilevel,ji2022will,li2022fully,dagreou2022framework,arbel2022amortized}. These studies notoriously rely on LLSC. 
This was extended in \cite{huang2021enhanced} to a general BLO problem with a nonsmooth UL objective via Bregman distance, but still assuming LLSC. Further, when the total UL objective is either strongly convex or convex, the suboptimality gap (as measured by the total UL objective values) is used in \cite{ghadimi2018approximation,hong2020two,ji2021lower} and the optimality gap w.r.t. UL variable is used in \cite{hong2020two}
. 
Recently, under the weak convexity of the total UL objective, the gradient norm of the standard Moreau envelope is used in \cite{hong2020two,chen2022single}. 

Without LL strongly convex assumption, the stationarity measure for BLO has been much less studied, from an algorithmic point of view.  
Several recent work \cite{sow2022constrained,NeurIPS2022-Liu} reformulate a BLO problem into a constrained problem by using value function approach, and develop non-asymptotic convergence to various notions of approximate KKT points, under stationarity measure involving value function \cite{NeurIPS2022-Liu} or regularized value function \cite{sow2022constrained}. Very recently, inspired by the recent developments of nonsmooth nonconvex optimization, the authors in \cite{chen2023bilevel} formalized the local optimality of BLO via the notion of Goldstein stationarity condition \cite{zhang2020complexity} of the hyper-objective (i.e., the value function of a parameterized simple BLO problem),  under some growth conditions. When LL problem is unconstrained and be convex, BLO problem is equivalent to an equality-constrained optimization problem by replacing LL problem by its first-order stationarity condition. This motivates us to use KKT condition of the resulting equality-constrained problem as a measure of stationarity of the solution returned by BLO algorithms. Note that this kind of reformulation is called MPEC formulation in the optimization theory for general BLO problems with LL inequality and equality constraints \cite{kim2020mpec}. 

\section{A Unified Proof Sketch of Theorems in Section~\ref{Convergence Analysis}}
\label{proofsketch}

Theorems~\ref{thmconvex1}-\ref{thm_Strong_convex} in Section~\ref{Convergence Analysis} provide various practical step size strategies for sl-BAMM, which enjoy favorably both extremely simpler implementation and theoretical convergence guarantee. 
In this section, we highlight the key steps of the proofs towards Theorems~\ref{thmconvex1}-\ref{thm_Strong_convex} {\bf without the gradient boundedness assumption (GBA) and beyond lower-level strong convexity (LLSC)}, and highlight the differences between our new convergence analysis and the existing ones. 

To analyze the convergence of sl-BAMM with or without LLSC, according to different step size strategies, we will identify different intrinsic Lyapunov (potential) functions, whose general form is given by
\begin{equation}\label{Lyapunov}
	V_k
	:=a_k \left[ F(\x_k,\y^*_{\mu_k}(\x_k)) - F_0\right]
	+b_k \|\y_k-\y^*_{\mu_k}(\x_k)\|^2
	+ c_k \|\v_k-\v^*_{\mu_k}(\x_k)\|^2,
\end{equation}
where $\y^*_{\mu}(\x)$ is the unique minimizer of the aggregation function $\psi_{\mu}(\x,\cdot)=\mu F(\x,\cdot)+(1-\mu)f(\x,\cdot)$, and  $\v^*_{\mu}(\x)=\left[\nabla_{\y\y}^2 \psi_{\mu}(\x,\y^*_{\mu}(\x))\right]^{-1}\nabla_{\y} F\left(\x,\y^*_{\mu}(\x)\right)$ is the correct dual multiplier. Here $\{ a_k, b_k, c_k \}_{k=1}^\infty$ is a sequence of positive constants chosen later depending on the step size strategy. In particular, when $f(\x, \cdot)$ is strongly convex, we take $\mu_{k}=0$ for all $k$ and then $F(\x, \y^*(\x))$ is the total UL objective. 
Generally, the first term $\Phi_{\mu}(\x)=F(\x,\y^*_{\mu}(\x))$ in the expression of $V_k$ quantifies the (approximate) overall UL objective functions, the second term $\|\y_k-\y^*_{\mu_k}(\x_k)\|^2$ characterizes LL solution errors, and the third term $\|\v_k-\v^*_{\mu_k}(\x_k)\|^2$ delineates the multiplier errors. 

The proofs of Theorem~\ref{thmconvex1}-\ref{thm_Strong_convex} contain three major steps, which include: (1) upper-bounding the descent of the total UL objective; (2) controlling both LL solution and multiplier errors in the descent of the total UL objective; and (3) combining all results in the previous steps and choosing suitable coefficients $\{ a_k, b_k, c_k \}_{k=1}^\infty$ of Lyapunov function to prove the convergence guarantee. More detailed steps can be found as follows.

{\bf Step 1: upper-bounding the descent of the total UL objective.}

We first analyze the descent of the approximate total UL objective as below, which does not require either LLSC or any assumption on the boundedness  of $\nabla_{\y} F(\x,\y)$.
\begin{equation}\label{step1}
	\begin{aligned}
	&\Phi_{\mu_{k+1}}(\x_{k+1})-\Phi_{\mu_{k}}(\x_k)  \\
	\leq & -\frac{\alpha_k}{2}\|\nabla \Phi_{\mu_{k}}(\x_k)\|^2
	-\frac{1}{2}\left(\frac{1}{\alpha_k}
	-\frac{C_{\Phi1}\|\v^*_{\mu_k}(\x_k)\|+C_{\Phi2}}{\sigma_{\psi_{\mu_k}}^2}\right)\|\x_{k+1}-\x_k\|^2 
	+\alpha_k L_{\psi_{\y1}}^2\|\v_{k}-\v^*_{\mu_{k}}(\x_k)\|^2 \\  
	&+\alpha_k \left( L_{\psi_{\x\y2}} \left\| \v^*_{\mu_{k}}(\x_k)\right\| + L_{F_{\x2}} \right)^2 
	\|\y_{k}-\y^*_{\mu_k}(\x_k)\|^2
	+\frac{2\big\|\nabla_{\y} F(\x_{k+1},\y^*_{\mu_{k+1}}(\x_{k+1}))\big\|^2}{\sigma_{F}} \left(\frac{\mu_k-\mu_{k+1}}{\mu_k}\right),
\end{aligned}
\end{equation}
where both LL solution and multiplier errors appear on the right hand side. Here $C_{\Phi1}$ and $C_{\Phi2}$ are some explicit positive constants given in Lemma~\ref{lem12}. Importantly, compared to the previous results using both LLSC and GBA, the coefficients of $\| \x_{k+1} - \x_{k} \|^2$ and $\|\y_{k}-\y^*_{\mu_k}(\x_k)\|^2$ depend on $\| \v^*_{\mu_{k}}(\x_k) \| $ since Lipschitz continuity of the hypergradient and its surrogate could not be guaranteed without any assumption on the boundedness  of $\nabla_{\y} F(\x,\y)$ and its variant. To address this issue, we characterize a weaker smoothness of $\Phi_{\mu}(\cdot)$ in Lemma \ref{lem12}. Note also that there is an extra term in Inequality \eqref{step1}  involving the descent of the averaging parameter. 

{\bf Step 2: controlling both LL solution and multiplier errors.}

Noticing the above key points, we then estimate carefully LL solution error as below,
\begin{equation}\label{step2-1}
\begin{aligned}
	&\| \y_{k+1}-\y^*_{\mu_{k+1}}(\x_{k+1}) \|^2
	- \| \y_k-\y^*_{\mu_k}(\x_k) \|^2 
	\leq 
	-\frac{1}{2} \beta_k \sigma_{\psi_{\mu_{k}}}
	\| \y_k-\y^*_{\mu_k}(\x_k) \|^2
	+\frac{6L_{\psi_{\y1}}^2}{\beta_{k}\sigma_{\psi_{\mu_k}}^3}
	\|\x_k-\x_{k+1}\|^2 \\
	&\qquad\qquad\qquad\qquad\qquad\qquad\qquad\qquad\qquad\qquad
	+\frac{24 \| \nabla_{\y} F(\x_{k+1},\y^*_{\mu_{k+1}}(\x_{k+1})) \|^2}{\sigma_F^2 \beta_{k} \sigma_{\psi_{\mu_{k}}} }
	\left(\frac{\mu_{k}-\mu_{k+1}}{\mu_k}\right)^2,
\end{aligned}
\end{equation}
and the multiplier error is bound by 
\begin{equation}\label{step2-1}
	\begin{aligned}
		& \| \v_{k+1}-\v^*_{\mu_{k+1}}(\x_{k+1}) \|^2
		- \| \v_k-\v^*_{\mu_{k}}(\x_k) \|^2 \\
		\leq &
		-\frac{1}{2}\eta_k\sigma_{\psi_{\mu_{k}}}
		\| \v_k-\v^*_{\mu_{k}}(\x_k) \|^2
		+\frac{6\left(L_{\v1} \|\v^*_{\mu_{k}}(\x_{k})\|+L_{\v2}\right)^2}{\eta_{k}\sigma_{\psi_{\mu_k}}^5}
		\|\x_k-\x_{k+1}\|^2  \\
		&+\frac{3\eta_{k}}{ \sigma_{\psi_{\mu_k}}}   
		\left(L_{\psi_{\y\y2}} \|\v^*_{\mu_k}(\x_k)\|
		+ L_{F_{\y2}}\right)^2
		\| \y_{k}-\y^*_{\mu_k}(\x_k) \|^2  \\
		&+\frac{6\left(C_{\v1} \|\v^*_{\mu_{k+1}}(\x_{k+1})\|
			+C_{\v2}\right)^2  
			\|\nabla_{\y} F(\x_{k+1},\y^*_{\mu_{k+1}}(\x_{k+1}))\|^2 }{\eta_k\sigma_{\psi_{\mu_{k}}}}
		\left(\frac{\mu_{k}-\mu_{k+1}}{\mu_{k}^2}\right)^2.
	\end{aligned}
\end{equation}
Here $C_{\v1}$ and $C_{\v2}$ are some explicit positive constants given in Lemma~\ref{lem13}. 

{\bf Step 3: choosing suitable coefficients such that the descent of Lyapunov function is well controlled.
 }

Combining the estimates in Steps 1 and 2, we can find suitable decreasing coefficients $\{ a_k, b_k, c_k \}_{k=1}^\infty$ and the averaging parameter $\mu_{k}$ such that the descent of Lyapunov function is well controlled by a summable series as follows: 
\begin{equation}
	\begin{aligned}
		V_{k+1}-V_k\leq & -\frac{a_{k+1} \alpha_k}{2}\|\nabla \Phi_{\mu_{k}}(\x_k)\|^2
		-\frac{\widehat{\alpha}_k}{2}\|\x_{k+1}-\x_k\|^2	
		-\frac{\widehat{\beta}_k}{2}\|\y_k-\y^*_{\mu_{k}}(\x_k)\|^2-\frac{\widehat{\eta}_k}{2}\|\v_k-\v^*_{\mu_{k}}(\x_k)\|^2 
		\\
		&+\frac{2 \big\|\nabla_{\y} F(\x_{k+1},\y^*_{\mu_{k+1}}(\x_{k+1}))\big\|^2 a_{k+1}}{\sigma_{F}} \left(\frac{\mu_k-\mu_{k+1}}{\mu_k}\right) \\
		&+\frac{24 \left\| \nabla_{\y} F(\x_{k+1},\y^*_{\mu_{k+1}}(\x_{k+1})) \right\|^2 b_{k+1}}{\sigma_F^2 \beta_k \sigma_{\psi_{\mu_{k}}}}
		\left(\frac{\mu_{k}-\mu_{k+1}}{\mu_k}\right)^2\\
		&+\frac{6\left(C_{\v1} \|\v^*_{\mu_{k+1}}(\x_{k+1})\|
			+C_{\v2}\right)^2  
			\|\nabla_{\y} F(\x_{k+1},\y^*_{\mu_{k+1}}(\x_{k+1}))\|^2 c_{k+1}}{\eta_k\sigma_{\psi_{\mu_{k}}}}
		\left(\frac{\mu_{k}-\mu_{k+1}}{\mu_{k}^2}\right)^2,
	\end{aligned}
\end{equation}
where the coefficients $\{\widehat{\alpha}_k, \widehat{\beta}_k, \widehat{\eta}_k\}$ explicitly given in \eqref{keycoofficients0} will be positive after choosing proper decreasing coefficients $\{ a_k, b_k, c_k \}_{k=1}^\infty$ and the averaging parameter $\mu_{k}$. Then telescoping leads to $V_K=O(1)$ as $K\rightarrow\infty$ and then 
\begin{align*}
	\sum_{k=K}^{2K+1}
	\frac{a_{k+1} \alpha_k}{2}\|\nabla \Phi_{\mu_{k}}(\x_k)\|^2
	+\frac{\widehat{\alpha}_k}{2}\|\x_{k+1}-\x_k\|^2	
	+\frac{\widehat{\beta}_k}{2}\|\y_k-\y^*_{\mu_{k}}(\x_k)\|^2
	+\frac{\widehat{\eta}_k}{2}\|\v_k-\v^*_{\mu_{k}}(\x_k)\|^2 =O(1).
\end{align*}
The desired convergence guarantee follows by analyzing the decay or grow rates of $a_{k+1} \alpha_k$, $\widehat{\alpha}_k, \widehat{\beta}_k, \widehat{\eta}_k$. We refer to next section for detailed proofs.

Note that Lyapunov function argument for BLO algorithms has been used in establishing the convergence rates of ALSET \cite{chen2021closing}, STABLE \cite{chen2022single}, FSLA \cite{li2022fully}, SOBA and SABA \cite{dagreou2022framework}, and so on. In contrast to the above mentioned works, without GBA and LLSC, Lyapunov function with flexible coefficients must be considered.  

\section{Detailed Proofs}

This section is devoted to proofs of Theorems~\ref{thmconvex1}-\ref{thm_Strong_convex} in Section~\ref{Convergence Analysis}. 

\subsection{Fundamental descent lemmas}\label{dlemma-section}

We first analyze the descent of the approximate overall UL objective $\Phi_{\mu_{k}}(\x_k)=F(\x_k,\y^*_{\mu_k}(\x_k))$  in the next lemma, which is the main result in Step 1 of proof sketch.
\begin{lemma}\label{dlemma1}
	 Suppose Assumptions \ref{Assump0}, and either Assumption \ref{assump_convex0} or Assumption \ref{assump_convex20} holds. Let $\mu_{k+1}\leq\mu_{k}\leq\frac{1}{2}$ for all $k$. Then the sequence of $\x_k,\y_k,\v_k$ generated by sl-BAMM satisfies
	 \begin{equation}\label{estimate1}
	\begin{aligned}
		&\Phi_{\mu_{k+1}}(\x_{k+1})-\Phi_{\mu_{k}}(\x_k)  \\
		\leq & -\frac{\alpha_k}{2}\|\nabla \Phi_{\mu_{k}}(\x_k)\|^2
		-\frac{1}{2}\left(\frac{1}{\alpha_k}
		-\frac{C_{\Phi1}\|\v^*_{\mu_k}(\x_k)\|+C_{\Phi2}}{\sigma_{\psi_{\mu_k}}^2}\right)\|\x_{k+1}-\x_k\|^2 
		+\alpha_k L_{\psi_{\y1}}^2\|\v_{k}-\v^*_{\mu_{k}}(\x_k)\|^2 \\  
		&+\alpha_k \left( L_{\psi_{\x\y2}} \left\| \v^*_{\mu_{k}}(\x_k)\right\| + L_{F_{\x2}} \right)^2 
		\|\y_{k}-\y^*_{\mu_k}(\x_k)\|^2
		+\frac{2\big\|\nabla_{\y} F(\x_{k+1},\y^*_{\mu_{k+1}}(\x_{k+1}))\big\|^2}{\sigma_{F}} \left(\frac{\mu_k-\mu_{k+1}}{\mu_k}\right), 
	\end{aligned}
\end{equation}
where both $C_{\Phi1}$ and $C_{\Phi2}$ are constants given by 
\begin{align*}
	C_{\Phi1}:=&L_{\psi_{\y1}}^2 L_{\psi_{\y\y2}} +L_{\psi_{\y1}}\left(L_{\psi_{\y\y1}}+L_{\psi_{\x\y2}}\right)\sigma_{\psi_\mu} 
	+ L_{\psi_{\x\y1}} \sigma_{\psi_\mu} ^2,\\
	C_{\Phi2}:=&L_{\psi_{\y1}}^2  L_{F_{\y2}} 
	+  L_{\psi_{\y1}} \left(L_{F_{\y1}}+L_{F_{\x2}}\right) \sigma_{\psi_\mu} + L_{F_{\x1}} \sigma_{\psi_\mu}^2.
\end{align*}
	In particular, if LL objective $f(\x,\cdot)$ is $\sigma_{f}$-strongly convex, i.e., Assumption \ref{assump_convex20} holds, we take $\mu_{k}=0$ for all $k$ and then the last term in Estimate \eqref{estimate1} is redundant.
\end{lemma}

Since the descent of the approximate overall UL objective depends on both LL solution and multiplier errors, we next analyze these errors, yielding the main result in Step 2 of proof sketch.
\begin{lemma}\label{dlemma2}
	Suppose Assumptions in Lemma~\ref{dlemma1} hold.  If we choose 
	\begin{equation*}
			\beta_k\leq\frac{2}{\sigma_{\psi_{\mu_{k}}}+L_{\psi_{\mu_{k}}}},\quad 
			\eta_k\leq \frac{1}{L_{\psi_{\mu_k}}}.
		\end{equation*}
	Then the sequence of $\x_k,\y_k,\v_k$ generated by  sl-BAMM satisfies
	\begin{equation}\label{estimate2}
	\begin{aligned}
			\| \y_{k+1}-\y^*_{\mu_{k+1}}(\x_{k+1}) \|^2
			- \| \y_k-\y^*_{\mu_k}(\x_k) \|^2
			\leq &
			-\frac{1}{2} \beta_k \sigma_{\psi_{\mu_{k}}}
			\| \y_k-\y^*_{\mu_k}(\x_k) \|^2
			+\frac{6L_{\psi_{\y1}}^2}{\beta_{k}\sigma_{\psi_{\mu_k}}^3}
			\|\x_k-\x_{k+1}\|^2\\
			&+\frac{24 \| \nabla_{\y} F(\x_{k+1},\y^*_{\mu_{k+1}}(\x_{k+1})) \|^2}{\sigma_F^2 \beta_{k} \sigma_{\psi_{\mu_{k}}} }
			\left(\frac{\mu_{k}-\mu_{k+1}}{\mu_k}\right)^2,
		\end{aligned}
	\end{equation}
	and 
	\begin{equation}\label{estimate3}
	\begin{aligned}
			& \| \v_{k+1}-\v^*_{\mu_{k+1}}(\x_{k+1}) \|^2
			- \| \v_k-\v^*_{\mu_{k}}(\x_k) \|^2 \\
			\leq &
			-\frac{1}{2}\eta_k\sigma_{\psi_{\mu_{k}}}
			\| \v_k-\v^*_{\mu_{k}}(\x_k) \|^2
			+\frac{6\left(L_{\v1} \|\v^*_{\mu_{k}}(\x_{k})\|+L_{\v2}\right)^2}{\eta_{k}\sigma_{\psi_{\mu_k}}^5}
			\|\x_k-\x_{k+1}\|^2 \\
			&+\frac{3\eta_{k}}{ \sigma_{\psi_{\mu_k}}}   
			\left(L_{\psi_{\y\y2}} \|\v^*_{\mu_k}(\x_k)\|
			+ L_{F_{\y2}}\right)^2
			\| \y_{k}-\y^*_{\mu_k}(\x_k) \|^2\\
			&+\frac{6\left(C_{\v1} \|\v^*_{\mu_{k+1}}(\x_{k+1})\|
					+C_{\v2}\right)^2  
					\|\nabla_{\y} F(\x_{k+1},\y^*_{\mu_{k+1}}(\x_{k+1}))\|^2 }{\eta_k\sigma_{\psi_{\mu_{k}}}}
			\left(\frac{\mu_{k}-\mu_{k+1}}{\mu_{k}^2}\right)^2,
		\end{aligned}
	\end{equation}
where both $C_{\v1}$ and $C_{\v2}$ are constants given by 
\begin{align*}
	C_{\v1}:=2\left(L_{F_{\y\y2}}+L_{f_{\y\y2}}\right)/\sigma_{F}^2,
	\quad
	C_{\v2}:=2\left(2L_{F_{\y2}}+L_{f_{\y2}}\right)/\sigma_{F}^2.
\end{align*}
In particular, if LL objective $f(\x,\cdot)$ is $\sigma_{f}$-strongly convex, we take $\mu_{k}=0$ for all $k$ and then the last terms in Estimates \eqref{estimate2}-\eqref{estimate3} are both redundant.
\end{lemma}

Combining the estimates in Steps 1 and 2, for any decreasing positive coefficients $\{ a_k, b_k, c_k \}_{k=1}^\infty$, we get the following descent of Lyapunov function, which is the key step to improving the existing result.  
\begin{lemma}\label{dlemma3}
	Suppose Assumptions in Lemma \ref{dlemma2} hold. 
	Let $\{ a_k, b_k, c_k \}_{k=1}^\infty$ be a sequence of decreasing positive constants, then the following descent of Lyapunov function holds: 
	\begin{equation}\label{estimate4}
			\begin{aligned}
					V_{k+1}-V_k\leq & -\frac{a_{k+1} \alpha_k}{2}\|\nabla \Phi_{\mu_{k}}(\x_k)\|^2-\frac{\widehat{\alpha}_k}{2}\|\x_{k+1}-\x_k\|^2	
					-\frac{\widehat{\beta}_k}{2}\|\y_k-\y^*_{\mu_{k}}(\x_k)\|^2-\frac{\widehat{\eta}_k}{2}\|\v_k-\v^*_{\mu_{k}}(\x_k)\|^2 
					\\
					&+\frac{2 \big\|\nabla_{\y} F(\x_{k+1},\y^*_{\mu_{k+1}}(\x_{k+1}))\big\|^2 a_{k+1}}{\sigma_{F}} \left(\frac{\mu_k-\mu_{k+1}}{\mu_k}\right) 
					\\
					&+\frac{24 \left\| \nabla_{\y} F(\x_{k+1},\y^*_{\mu_{k+1}}(\x_{k+1})) \right\|^2 b_{k+1}}{\sigma_F^2 \beta_k \sigma_{\psi_{\mu_{k}}}}
					\left(\frac{\mu_{k}-\mu_{k+1}}{\mu_k}\right)^2\\
					&+\frac{6\left(C_{\v1} \|\v^*_{\mu_{k+1}}(\x_{k+1})\|
							+C_{\v2}\right)^2  
							\|\nabla_{\y} F(\x_{k+1},\y^*_{\mu_{k+1}}(\x_{k+1}))\|^2 c_{k+1}}{\eta_k\sigma_{\psi_{\mu_{k}}}}
					\left(\frac{\mu_{k}-\mu_{k+1}}{\mu_{k}^2}\right)^2,
				\end{aligned}
		\end{equation}
	where the coefficients $\{ \widehat{\alpha}_k, \widehat{\beta}_k, \widehat{\eta}_k \}$ are given as below.
	\begin{equation}\label{keycoofficients0}
		\begin{aligned}
			\widehat{\alpha}_k
			:=&\frac{a_{k+1}}{\alpha_k}
			-\frac{ a_{k+1} \left(C_{\Phi1}\|\v^*_{\mu_k}(\x_k)\| + C_{\Phi2}\right) }{\sigma_{\psi_{\mu_{k}}}^2}
			-\frac{12L_{\psi_{\y1}}^2 b_{k+1} }{\beta_{k}\sigma_{\psi_{\mu_k}}^3}
			-\frac{ 12 \left(L_{\v1} \|\v^*_{\mu_{k}}(\x_{k})\|+L_{\v2}\right)^2 c_{k+1} }{ \eta_{k}\sigma_{\psi_{\mu_k}}^5 } ,\\
			\widehat{\beta}_k
			:=&
			b_{k+1}\beta_k \sigma_{\psi_{\mu_{k}}}
			-2 a_{k+1}\alpha_k \left( L_{\psi_{\x\y2}} \| \v^*_{\mu_k}(\x_k) \| + L_{F_{\x2}} \right)^2 
			-6 
			\left(L_{\psi_{\y\y2}} \| \v^*_{\mu_k}(\x_k) \| 
			+ L_{F_{\y2}}\right)^2
			\frac{ c_{k+1}\eta_{k}}{ \sigma_{\psi_{\mu_k}}} , \\
			\widehat{\eta}_k
			:=&
			c_{k+1}\eta_k\sigma_{\psi_{\mu_{k}}}
			-2 a_{k+1} \alpha_k  L_{\psi_{\y1}}^2.
		\end{aligned}
	\end{equation}
	Here $L_{\v1}:=L_{\psi_{\y\y2}} L_{\psi_{\y1}}+L_{\psi_{\y\y1}} \sigma_{\psi_{\mu}} $ and $L_{\v2}:=L_{F_{\y2}} L_{\psi_{\y1}}
	+L_{F_{\y1}} \sigma_{\psi_{\mu}} $. In particular, if LL objective $f(\x,\cdot)$ is $\sigma_{f}$-strongly convex, we take $\mu_{k}=0$ for all $k$ and then the terms involving $\mu_{k}-\mu_{k+1}$ in Estimate \eqref{estimate4} are all redundant. 
\end{lemma}

\subsection{A general setting of step sizes}

To simplify and unify the proofs of Theorems~\ref{thmconvex1}-\ref{thm_Strong_convex}, we introduce our framework in which all of the step size strategies in Theorems~\ref{thmconvex1}-\ref{thm_Strong_convex} have a unified form as below,
\begin{equation}\label{unified-stepsizes}
	\beta_{k}\in \left[ \bar{\beta}(k+1)^{-\tau/2}, 1/ \left(L_{F_{\y2}} + L_{f_{\y2}}\right) \right]=:I_\beta,
	\quad
	\eta_k= c_\eta \frac{b_{k+1}}{c_{k+1}}\beta_k \sigma_{\psi_{\mu_{k}}}^2 ,
	\quad
	\alpha_k=c_\alpha \frac{a_{k+1}}{c_{k+1}} \eta_k \sigma_{\psi_{\mu_{k}}}^5,
\end{equation}
where $a_{k+1}, b_{k+1}, c_{k+1}$ are the coefficients of Lyapunov function \eqref{Lyapunov}, and the positive constants $c_\eta, c_\alpha$ may depend on $k$. Here and in what follows, we drop their dependency in $k$ for clarity. Actually, for any fixed step size strategy $(\beta_k,\eta_k,\alpha_k)$ with $\beta_k\in I_\beta$, after defining
\begin{equation}\label{ceta-calpha}
	c_\eta=\left(\frac{c_{k+1}}{b_{k+1}}\right) \left(\frac{\eta_{k}}{\beta_{k}}\right) \sigma_{\psi_{\mu_{k}}}^{-2} ,
	\quad
	c_\alpha=\left(\frac{c_{k+1}}{a_{k+1}}\right) \left(\frac{\alpha_{k}}{\eta_{k}}\right) \sigma_{\psi_{\mu_{k}}}^{-5},
\end{equation}
the step size strategy $(\beta_k,\eta_k,\alpha_k)$ has the form \eqref{unified-stepsizes}. 

Thanks to the setting \eqref{unified-stepsizes}, the key coefficients $\{ \widehat{\alpha}_k, \widehat{\beta}_k, \widehat{\eta}_k \}$ in Estimate \eqref{estimate4} have the following simplified form: 
\begin{align}\label{simpleform1}
	\widehat{\alpha}_k
	=&\left(1- c_\alpha \left( 12 \left(L_{\v1} \|\v^*_{\mu_{k}}(\x_{k})\|+L_{\v2}\right)^2 
	+\frac{12L_{\psi_{\y1}}^2 b_{k+1} \alpha_k }{c_\alpha a_{k+1}\beta_{k}\sigma_{\psi_{\mu_k}}^3} 
	+\frac{\alpha_k}{c_\alpha \sigma_{\psi_{\mu_k}}^2}
	\left(C_{\Phi1}\|\v^*_{\mu_k}(\x_k)\| + C_{\Phi2}\right)\right)
	\right)\frac{a_{k+1}}{\alpha_k}, \\
	\widehat{\beta}_k
	= &
	\left(1-c_\eta
	\left(
	6\left(L_{\psi_{\y\y2}} \| \v^*_{\mu_k}(\x_k) \| 
	+ L_{F_{\y2}}\right)^2
	+ \frac{2 a_{k+1}\alpha_k}{c_\eta b_{k+1}\beta_k \sigma_{\psi_{\mu_{k}}}}\left( L_{\psi_{\x\y2}}^2 \| \v^*_{\mu_k}(\x_k) \| + L_{F_{\x2}} \right)^2
	\right)\right)
	b_{k+1}\beta_k \sigma_{\psi_{\mu_{k}}}, \\
	\widehat{\eta}_k
	=&\left( 1- 2 L_{\psi_{\y1}}^2 \frac{a_{k+1} \alpha_k   }{c_{k+1} \eta_k\sigma_{\psi_{\mu_{k}}}} \right)c_{k+1}\eta_k\sigma_{\psi_{\mu_{k}}}.
\end{align}
For any fixed step size strategy $(\beta_k,\eta_k,\alpha_k)$, to control the descent of Lyapunov function very well, it is helpful to find suitable coefficients $a_{k+1}, b_{k+1}, c_{k+1}$ such that 
\begin{equation}\label{keyestimate1}
	\widehat{\alpha}_k
	\geq \frac{a_{k+1}}{2\alpha_k},
	\quad
	\widehat{\beta}_k
	\geq 
	\frac{1}{2} b_{k+1}\beta_k \sigma_{\psi_{\mu_{k}}},
	\quad
	\widehat{\eta}_k
	\geq \frac{1}{2} c_{k+1}\eta_k\sigma_{\psi_{\mu_{k}}},
\end{equation}
and all of the sequences involving $\mu_k-\mu_{k+1}$ in Estimate \eqref{estimate4} (denoted by $A_\mu^k, B_\mu^k, C_\mu^k$ respectively ) are summable. If this is well done, from Estimate \eqref{estimate4}, one first get 
\begin{equation}\label{estimate5}
	\begin{aligned}
		V_{k+1}-V_k\leq A_\mu^k + B_\mu^k + C_\mu^k.
	\end{aligned}
\end{equation}
By the definition \eqref{Lyapunov}, we get $V_k\geq0$ for all $k$ and then
\begin{equation}
	V_{K}
	\leq V_{K_0} +\sum_{k=1}^\infty\left(A_\mu^k + B_\mu^k + C_\mu^k\right) 
	\implies 
	V_K=O(1)\ \mathrm{as}\ K\rightarrow\infty.
\end{equation}
Further, by using Estimate \eqref{estimate4} again, one can get 
\begin{equation}
	\sum_{k=K}^{2K+1}
	\left(\frac{a_{k+1} \alpha_k}{2}\|\nabla \Phi_{\mu_{k}}(\x_k)\|^2
	+\frac{\widehat{\alpha}_k}{2}\|\x_{k+1}-\x_k\|^2	
	+\frac{\widehat{\beta}_k}{2}\|\y_k-\y^*_{\mu_{k}}(\x_k)\|^2
	+\frac{\widehat{\eta}_k}{2}\|\v_k-\v^*_{\mu_{k}}(\x_k)\|^2
	\right)
	=O(1). 
\end{equation}
Then Estimate \eqref{keyestimate1} yields 
\begin{equation}\label{keyestimate2}
	\begin{aligned}
	\sum_{k=K}^{2K+1}
	\Big( 
	2a_{k+1} \alpha_k\|\nabla \Phi_{\mu_{k}}(\x_k)\|^2
	&+\frac{a_{k+1}}{\alpha_k} \|\x_{k+1}-\x_k\|^2	\\
	&+ b_{k+1}\beta_k \sigma_{\psi_{\mu_{k}}}
	\|\y_k-\y^*_{\mu_{k}}(\x_k)\|^2
	+c_{k+1}\eta_k\sigma_{\psi_{\mu_{k}}}
	\|\v_k-\v^*_{\mu_{k}}(\x_k)\|^2
	\Big)
	=O(1). 
\end{aligned}
\end{equation}
Now, by lower-bounding the coefficients in Estimate \eqref{keyestimate2}, according to different step size strategies and the coefficients of Lyapunov function, we can get the estimates on 
\begin{equation*}
	\min_{0\leq k\leq K} \| \nabla \Phi_{\mu_{k}}(\x_k) \|^2
	+ \|\x_{k+1}-\x_k\|^2	
	+  \|\y_k-\y_{\mu_k}^*(\x_k)\|^2 
	+ \|\v_k-\v_{\mu_k}^*(\x_k)\|^2,
\end{equation*}
by choosing different upper bounds on $p$ and $\tau$. Finally, by the optimality of $\y^*_{\mu_k}(\x_k)$ and the definition of $\v^*_{\mu_{k}}(\x_k)$, there exists a positive constant $C$ independent of $k$ such that 
\begin{align}
	\mathrm{KKT}(\x_{k}, \y_{k}, \v_{k})
	\leq 
	C\left(\|\nabla \Phi_{\mu_{k}}(x_{k})\|^2 + \|\y_{k}-\y^*_{\mu_{k}}(\x_{k})\|^2 + \|\v_{k}-\v^*_{\mu_{k}}(\x_{k})\|^2+ \mu_{k}^2\right),
\end{align}
which yields the convergence rates measured by KKT residual. 

\subsection{Proofs of Theorems~\ref{thmconvex1}-\ref{thm_Strong_convex}
}

For the convenience of the reader, we first give an overview of the step sizes strategies in the main theorems and of the well-chosen coefficients of Lyapunov function in Table \ref{table-stepsizes}.
\begin{table}[h!]
	\centering 
	\makeatletter\def\@captype{table}\makeatother \caption{Summary of the step size strategies and the well-chosen Lyapunov coefficients.} 
	\vspace{0.2cm}
	\setlength{\tabcolsep}{1.4mm}{
		\begin{tabular}{ccccccccc}
			\toprule
			\multirow{2}{*}{Strategy}&\multicolumn{3}{c}{Step sizes}&&&\multicolumn{3}{c}{Coefficients of Lyapunov function} \\
			&\multirow{1}{*}{$\beta_{k}$}&\multirow{1}{*}{$\eta_{k}$}&\multirow{1}{*}{$\alpha_{k}$}&&& $a_{k+1}$&$b_{k+1}$ & $c_{k+1}$ \\
			\midrule
			Unified form
			& $\beta_{k} $ & $c_\eta \frac{b_{k+1}}{c_{k+1}}\beta_k \sigma_{\psi_{\mu_{k}}}^2 $&$ c_\alpha \frac{a_{k+1}}{c_{k+1}} \eta_k \sigma_{\psi_{\mu_{k}}}^5 $&&&  $a_{k+1}$ &  $ b_{k+1} $   & $ c_{k+1} $   \\
			\midrule
			S1 (Thm 3.4) &$\beta_{k}\in I_\beta$& $ (k+1)^{-\tau/2} \beta_{k} \mu_{k}^2 $ &$ (k+1)^{-3\tau/2} \beta_{k} \mu_{k}^7 $&&&  $(k+1)^{-\tau} \sigma_{\psi_{\mu_{k}}}^2 $ &  $(k+1)^{-\tau/2} \sigma_{\psi_{\mu_{k}}}^2 $   & $(k+1)^{-\tau} \sigma_{\psi_{\mu_{k}}}^6 $  \\
			\midrule
			S2 (Thm 3.5) &$\beta_{k}\in I_\beta$& $ (k+1)^{-\tau/2} \beta_{k} \mu_{k} $ & $ (k+1)^{-3\tau/2} \beta_{k} \mu_{k}^5 $ &&&  $(k+1)^{-\tau} $ &  $(k+1)^{-\tau/2}  $   & $(k+1)^{-\tau} \sigma_{\psi_{\mu_{k}}}^3 $  \\
			\midrule
			S3 (Thm 3.6)&$\beta_{k}\in I_\beta$& $ (k+1)^{-\tau/2} \beta_{k} $ & $ (k+1)^{-3\tau/2} \beta_{k} \mu_{k}^3 $ &&&  $(k+1)^{-\tau}  $ &  $(k+1)^{-\tau/2}  $   & $(k+1)^{-\tau} \sigma_{\psi_{\mu_{k}}}^2 $   \\
			\midrule
			Thm 3.10 &$\beta_{k}\in I_\beta$& $ \bar{\eta}(k+1)^{-\tau/2} \beta_{k} $ & $ \bar{\alpha}(k+1)^{-\tau} \beta_{k} $ &&&  $1 $ &  $ 1$   & $ 1 $ \\ 
			\bottomrule		
		\end{tabular}
	}
\label{table-stepsizes}
	\vspace{-0.4cm} 
\end{table}

Using the simplified form \eqref{simpleform1} of $ \widehat{\alpha}_k, \widehat{\beta}_k, \widehat{\eta}_k $, note that $\frac{a_{k+1} \alpha_k}{c_{k+1} \eta_k\sigma_{\psi_{\mu_{k}}}}
=\frac{a_{k+1} \alpha_k}{c_{\eta} b_{k+1}\beta_k \sigma_{\psi_{\mu_{k}}}^3}
= c_\alpha\left(\frac{ a_{k+1} }{c_{k+1} } \right)^2 \sigma_{\psi_{\mu_{k}}}^4$ and 
\begin{align*}
	\frac{ b_{k+1} \alpha_k }{c_\alpha a_{k+1}\beta_{k}\sigma_{\psi_{\mu_k}}^3} 
	= c_\eta \left(\frac{ b_{k+1} }{c_{k+1} } \right)^2 \sigma_{\psi_{\mu_{k}}}^4,
	\quad 
	\frac{\alpha_k}{c_\alpha \sigma_{\psi_{\mu_k}}^2}
	= \left(\frac{a_{k+1}}{c_{k+1}}\right) \eta_k \sigma_{\psi_{\mu_{k}}}^3,
	\quad
	\frac{a_{k+1}\alpha_k}{c_\eta b_{k+1}\beta_k \sigma_{\psi_{\mu_{k}}}} 
	=c_\alpha \left(\frac{a_{k+1}}{c_{k+1}}\right)^2  \sigma_{\psi_{\mu_{k}}}^6,
\end{align*}
we can prove that Estimate \eqref{keyestimate1} holds in the setting of Theorems~\ref{thmconvex1}-\ref{thm_Strong_convex}, with the help of the estimates in Table \ref{table-estimates1}.
\begin{table}[h!]
	\centering 
	\makeatletter\def\@captype{table}\makeatother \caption{Estimates for the proof of Estimate \eqref{keyestimate1}.} 
	\vspace{0.2cm}
	\setlength{\tabcolsep}{1.4mm}{
		\begin{tabular}{ccccccccc}
			\toprule
			Strategy&{$c_\alpha$}&\multirow{1}{*}{$c_\eta$}&$\| \v_{\mu_{k}}^*(\x_k) \|$& $\frac{ b_{k+1} \alpha_k }{c_\alpha a_{k+1}\beta_{k}\sigma_{\psi_{\mu_k}}^3} $&$\frac{\alpha_k}{c_\alpha \sigma_{\psi_{\mu_k}}^2}$ & $ \frac{a_{k+1} \alpha_k}{c_{k+1} \eta_k\sigma_{\psi_{\mu_{k}}}}$ \\
			\midrule
			S1 (Thm 3.4) &$(k+1)^{-\tau} \mu_{k}^4 $& $(k+1)^{-\tau} \mu_{k}^4 $ &$ \mu_{k}^{-2} $&  $ 1$ &  $(k+1)^{-\tau/2} \beta_k \mu_{k} $   & $(k+1)^{-\tau} $  \\
			\midrule
			S2 (Thm 3.5) &$(k+1)^{-\tau} \mu_{k}^2 $& $(k+1)^{-\tau} \mu_{k}^2 $ & $ \mu_{k}^{-1} $ & $1  $   & $(k+1)^{-\tau/2} \beta_k \mu_{k} $  & $(k+1)^{-\tau}$\\
			\midrule
			S3 (Thm 3.6)&$(k+1)^{-\tau} $& $(k+1)^{-\tau} $ & $ 1 $ & $1  $   & $(k+1)^{-\tau/2} \beta_k \mu_{k} $  & $(k+1)^{-\tau/2}  $     \\
			\midrule
			Thm 3.10 &$(k+1)^{-\tau/2} $& $(k+1)^{-\tau/2} $ & $ 1 $ & $ (k+1)^{-\tau/2} $   & $ (k+1)^{-\tau/2} $   & $ (k+1)^{-\tau/2} $      \\
			\bottomrule		
		\end{tabular}
	}
	\label{table-estimates1}
	\vspace{-0.4cm} 
\end{table}

Next we show that all of the sequences involving $\mu_k-\mu_{k+1}$ in Estimate \eqref{estimate4} (denoted by $A_\mu^k, B_\mu^k, C_\mu^k$ respectively ) are summable. Recall that 
\begin{align*}
	A_\mu^k:=&\frac{2 \big\|\nabla_{\y} F(\x_{k+1},\y^*_{\mu_{k+1}}(\x_{k+1}))\big\|^2 a_{k+1}}{\sigma_{F}} \left(\frac{\mu_k-\mu_{k+1}}{\mu_k}\right),\\
	 B_\mu^k:=&\frac{24 \left\| \nabla_{\y} F(\x_{k+1},\y^*_{\mu_{k+1}}(\x_{k+1})) \right\|^2 b_{k+1}}{\sigma_F^2 \beta_k \sigma_{\psi_{\mu_{k}}}}
	 \left(\frac{\mu_{k}-\mu_{k+1}}{\mu_k}\right)^2,\\
	 C_\mu^k:=&\frac{6\left(C_{\v1} \|\v^*_{\mu_{k+1}}(\x_{k+1})\|
		 	+C_{\v2}\right)^2  
		 	\|\nabla_{\y} F(\x_{k+1},\y^*_{\mu_{k+1}}(\x_{k+1}))\|^2 c_{k+1}}{\eta_k\sigma_{\psi_{\mu_{k}}}}
	 \left(\frac{\mu_{k}-\mu_{k+1}}{\mu_{k}^2}\right)^2.
\end{align*}
Using the step sizes strategies and the well-chosen coefficients of Lyapunov function in Table \ref{table-stepsizes}, we can get the explicit estimates of $A_\mu^k, B_\mu^k, C_\mu^k$ and the ones for the coefficients in Estimate \eqref{keyestimate2}, which is summarized in Table \ref{table-estimates2}.
\begin{table}[h!]
	\centering 
	\makeatletter\def\@captype{table}\makeatother \caption{Important estimates used in the proofs of Theorems~\ref{thmconvex1}-\ref{thm_Strong_convex}.} 
	\vspace{0.2cm}
	\setlength{\tabcolsep}{1.4mm}{
		\begin{tabular}{cccccccccc}
			\toprule
			Strategy&$A_\mu^k$&$B_\mu^k$&$C_\mu^k$& $a_{k+1} \alpha_{k}$ & $\frac{ 1 }{ \alpha_{k}^2} $&$\frac{ b_{k+1}\beta_k \sigma_{\psi_{\mu_{k}}}}{ a_{k+1} \alpha_{k} }$ & $ \frac{c_{k+1} \eta_k \sigma_{\psi_{\mu_{k}}}}{a_{k+1} \alpha_{k}}$ & $\frac{a_{k+1}}{\alpha_{k}}$ \\
			\midrule
			S1 (Thm 3.4) &$k^{-1-\tau}  $& $k^{-2+p} $ & $k^{-2+5p} $ &  $ k^{-9p-3\tau}$ &  $k^{ 14p+3\tau}$   & $k^{6p+2\tau} $ & $k^{\tau}$ & $k^{5p+\tau/2}$ \\
			\midrule
			S2 (Thm 3.5) &$k^{-1-\tau} $& $k^{-2+p} $ & $k^{-2+3p} $ & $k^{-5p-3\tau}  $   & $k^{10p+3\tau}  $  & $k^{4p+2\tau}$ & $k^{\tau}$ & $k^{5p+\tau/2}$ \\
			\midrule
			S3 (Thm 3.6)&$k^{-1-\tau} $& $k^{-2-\tau} $ & $k^{-2-\tau} $ & $k^{-3p-3\tau}$   & $k^{4p+2\tau} $  & $k^{\tau} $ &$ k^{\tau} $ & $k^{3p+\tau/2}$   \\
			\midrule
			Thm 3.10 &$0$& $0$ & $ 0 $ & $ k^{-3\tau/2} $   & $ k^{2\tau} $   & $ k^{\tau} $  & $ k^{\tau/2} $ & $k^{\tau}$    \\
			\bottomrule		
		\end{tabular}
	}
	\label{table-estimates2}
	\vspace{-0.4cm} 
\end{table}
Finally, by using the following inequality
\begin{equation}
	\sum_{k=K}^{2K+1} (k+1)^{-s}\geq \int_{K+1}^{2K+2} t^{-s} dt=\frac{2^{1-s}-1}{1-s}(K+1)^{1-s},
\end{equation}
we can get the estimates on 
\begin{equation*}
	\min_{K\leq k\leq 2K+1} \| \nabla \Phi_{\mu_{k}}(\x_k) \|^2
	+ \|\x_{k+1}-\x_k\|^2	
	+  \|\y_k-\y_{\mu_k}^*(\x_k)\|^2 
	+ \|\v_k-\v_{\mu_k}^*(\x_k)\|^2,
\end{equation*}
and then the convergence rate results in Theorems~\ref{thmconvex1}-\ref{thm_Strong_convex} hold.

\section{Proofs of fundamental descent lemmas}
\label{lemmas}

This section is devoted to the proofs of fundamental descent lemmas in Section \ref{dlemma-section}. Note that, unlike existing analyses, we do not require a strong assumption on the boundedness of $\nabla_{\y} F(\x,\y)$. Hence, our theoretical analysis can be applied to a wider variety of applications in deep learning. To handle general BLOs, e.g., BLOs with multiple LL minimizers, we have taken the following standard assumption as in \cite{liu2020generic}.

\begin{assumption}[\bf Assumption \ref{asump-strongconvex}]\label{asump-strongconvexa}
	For any $\x$, the aggregation function $\psi_{\mu}(\x,\cdot)=\mu  F(\x,\cdot)+(1-\mu)f(\x,\cdot)$ is strongly convex for $\mu=0$ or $\mu>0$.
\end{assumption}
Under this assumption, the aggregation function $\psi_{\mu}(\x,\y)=\mu F(\x,\y)+(1-\mu) f(\x,\y)$ is $\sigma_{\psi_{\mu}}$-strongly convex in $\y$ with $\sigma_{\psi_{\mu}}=\mu \sigma_{F} +(1-\mu)\sigma_{f}$. Hence $\psi_{\mu}(\x,\y)$ has a unique minimal point for every $\x$, denoted by $\y^*_{\mu}(\x)$. Furthermore, under the smoothness assumption of $f$ and $F$, $\psi_{\mu}(\x,\y)$ is also $L_{\psi_{\mu}}$-smooth in $\y$ for any fixed $\x$. 

\subsection{Auxiliary Lemmas}


\begin{lemma}\label{lemma*}
	Suppose Assumptions \ref{Assump0}, and either Assumption \ref{assump_convex0} or Assumption \ref{assump_convex20} holds, the following statements hold.
	\begin{enumerate}[(i)]
		\item The function $\y^*_{\mu}(\x)$ is Lipschitz continuous in $\x$, i.e., for all $\x,\x'$, 
		\begin{equation}
			\left\| \y^*_{\mu}(\x)- \y^*_{\mu}(\x')\right\|
			\leq \frac{L_{\psi_{\y1}}}{\sigma_{\psi_\mu}}
			\|\x-\x'\|,
		\end{equation}
		where both $L_{\psi_{\y1}}$ and $\sigma_{\psi_\mu}$ are constants given by 
		\begin{equation*}
			\sigma_{\psi_\mu}=\mu \sigma_F+(1-\mu)\sigma_f,\quad L_{\psi_{\y1}}= \mu L_{F_{\y1}}+(1-\mu)L_{f_{\y1}}.
		\end{equation*}
		
		\item The function $\v^*_{\mu}(\x)$ satisfies 
		\begin{equation}\label{boundv}
			\|\v^*_{\mu}(\x)\|=\left\|\left[\nabla_{\y\y}^2 \psi_{\mu}(\x,\y^*_{\mu}(\x))\right]^{-1}\nabla_{\y} F(\x,\y^*_{\mu}(\x))\right\|
			\leq 
			\frac{\|\nabla_{\y} F(\x,\y^*_{\mu}(\x))\|}{\sigma_{\psi_{\mu}}},
		\end{equation}
		and for all $\x,\x'$, 
		\begin{align}
			\left\| \v^*_{\mu}(\x)- \v^*_{\mu}(\x')\right\|
			\leq \frac{ \left(L_{\v1} \|\v^*_{\mu}(\x)\| +L_{\v2}\right) }{\sigma_{\psi_{\mu}}^2}
			\|\x-\x'\| , 
		\end{align}
		where both $L_{\v1}:=L_{\psi_{\y\y2}} L_{\psi_{\y1}}+L_{\psi_{\y\y1}} \sigma_{\psi_{\mu}} $ and $L_{\v2}:=L_{F_{\y2}} L_{\psi_{\y1}}
		+L_{F_{\y1}} \sigma_{\psi_{\mu}} $ are constants. 
	\end{enumerate}
\end{lemma}
\begin{proof}
	\begin{enumerate}[(i)]
		\item By the optimality of $\y^*_\mu$ to LL problem, we have $\nabla_{\y}\psi_{\mu}(\x,\y^*_{\mu}(\x))=0$. Thus
		\begin{equation*}
			\left(\nabla_{\y}\psi_{\mu}(\x,\y^*_{\mu}(\x))
			-\nabla_{\y}\psi_{\mu}(\x,\y^*_{\mu}(\x'))\right)+\left(\nabla_{\y}\psi_{\mu}(\x,\y^*_{\mu}(\x'))
			-\nabla_{\y}\psi_{\mu}(\x',\y^*_{\mu}(\x'))\right)=0.
		\end{equation*}
		Multiplying the above equation by $\y^*_{\mu}(\x)-\y^*_{\mu}(\x')$, by the $\sigma_{\psi_\mu}$-strong convexity of $\psi_{\mu}(\x,\cdot)$ and $L_{\psi_{y1}}$-smoothness of $\psi_{\mu}(\cdot,\y)$, we have 
		\begin{align*}
			\sigma_{\psi_\mu}\|\y^*_{\mu}(\x)-\y^*_{\mu}(\x')\|^2\leq& \|\nabla_{\y}\psi_{\mu}(\x,\y^*_{\mu}(\x'))
			-\nabla_{\y}\psi_{\mu}(\x',\y^*_{\mu}(\x'))\|\cdot\|\y^*_{\mu}(\x)-\y^*_{\mu}(\x')\|\\
			\leq& L_{\psi_{\y1}} \|\x-\x'\|\cdot \|\y^*_{\mu}(\x)-\y^*_{\mu}(\x')\|.
		\end{align*}
		Then the conclusion follows immediately from the above inequality.
		
		\item First, by the $\sigma_{\psi_\mu}$-strong convexity of $\psi_{\mu}(\x,\cdot)$ and Cauchy-Schwarz inequality, we have 
		\begin{align*}
			\sigma_{\psi_\mu} \|\v^*_{\mu}(\x) \|^2\leq& \left\langle\v^*_{\mu}(\x),  \nabla_{\y\y}^2\psi_{\mu}(\x,\y^*_{\mu}(\x))
			\v^*_{\mu}(\x) \right\rangle
			=\left\langle\v^*_{\mu}(\x),  \nabla_{\y} F(\x,\y^*_{\mu}(\x)) \right\rangle
			\leq \|\v^*_{\mu}(\x) \|\cdot\|\nabla_{\y} F(\x,\y^*_{\mu}(\x))\|, 
		\end{align*}
		which implies Inequality~\eqref{boundv}. 
		
		Second, the definition of $\v^*_{\mu}(\x)$ says that $\left[\nabla_{\y\y}^2 \psi_{\mu}(\x,\y^*_{\mu}(\x))\right]\v^*_{\mu}(\x)=\nabla_{\y} F(\x,\y^*_{\mu}(\x))$. Thus
		\begin{align*}
			&\left[\nabla_{\y\y}^2 \psi_{\mu}(\x',\y^*_{\mu}(\x'))\right]
			\left[\v^*_{\mu}(\x')
			- \v^*_{\mu}(\x)\right]\\
			=&\left[\nabla_{\y} F(\x',\y^*_{\mu}(\x'))
			-\nabla_{\y} F(\x,\y^*_{\mu}(\x))\right]
			+\left[\nabla_{\y\y}^2 \psi_{\mu}(\x,\y^*_{\mu}(\x))
			-\nabla_{\y\y}^2 \psi_{\mu}(\x',\y^*_{\mu}(\x'))\right]\v^*_{\mu}(\x).
		\end{align*}
		By (i), the $\sigma_{\psi_\mu}$-strong convexity of $\psi_{\mu}(\x',\cdot)$ and Lipschitz continuity of $\nabla_{\y}\psi_{\mu}$ and $\nabla_{\y} F$, we have 
		\begin{align*}
			\sigma_{\psi_{\mu}} \|\v^*_{\mu}(\x)- \v^*_{\mu}(\x')\|
			\leq& \left(L_{F_{\y1}}\|\x-\x'\|+L_{F_{\y2}}\|\y^*_{\mu}(\x)-\y^*_{\mu}(\x')\|\right)\\
			&
			+ \|\v^*_{\mu}(\x)\|
			\left(L_{\psi_{\y\y1}}\|\x-\x'\|+L_{\psi_{\y\y2}}\|\y^*_{\mu}(\x)-\y^*_{\mu}(\x')\|\right)\\
			\leq & \frac{\|\x-\x'\|}{\sigma_{\psi_{\mu}}}
			\Big(\|\v^*_{\mu}(\x)\| \left(L_{\psi_{\y\y2}} L_{\psi_{\y1}}+\sigma_{\psi_{\mu}} L_{\psi_{\y\y1}}\right)
			+L_{F_{\y2}} L_{\psi_{\y1}}
			+\sigma_{\psi_{\mu}} L_{F_{\y1}}\Big),
		\end{align*}
		which implies the desired result.
	\end{enumerate}
\end{proof}

Next, we prove a lemma that characterizes the smoothness of $\Phi_\mu(\x)$ in $\x$, without any boundedness assumption of $\nabla_{\y} F(\x, \y)$. 

\begin{lemma}\label{lem12}
	Suppose Assumptions \ref{Assump0}, and either Assumption \ref{assump_convex0} or Assumption \ref{assump_convex20} holds, we have
	\begin{equation}
		\|\nabla \Phi_{\mu}(\x) -\nabla \Phi_{\mu}(\x')\|
		\leq 
		\frac{\left(C_{\Phi1}\|\v^*_{\mu}(\x)\| + C_{\Phi2}\right)}{\sigma_{\psi_{\mu}}^2}
		\|\x-\x'\|,
	\end{equation}
	where both $C_{\Phi1}$ and $C_{\Phi2}$ are constants given by 
	\begin{align*}
		C_{\Phi1}:=&L_{\psi_{\y1}}^2 L_{\psi_{\y\y2}} +L_{\psi_{\y1}}\left(L_{\psi_{\y\y1}}+L_{\psi_{\x\y2}}\right)\sigma_{\psi_\mu} 
		+ L_{\psi_{\x\y1}} \sigma_{\psi_\mu} ^2,\\
		C_{\Phi2}:=&L_{\psi_{\y1}}^2  L_{F_{\y2}} 
		+  L_{\psi_{\y1}} \left(L_{F_{\y1}}+L_{F_{\x2}}\right) \sigma_{\psi_\mu} + L_{F_{\x1}} \sigma_{\psi_\mu}^2.
	\end{align*}
\end{lemma}
\begin{proof}
	Recall that
	$
	\nabla \Phi_{\mu}(\x)
	=\nabla_{\x} F(\x,\y^*_{\mu}(\x))
	-\left[\nabla_{\x\y}^2\psi_{\mu}(\x,\y^*_{\mu}(\x))\right]\v^*_{\mu}(\x).
	$
	Thus
	\begin{align*}
		\nabla \Phi_{\mu}(\x)-\nabla \Phi_{\mu}(\x')
		=& \nabla_{\x} F(\x,\y^*_{\mu}(\x))-\nabla_{\x} F(\x',\y^*_{\mu}(\x'))\\
		&+\left[\nabla_{\x\y}^2\psi_{\mu}(\x',\y^*_{\mu}(\x'))\right]\v^*_{\mu}(\x')
		-\left[\nabla_{\x\y}^2\psi_{\mu}(\x,\y^*_{\mu}(\x))\right]\v^*_{\mu}(\x).
	\end{align*}
	First, by Lipschitz continuity of $\nabla_{\x} F$, we have 
	\begin{align*}
		\|\nabla_{\x} F(\x,\y^*_{\mu}(\x))-\nabla_{\x} F(\x',\y^*_{\mu}(\x'))\|\leq & L_{F_{\x1}}\|\x-\x'\|+L_{F_{\x2}} \|\y^*_{\mu}(\x)-\y^*_{\mu}(\x')\|.
	\end{align*}
	Second, note that 
	\begin{align*}
		&\left[\nabla_{\x\y}^2\psi_{\mu}(\x',\y^*_{\mu}(\x'))\right]\v^*_{\mu}(\x')
		-\left[\nabla_{\x\y}^2\psi_{\mu}(\x,\y^*_{\mu}(\x))\right]\v^*_{\mu}(\x)\\
		=&[\nabla_{\x\y}^2\psi_{\mu}(\x',\y^*_{\mu}(\x'))]\left[\v^*_{\mu}(\x')-\v^*_{\mu}(\x)\right]
		+\left[\nabla_{\x\y}^2\psi_{\mu}(\x',\y^*_{\mu}(\x'))-\nabla_{\x\y}^2\psi_{\mu}(\x,\y^*_{\mu}(\x))\right]\v^*_{\mu}(\x).
	\end{align*}
	Thus, by Lipschitz continuity of $\nabla_{\y\y}^2 \psi_\mu$ and Lemma~\ref{lemma*}, we have 
	\begin{align*}
		&\left\|\left[\nabla_{\x\y}^2\psi_{\mu}(\x',\y^*_{\mu}(\x'))\right]\v^*_{\mu}(\x')
		-\left[\nabla_{\x\y}^2\psi_{\mu}(\x,\y^*_{\mu}(\x))\right]\v^*_{\mu}(\x)\right\|\\
		\leq& L_{\psi_{\y1}}\left\| \v^*_{\mu}(\x')-\v^*_{\mu}(\x)\right\| 
		+\|\v^*_{\mu}(\x)\| \left(L_{\psi_{\x\y1}}\|\x'-\x\|+L_{\psi_{\x\y2}}\|\y^*_{\mu}(\x)-\y^*_{\mu}(\x')\|\right)\\
		\leq& \frac{ L_{\psi_{\y1}} \left(L_{\v1} \|\v^*_{\mu}(\x)\| +L_{\v2}\right) }{\sigma_{\psi_{\mu}}^2}
		\|\x-\x'\|
		+\|\v^*_{\mu}(\x)\| \left(L_{\psi_{\x\y1}}\|\x'-\x\|+L_{\psi_{\x\y2}}\|\y^*_{\mu}(\x)-\y^*_{\mu}(\x')\|\right).
	\end{align*}
	Hence, by using Lemma~\ref{lemma*}, we have 
	\begin{align*}
		&\|\nabla \Phi_{\mu}(\x) -\nabla \Phi_{\mu}(\x')\| \\
		\leq &
		\left(L_{F_{\x1}}
		+L_{\psi_{\x\y1}} \|\v^*_{\mu}(\x)\|
		+\frac{L_{\psi_{\y1}} \left(L_{\v1} \|\v^*_{\mu}(\x)\| +L_{\v2}\right) }{\sigma_{\psi_{\mu}}^2}\right)\|\x-\x'\|
		+\left(L_{F_{\x2}}
		+L_{\psi_{\x\y2}} \|\v^*_{\mu}(\x)\|\right)\|\y^*_{\mu}(\x)-\y^*_{\mu}(\x')\|\\
		\leq & 
		\left(L_{F_{\x1}}
		+L_{\psi_{\x\y1}} \|\v^*_{\mu}(\x)\|
		+\frac{L_{\psi_{\y1}} \left(L_{\v1} \|\v^*_{\mu}(\x)\| +L_{\v2}\right) }{\sigma_{\psi_{\mu}}^2}\right)\|\x-\x'\|
		+\left( L_{F_{\x2}}
		+L_{\psi_{\x\y2}} \|\v^*_{\mu}(\x)\|\right)
		\frac{L_{\psi_{\y1}}}{\sigma_{\psi_\mu}}
		\|\x-\x'\|,
	\end{align*}
	which implies the desired result since
	$$L_{\psi_{\y1}} L_{\v1}+L_{\psi_{\x\y2}} L_{\psi_{\y1}} \sigma_{\psi_\mu}+L_{\psi_{\x\y1}} \sigma_{\psi_\mu}^2
	=L_{\psi_{\y1}}^2 L_{\psi_{\y\y2}} +L_{\psi_{\y1}}\left(L_{\psi_{\y\y1}}+L_{\psi_{\x\y2}}\right)\sigma_{\psi_\mu} 
	+ L_{\psi_{\x\y1}}\sigma_{\psi_\mu} ^2$$ 
	and
	$$L_{\psi_{\y1}} L_{\v2} + L_{F_{\x2}} L_{\psi_{\y1}} \sigma_{\psi_\mu} + L_{F_{\x1}} \sigma_{\psi_\mu}^2
	=  L_{\psi_{\y1}}^2  L_{F_{\y2}} 
	+  L_{\psi_{\y1}} \left(L_{F_{\y1}}+L_{F_{\x2}}\right) \sigma_{\psi_\mu} + L_{F_{\x1}} \sigma_{\psi_\mu}^2.$$
\end{proof}

Next lemma characterizes the smoothness of $\y^*_\mu(\x)$ and $\v^*_\mu(\x)$ in $\mu$. It will plays an important role in the non-asymptotic analysis for sl-BAMM when BLO has multiple LL minimizers.
\begin{lemma}\label{lem13}
	Suppose Assumptions \ref{Assump0}, and either Assumption \ref{assump_convex0} or Assumption \ref{assump_convex20} holds, let $0<\mu\leq 1/2$ and $\mu'\leq 2\mu$, 
	\begin{enumerate}[(i)]
		\item 
		\begin{equation}\label{y_est}
			\|\y^*_{\mu}(\x)-\y^*_{\mu'}(\x)\|
			\leq  \frac{2\|\nabla_{\y} F(\x,\y^*_{\mu}(\x))\|}{\sigma_{F}}\cdot\frac{| \mu-\mu' |}{\mu'}.
		\end{equation}
		\item 
		\begin{equation}
			\|\v^*_{\mu}(\x)-\v^*_{\mu'}(\x)\|
			\leq  
			\left(C_{\v1}  \|\v^*_{\mu}(\x)\| 
			+ C_{\v2}\right)\|\nabla_{\y} F(\x,\y^*_{\mu}(\x))\|
			\cdot
			\frac{ |\mu-\mu' |}{\mu'^2},
		\end{equation}
		where both $C_{\v1}$ and $C_{\v2}$ are constants given by 
		\begin{align*}
			C_{\v1}:=2\left(L_{F_{\y\y2}}+L_{f_{\y\y2}}\right)/\sigma_{F}^2,
			\quad
			C_{\v2}:=2\left(2L_{F_{\y2}}+L_{f_{\y2}}\right)/\sigma_{F}^2.
		\end{align*}
	\end{enumerate}
\end{lemma}
\begin{proof}
	\begin{enumerate}[(i)]
		\item
		Recall that $\y^*_{\mu}=\y^*_{\mu}(\x)$ satisfies
		$
		\mu  \nabla_{\y} F(\x,\y^*_{\mu})+(1-\mu)\nabla_{\y} f(\x,\y^*_{\mu})=0.
		$
		Thus
		\begin{align*}
			&(\mu-\mu') \nabla_{\y} F(\x,\y^*_{\mu})+\mu' \big[\nabla_{\y} F(\x,\y^*_{\mu})-\nabla_{\y} F(\x,\y^*_{\mu'})\big]\\
			&+(\mu'-\mu)\nabla_{\y} f(\x,\y^*_{\mu})
			+(1-\mu')\big[\nabla_{\y} f(\x,\y^*_{\mu})-\nabla_{\y} f(\x,\y^*_{\mu'})\big]=0,
		\end{align*}
		which implies that
		\begin{align*}
			&\mu' \big[\nabla_{\y} F(\x,\y^*_{\mu})-\nabla_{\y} F(\x,\y^*_{\mu'})\big]
			+(1-\mu')\big[\nabla_{\y} f(\x,\y^*_{\mu})-\nabla_{\y} f(\x,\y^*_{\mu'})\big] \\
			=&
			(\mu'-\mu)\nabla_{\y} F(\x,\y^*_{\mu})
			+(\mu'-\mu)\frac{\mu}{1-\mu}\nabla_{\y} F(\x,\y^*_{\mu})
			=(\mu'-\mu)\frac{ \nabla_{\y} F(\x,\y^*_{\mu})}{1-\mu}.
		\end{align*}
		Since $F(\x, \cdot)$ is $\sigma_{F}$-strongly convex, we have 
		\begin{align*}
			\big\langle \nabla_{\y} F(\x,\y^*_{\mu})-\nabla_{\y} F(\x,\y^*_{\mu'}), \ \y^*_{\mu}-\y^*_{\mu'}\big\rangle\geq \sigma_F \|\y^*_{\mu}-\y^*_{\mu'}\|^2.
		\end{align*}
		Similarly, since $f(\x, \cdot)$ is convex, we get 
		\begin{align*}
			\big\langle \nabla_{\y} f(\x,\y^*_{\mu})-\nabla_{\y} f(\x,\y^*_{\mu'}),\  \y^*_{\mu}-\y^*_{\mu'}\big\rangle\geq 0.
		\end{align*}
		Combining the above inequalities, we have 
		\begin{align*}
			0\leq \mu' \sigma_F \|\y^*_{\mu}-\y^*_{\mu'}\|^2
			\leq 
			\left(\frac{\mu'-\mu}{1-\mu}\right)\big\langle \nabla_{\y} F(\x,\y^*_{\mu}), \y^*_{\mu}-\y^*_{\mu'}\big\rangle.
		\end{align*}
		Since $0<\mu\leq 1/2$, we get $1/(1-\mu)\leq 2$ and then 
		\begin{align*}
			\mu' \sigma_F \|\y^*_{\mu}-\y^*_{\mu'}\|^2
			\leq
			2 \|\nabla_{\y} F(\x,\y^*_{\mu})\|
			\cdot|\mu-\mu'| \cdot\|\y^*_{\mu}-\y^*_{\mu'}\|,
		\end{align*}
		which implies the desired result.
		
		\item The definition of $\v^*_{\mu}=\v^*_{\mu}(\x)$ says that  $\big[\nabla_{\y\y}^2\psi_{\mu}(\x,\y^*_{\mu})\big]\v^*_{\mu}=\nabla_{\y} F(\x,\y^*_{\mu})$. Thus we have 
		\begin{align*}
			[\nabla_{\y\y}^2\psi_{\mu'}(\x,\y^*_{\mu'})](\v^*_{\mu}-\v^*_{\mu'})
			+\left[\nabla_{\y\y}^2\psi_{\mu}(\x,\y^*_{\mu})-\nabla_{\y\y}^2\psi_{\mu'}(\x,\y^*_{\mu'})\right]\v^*_{\mu}
			=\nabla_{\y} F(\x,\y^*_{\mu})-\nabla_{\y} F(\x,\y^*_{\mu'}).
		\end{align*}
		Multiplying the above equation by $\v^*_{\mu}-\v^*_{\mu'}$, by the $\sigma_{\psi_{\mu}}$-strong convexity of $\psi_{\mu}(\x,\cdot)$, we get 
		\begin{align*}
			\sigma_{\psi_{\mu'}}\|\v^*_{\mu}-\v^*_{\mu'}\|^2
			\leq&\|\nabla_{\y\y}^2\psi_{\mu}(\x,\y^*_{\mu})-\nabla_{\y\y}^2\psi_{\mu'}(\x,\y^*_{\mu'})\|
			\cdot\|\v^*_{\mu}\|\cdot\|\v^*_{\mu}-\v^*_{\mu'}\|\\
			&+\|\nabla_{\y} F(\x,\y^*_{\mu})-\nabla_{\y} F(\x,\y^*_{\mu'})\|\cdot\|\v^*_{\mu}-\v^*_{\mu'}\|.
		\end{align*}
		Note that by the definition of $\psi_{\mu}$ we have 
		\begin{align*}
			\nabla_{\y\y}^2\psi_{\mu}(\x,\y^*_{\mu})-\nabla_{\y\y}^2\psi_{\mu'}(\x,\y^*_{\mu'})
			=&(\mu-\mu') \nabla_{\y\y}^2 F(\x,\y^*_{\mu})+\mu' \left[\nabla_{\y\y}^2 F(\x,\y^*_{\mu})-\nabla_{\y\y}^2 F(\x,\y^*_{\mu'})\right]\\
			&+(\mu'-\mu)\nabla_{\y\y}^2 f(\x,\y^*_{\mu})	
			+(1-\mu')\left[\nabla_{\y\y}^2 f(\x,\y^*_{\mu})
			-\nabla_{\y\y}^2 f(\x,\y^*_{\mu'})\right].
		\end{align*}
		Since $\mu'\leq 2\mu\leq 1$, we get 
		\begin{align*}
			\sigma_{\psi_{\mu'}}\|\v^*_{\mu} -\v^*_{\mu'} \|
			\leq&\big[( L_{F_{\y2}}+L_{f_{\y2}} ) |\mu-\mu' |
			+(L_{F_{\y\y2}}+L_{f_{\y\y2}})\|\y^*_{\mu}-\y^*_{\mu'}\|\big]
			\|\v^*_{\mu}\| 
			+L_{F_{\y2}}\|\y^*_{\mu}-\y^*_{\mu'}\|\\
			\leq & \left[\left(L_{F_{\y\y2}}+L_{f_{\y\y2}}\right)\|\v^*_{\mu}\|
			+L_{F_{\y2}}\right]
			\left\| \y^*_{\mu}-\y^*_{\mu'} \right\|
			+( L_{F_{\y2}}+L_{f_{\y2}} ) \|\v^*_{\mu}\| \cdot |\mu-\mu' |.
		\end{align*}
		By (i), since $\sigma_{\psi_{\mu'}}= \mu' \sigma_{F}$ when $\sigma_{f}=0$, we have 
		\begin{align*}
			\|\v^*_{\mu}-\v^*_{\mu'}\|
			\leq&\left[\left(L_{F_{\y\y2}}+L_{f_{\y\y2}}\right)\|\v^*_{\mu}\|
			+L_{F_{\y2}}\right]
			\frac{2 \|\nabla_{\y} F(\x,\y^*_{\mu}(\x))\|}{\sigma_{F}^2}
			\cdot\frac{| \mu-\mu' |}{\mu'^2 }\\
			&+\frac{( L_{F_{\y2}}+L_{f_{\y2}} )}{\sigma_{F}} \|\v^*_{\mu}\| \cdot \frac{|\mu-\mu' |}{\mu'}.
		\end{align*}
		Since $\mu'\leq 2\mu$, Lemma~\ref{lemma*} implies that 
		\begin{align*}
			\mu' \|\v^*_{\mu}(\x)\|
			\leq \mu' \frac{\|\nabla_{\y} F(\x,\y^*_{\mu}(\x))\|}{\mu \sigma_{F}}
			\leq 
			\frac{2 \|\nabla_{\y} F(\x,\y^*_{\mu}(\x))\|}{\sigma_{F}}.
		\end{align*}
		Hence we get the desired result. 
	\end{enumerate}
\end{proof}

Note that some of the coefficients of terms in Inequality \eqref{estimate4} depends on $\| \v^*_{\mu_k}(\x_k) \|$ and $\| \nabla_{\y} F(\x_{k+1},\y^*_{\mu_{k+1}}(\x_{k+1})) \|$. We next upper-bound such two types of terms. 
\begin{lemma}\label{lemma-boundvbyxy}
	Assume that $\nabla_{\y}F$ and $\nabla_{\y} f$ are Lipschitz continuous in $(\x,\y)$. If $(\x_k, \y_k)$ is bounded, that is, $\| \x_{k}\| +\|\y_{k}\|\leq D$. Then
	\begin{equation}
		\|\nabla_{\y} F(\x_k,\y^*_{\mu_k}(\x_k))\|=O( 1/\sigma_{\psi_{\mu_{k}}} ),
		\quad 
		\|\v^*_{\mu_k}(\x_k)\|=O\left(1/\sigma_{\psi_{\mu_k}}^2\right).
	\end{equation}
	In detail, 
	\begin{equation}
		\|\nabla_{\y} F(\x_k,\y^*_{\mu_k}(\x_k))\|
		\leq \left(C_{d1} D + C_{d2}\right)/ \sigma_{\psi_{\mu_k}},
		\quad
		\|\v^*_{\mu_k}(\x_k)\|
		\leq \left(C_{d1} D + C_{d2}\right)/ \sigma_{\psi_{\mu_k}}^2,
	\end{equation}
	where
	\begin{align*}
		C_{d1}:=& L_{F_{\y_2}} \left(L_{f_{\y1}} + L_{f_{\y2}} \right)
		+ \left( \sigma_{\psi_{\mu_k}}+ \mu_{k} L_{F_{\y_2}} \right) \left(L_{F_{\y1}} + L_{F_{\y2}} \right) ,\\
		C_{d2}:=& L_{F_{\y_2}} \| \nabla_{\y} f(0,0) \| 
		+ \left( \sigma_{\psi_{\mu_k}}+ \mu_{k} L_{F_{\y_2}} \right) \| \nabla_{\y} F(0,0) \| ,    
	\end{align*}
	are constants independent of $D$.
\end{lemma}
\begin{proof}
	First, since $\| \x_{k}\| +\|\y_{k}\|\leq D$, by Lipschitz continuity of $\nabla_{\y}F$ and $\nabla_{\y} f$, we have 
	\begin{align*}
		\| \nabla_{\y} F(\x_{k},\y_k) \|
		\leq & \| \nabla_{\y} F(0,0) \| + \| \nabla_{\y} F(\x_{k},\y_k) - \nabla_{\y} F(0,0) \|\\
		\leq & \| \nabla_{\y} F(0,0) \| + L_{F_{\y1}} \|\x_{k}\| +L_{F_{\y2}} \|\y_{k} \|  \\
		\leq &  \| \nabla_{\y} F(0,0) \| + \left(L_{F_{\y1}} + L_{F_{\y2}} \right) D .
	\end{align*}
	Similarly, we also have $\| \nabla_{\y} f(\x_{k},\y_k) \| \leq \| \nabla_{\y} f(0,0) \| + \left(L_{f_{\y1}} + L_{f_{\y2}} \right) D$. 
	
	Next, since $\nabla_{\y}\psi_{\mu_k}(\x_{k},\y_{\mu_k}^*(\x_k))=0$, by the $\sigma_{\psi_{\mu_k}}$-strong convexity of $\psi_{\mu_k}(\x,\cdot)$,  we get
	\begin{align*}
		\sigma_{\psi_{\mu_k}}\left\| \y_k-\y_{\mu_k}^*(\x_k) \right\|
		\leq & \|\nabla_{\y}\psi_{\mu_k}(\x_{k},\y_k)-\nabla_{\y}\psi_{\mu_k}(\x_{k},\y_{\mu_k}^*(\x_k))\|
		= \|\nabla_{\y}\psi_{\mu_k}(\x_{k},\y_k)\|.
	\end{align*} 
	This implies that $\left\| \y_k-\y_{\mu_k}^*(\x_k) \right\|\leq  \|\nabla_{\y}\psi_{\mu_k}(\x_{k},\y_k)\|/\sigma_{\psi_{\mu_k}}$ and then 
	\begin{align*}
		\left\| \nabla_{\y} F(\x_{k},\y^*_{\mu_{k}}(\x_{k})) \right\|
		\leq & \| \nabla_{\y} F(\x_{k},\y_k) \|
		+ \| \nabla_{\y} F(\x_{k},\y^*_{\mu_{k}}(\x_{k})) 
		-\nabla_{\y} F(\x_{k},\y_k) \|\\
		\leq & \| \nabla_{\y} F(\x_{k},\y_k) \|
		+ L_{F_{\y_2}} \left\| \y_k-\y_{\mu_k}^*(\x_k) \right\| \\
		\leq & \| \nabla_{\y} F(\x_{k},\y_k) \|
		+\frac{ L_{F_{\y_2}}  \|\nabla_{\y}\psi_{\mu_k}(\x_{k},\y_k)\|}{\sigma_{\psi_{\mu_k}}}\\
		\leq &
		\left(1+\frac{ \mu_k L_{F_{\y_2}} }{ \sigma_{\psi_{\mu_k}} }  \right)\| \nabla_{\y} F(\x_{k},\y_k) \|
		+\frac{L_{F_{\y_2}} }{\sigma_{\psi_{\mu_k}}}\| \nabla_{\y} f(\x_{k},\y_k) \| , 
	\end{align*}
	since $\nabla_{\y}\psi_{\mu_k}=\mu_k \nabla_{\y} F+(1-\mu_{k}) \nabla_{\y} f$. Thus, by Lemma~\ref{lemma*}, we have 
	\begin{align*}
		\|\v^*_{\mu_k}(\x_k)\|
		\leq 
		\left(\frac{1}{\sigma_{\psi_{\mu_k}}}+\frac{ \mu_k L_{F_{\y_2}} }{\sigma_{\psi_{\mu_k}}^2 } \right)\| \nabla_{\y} F(\x_{k},\y_k) \|
		+\frac{L_{F_{\y_2}} }{\sigma_{\psi_{\mu_k}}^2}\| \nabla_{\y} f(\x_{k},\y_k) \| 
		\leq  \frac{\left(C_{d1} D + C_{d2}\right)}{\sigma_{\psi_{\mu_k}}^2}.
	\end{align*}
\end{proof}

\subsection{Proof of Lemma \ref{dlemma1}}

By the above lemmas, we can analyze the descent of the approximate overall  UL objective.

\begin{lemma}[\bf Lemma \ref{dlemma1}]\label{dlemma1a}
	Suppose Assumptions \ref{Assump0}, and either Assumption \ref{assump_convex0} or Assumption \ref{assump_convex20} holds. Let $\mu_{k+1}\leq\mu_{k}\leq\frac{1}{2}$ for all $k$. Then the sequence of $\x_k,\y_k,\v_k$ generated by sl-BAMM satisfies
	\begin{equation}\label{estimate1a}
		\begin{aligned}
			&\Phi_{\mu_{k+1}}(\x_{k+1})-\Phi_{\mu_{k}}(\x_k)  \\
			\leq & -\frac{\alpha_k}{2}\|\nabla \Phi_{\mu_{k}}(\x_k)\|^2
			-\frac{1}{2}\left(\frac{1}{\alpha_k}
			-\frac{C_{\Phi1}\|\v^*_{\mu_k}(\x_k)\|+C_{\Phi2}}{\sigma_{\psi_{\mu_k}}^2}\right)\|\x_{k+1}-\x_k\|^2 
			+\alpha_k L_{\psi_{\y1}}^2\|\v_{k}-\v^*_{\mu_{k}}(\x_k)\|^2 \\  
			&+\alpha_k \left( L_{\psi_{\x\y2}} \left\| \v^*_{\mu_{k}}(\x_k)\right\| + L_{F_{\x2}} \right)^2 
			\|\y_{k}-\y^*_{\mu_k}(\x_k)\|^2
			+\frac{2\big\|\nabla_{\y} F(\x_{k+1},\y^*_{\mu_{k+1}}(\x_{k+1}))\big\|^2}{\sigma_{F}} \left(\frac{\mu_k-\mu_{k+1}}{\mu_k}\right), 
		\end{aligned}
	\end{equation}
	where both $C_{\Phi1}$ and $C_{\Phi2}$ are constants given by 
	\begin{align*}
		C_{\Phi1}:=&L_{\psi_{\y1}}^2 L_{\psi_{\y\y2}} +L_{\psi_{\y1}}\left(L_{\psi_{\y\y1}}+L_{\psi_{\x\y2}}\right)\sigma_{\psi_\mu} 
		+ L_{\psi_{\x\y1}} \sigma_{\psi_\mu} ^2,\\
		C_{\Phi2}:=&L_{\psi_{\y1}}^2  L_{F_{\y2}} 
		+  L_{\psi_{\y1}} \left(L_{F_{\y1}}+L_{F_{\x2}}\right) \sigma_{\psi_\mu} + L_{F_{\x1}} \sigma_{\psi_\mu}^2.
	\end{align*}
	In particular, if LL objective $f(\x,\cdot)$ is $\sigma_{f}$-strongly convex, i.e., Assumption \ref{assump_convex20} holds, we take $\mu_{k}=0$ for all $k$ and then the last term in Estimate \eqref{estimate1} is redundant.
\end{lemma}
\begin{proof}
	First, note that
	\begin{align*}
		F(\x_{k+1},\y^*_{\mu_{k+1}}(\x_{k+1}))-F(\x_{k},\y^*_{\mu_k}(\x_k))=&F(\x_{k+1},\y^*_{\mu_{k+1}}(\x_{k+1}))-F(\x_{k+1},\y^*_{\mu_{k}}(\x_{k+1}))\\
		&+F(\x_{k+1},\y^*_{\mu_{k}}(\x_{k+1}))-F(\x_{k},\y^*_{\mu_{k}}(\x_{k})).
	\end{align*}
	By the convexity of $F(\x,\cdot)$ and Cauchy-Schwarz inequality, we have 
	\begin{align*}
		F(\x_{k+1},\y^*_{\mu_{k+1}}(\x_{k+1}))-F(\x_{k+1},\y^*_{\mu_{k}}(\x_{k+1}))
		\leq&\big\langle\nabla_{\y} F(\x_{k+1},\y^*_{\mu_{k+1}}(\x_{k+1})), \y^*_{\mu_{k+1}}(\x_{k+1})-\y^*_{\mu_{k}}(\x_{k+1})\big\rangle\\
		\leq&\big\|\nabla_{\y} F(\x_{k+1},\y^*_{\mu_{k+1}}(\x_{k+1}))\big\|
		\cdot\| \y^*_{\mu_{k+1}}(\x_{k+1})-\y^*_{\mu_{k}}(\x_{k+1})\|\\
		\leq &\frac{2\big\|\nabla_{\y} F(\x_{k+1},\y^*_{\mu_{k+1}}(\x_{k+1}))\big\|^2}{\sigma_{F}} \left(\frac{\mu_k-\mu_{k+1}}{\mu_k}\right),
	\end{align*}
	where the last inequality follows from \eqref{y_est} with $\mu'=\mu_{k}$ and $\mu=\mu_{k+1}$. 
	
	Next, denote $L_{\Phi\mu_{k}}:=C_{\Phi1}\|\v^*_{\mu_k}(\x_k)\| + C_{\Phi2}$. Since $L_{\Phi\mu}$ is independent of $\x'$ in Lemma \ref{lem12}, similar to the proof of the classical descent lemma (cf. Lemma 5.7 in \cite{beck2017first}), one can get 
	\begin{align*}
		&F(\x_{k+1},\y^*_{\mu_{k}}(\x_{k+1}))-F(\x_{k},\y^*_{\mu_{k}}(\x_{k})) 
		= \Phi_{\mu_{k}}(\x_{k+1})-\Phi_{\mu_{k}}(\x_k) \\
		\leq& \big\langle \nabla \Phi_{\mu_{k}}(\x_k), \x_{k+1}-\x_k \big\rangle+\frac{L_{\Phi \mu_{k}}}{2\sigma_{\psi_{\mu_k}}^2}\|\x_{k+1}-\x_k\|^2\\
		\leq & -\alpha_k \left\langle \nabla \Phi_{\mu_{k}}(\x_k), \nabla_{\x} F(\x_k,\y_{k})-\nabla_{\x\y}^2 \psi_{\mu_k}(\x_k,\y_{k})\v_{k} \right\rangle +\frac{L_{\Phi\mu_{k}}}{2\sigma_{\psi_{\mu_k}}^2}
		\|\x_{k+1}-\x_k\|^2\\
		\leq & -\frac{\alpha_k}{2}\|\nabla \Phi_{\mu_{k}}(\x_k)\|^2
		-\frac{\alpha_k}{2}
		\left\|\nabla_{\x} F(\x_k,\y_{k})-\nabla_{\x\y}^2 \psi_{\mu_k}(\x_k,\y_{k})\v_{k}\right\|^2\\
		& +\frac{\alpha_k}{2}\left\|\nabla \Phi_{\mu_{k}}(\x_k)-\nabla_{\x} F(\x_k,\y_{k})+\nabla_{\x\y}^2 \psi_{\mu_k}(\x_k,\y_{k})\v_{k}\right\|^2
		+\frac{L_{\Phi\mu_{k}}}{2\sigma_{\psi_{\mu_k}}^2}\|\x_{k+1}-\x_k\|^2,
	\end{align*}
	where the second inequality used the definition of $\x_{k+1}$ and the third one follows from $2\langle a, b\rangle=\| a \|^2 + \| b \|^2 - \| a-b \|^2$. Note that the update of $\x_{k+1}$ implies that 
	\begin{equation*}
		\left\|\nabla_{\x} F(\x_k,\y_{k})-\nabla_{\x\y}^2 \psi_{\mu_k}(\x_k,\y_{k}) \v_{k}\right\|^2=\frac{1}{\alpha_k^2}\|\x_{k+1}-\x_k\|^2,
	\end{equation*}
	and by the definition of $\Phi_{\mu_{k}}(\x)$ we get
	\begin{equation*}
		\nabla \Phi_{\mu_k}(\x_k)=\nabla_{\x} F(\x_k,\y^*_{\mu_k}(\x_k))-\nabla_{\x\y}^2\psi_{\mu_k}(\x_k,\y^*_{\mu_k}(\x_k)) \v^*_{\mu_{k}}(\x_k),
	\end{equation*}
	which implies that 
	\begin{align*}
		&\left\|\nabla \Phi_{\mu_{k}}(\x_k)-\nabla_{\x} F(\x_k,\y_{k})+\nabla_{\x\y}^2 \psi_{\mu_k}(\x_k,\y_{k})\v_{k}\right\|\\
		\leq  & \|\nabla_{\x} F(\x_k, \y^*_{\mu_k}(\x_k))-\nabla_{\x} F(\x_k,\y_{k})\|
		+\left\|\left[\nabla_{\x\y}^2 \psi_{\mu_{k}}(\x_k,\y_k)
		-\nabla_{\x\y}^2 \psi_{\mu_{k}}(\x_k, \y^*_{\mu_{k}}(\x_k))\right]  \v^*_{\mu_{k}}(\x_k)\right\|\\
		&+\left\|\nabla_{\x\y}^2 \psi_{\mu_{k}}(\x_k, \y_k)[\v_{k}-\v^*_{\mu_k}(\x_k)]\right\|\\
		\leq & 
		\left(L_{F_{\x2}} + L_{\psi_{\x\y2}} \left\| \v^*_{\mu_{k}}(\x_k)\right\|\right) \|\y_{k}-\y^*_{\mu_k}(\x_k)\|
		+L_{\psi_{\y1}} \|\v_{k}-\v^*_{\mu_k}(\x_k)\| .
	\end{align*}
	The desired result follows from the inequality $\left(\sum_{i=1}^{r} a_i\right)^2\leq r\sum_{i=1}^{r} a_i^2$. 
\end{proof}

\subsection{Proof of Lemma \ref{dlemma2}}

Lemma~\ref{dlemma1a} implies that the descent of the approximate overall  UL objective depends on the errors of LL and multiplier variables. We next analyze the errors of LL and multiplier variables. 

\begin{lemma}\label{lemmayv}
	Suppose Assumptions \ref{Assump0}, and either Assumption \ref{assump_convex0} or Assumption \ref{assump_convex20} holds.  If we choose 
	\begin{equation}
		\beta_k\leq\frac{2}{\sigma_{\psi_{\mu_k}}+L_{\psi_{\mu_k}}},\quad 
		\eta_k\leq \frac{1}{L_{\psi_{\mu_k}}}.
	\end{equation}
	Then the sequence of $\x_k,\y_k,\v_k$ generated by sl-BAMM satisfies
	\begin{align}
		\left\|\y_{k+1}-\y^*_{\mu_k}(\x_k)\right\|^2
		\leq 
		\left(1- \beta_k \sigma_{\psi_{\mu_k}} \right)
		\left\|\y_k-\y^*_{\mu_k}(\x_k)\right\|^2,
	\end{align}
	and 
	\begin{align*}
		\left\|\v_{k+1}-\v^*_{\mu_k}(\x_k)\right\|^2
		\leq &
		\left(1-\eta_k \sigma_{\psi_{\mu_k}}\right) \left\|\v_k-\v^*_{\mu_k}(\x_k)\right\|^2
		+\frac{2\eta_{k}}{ \sigma_{\psi_{\mu_k}}}
		\left(L_{\psi_{\y\y2}} \|\v^*_{\mu_k}(\x_k)\| 
		+ L_{F_{\y2}}\right)^2
		\left\| \y_{k}-\y^*_{\mu_k}(\x_k) \right\|^2 .
	\end{align*}
\end{lemma}
\begin{proof}
	First, by the update of $\y_{k+1}$, we have 
	\begin{align*}
		\left\|\y_{k+1}-\y^*_{\mu_k}(\x_k)\right\|^2
		=&\left\|\y_k-\beta_k\nabla_{\y} \psi_{\mu_k}(\x_k,\y_k)-\y^*_{\mu_k}(\x_k)\right\|^2\\
		=&\left\|\y_k-\y^*_{\mu_k}(\x_k)\right\|^2
		-2\beta_k\left\langle \y_k-\y^*_{\mu_k}(\x_k), \nabla_{\y} \psi_{\mu_k}(\x_k,\y_k)\right\rangle
		+\beta_k^2 \left\|\nabla_{\y} \psi_{\mu_k}(\x_k,\y_k)\right\|^2.
	\end{align*}
	By the $\sigma_{\psi_{\mu_k}}$-strong convexity and $L_{\psi_{\mu_k}}$-smoothness of $\psi_{\mu_k}(\x,\cdot)$, since $\nabla_{\y} \psi_{\mu_k}(\x_k,\y^*_{\mu_k}(\x_k))=0$, Theorem 2.1.12 in \cite{nesterov2018lectures} implies that 
	\begin{align*}
		\left\langle \y_k-\y^*_{\mu_k}(\x_k), \nabla_{\y} \psi_{\mu_k}(\x_k,\y_k)\right\rangle
		=&\left\langle \y_k-\y^*_{\mu_k}(\x_k), \nabla_{\y} \psi_{\mu_k}(\x_k,\y_k)-\nabla_{\y} \psi_{\mu_k}(\x_k,\y^*_{\mu_k}(\x_k))\right\rangle\\
		\geq&\frac{\sigma_{\psi_{\mu_k}} L_{\psi_{\mu_k}}}{\sigma_{\psi_{\mu_k}}+L_{\psi_{\mu_k}}}\left\|\y_k-\y^*_{\mu_k}(\x_k)\right\|^2
		+\frac{1}{\sigma_{\psi_{\mu_k}}+L_{\psi_{\mu_k}}}
		\left\|\nabla_{\y} \psi_{\mu_k}(\x_k,\y_k)\right\|^2.
	\end{align*}
	Hence, when $0<\beta_k\leq\frac{2}{\sigma_{\psi_{\mu_k}}+L_{\psi_{\mu_k}}}$, we have  
	\begin{align*}
		\left\|\y_{k+1}-\y^*_{\mu_k}(\x_k)\right\|^2
		\leq & \left(1-\frac{2\sigma_{\psi_{\mu_k}} L_{\psi_{\mu_k}}}{\sigma_{\psi_{\mu_k}}+L_{\psi_{\mu_k}}} \beta_k\right)
		\left\|\y_k-\y^*_{\mu_k}(\x_k)\right\|^2,
	\end{align*}
	which implies the desired result since $\frac{1}{2}\leq\frac{ L_{\psi_{\mu_{k}}}}{\sigma_{\psi_{\mu_{k}}}+L_{\psi_{\mu_{k}}}} \leq 1$.
	
	Second, by the update of $\v_{k+1}$, using $\left[\nabla_{\y\y}^2 \psi_{\mu_k}(\x_k,\y^*_{\mu_k}(\x_k))\right]\v^*_{\mu_k}(\x_k)=\nabla_{\y} F(\x_k,\y^*_{\mu_k}(\x_k))$, 
	\begin{align*}
		\v_{k+1}-\v^*_{\mu_k}(\x_k)
		=&\v_k-\v^*_{\mu_k}(\x_k)-\eta_k\left([\nabla_{\y\y}^2 \psi_{\mu_k}(\x_k,\y_{k})]\v_k-\nabla_{\y} F(\x_k,\y_{k})\right)\\
		=&\v_k-\v^*_{\mu_k}(\x_k)
		-\eta_k\left[\nabla_{\y\y}^2 \psi_{\mu_k}(\x_k,\y_{k})\right]\left[\v_k-\v^*_{\mu_k}(\x_k)\right]\\
		&-\eta_k \left[\nabla_{\y\y}^2 \psi_{\mu_k}(\x_k,\y_{k})-\nabla_{\y\y}^2 \psi_{\mu_k}(\x_k,\y^*_{\mu_k}(\x_k))\right] \v^*_{\mu_k}(\x_k)\\
		&-\eta_k \left[\nabla_{\y} F(\x_k,\y^*_{\mu_k}(\x_k))-\nabla_{\y} F(\x_k,\y_{k})\right].
	\end{align*}
	Hence, by Cauchy-Schwarz inequality, for all $\varepsilon>0$, we have 
	\begin{align*}
		& \left\|\v_{k+1}-\v^*_{\mu_k}(\x_k)\right\|^2\\
		\leq &
		\left(1+\varepsilon\right)\left\| \left[I-\eta_k \nabla_{\y\y}^2 \psi_{\mu_k}(\x_k,\y_{k})\right]\left[\v_k-\v^*_{\mu_k}(\x_k)\right]\right\|^2\\
		&+ \left(1+\frac{1}{\varepsilon}\right) \eta_{k}^2 \Big\| \left[\nabla_{\y\y}^2 \psi_{\mu_k}(\x_k,\y_{k})-\nabla_{\y\y}^2 \psi_{\mu_k}(\x_k,\y^*_{\mu_k}(\x_k))\right] \v^*_{\mu_k}(\x_k)
		+\nabla_{\y} F(\x_k,\y^*_{\mu_k}(\x_k))-\nabla_{\y} F(\x_k,\y_{k}) \Big\|^2 .
	\end{align*}	
	When $\eta_k\leq 1/L_{\psi_{\mu_k}}$, by the $\sigma_{\psi_{\mu_k}}$-strong convexity of $\psi_{\mu_k}(\x,\cdot)$, since the spectral norm is monotone, we have 
	\begin{align*}
		\left\| \left[I-\eta_k \nabla_{\y\y}^2 \psi_{\mu_k}(\x_k,\y_{k+1})\right]\left[\v_k-\v^*_{\mu_k}(\x_k)\right]\right\|
		\leq \left(1-\eta_k \sigma_{\psi_{\mu_k}}\right) \left\|\v_k-\v^*_{\mu_k}(\x_k)\right\|.
	\end{align*}
	Next, by Lipschitz continuity of $\nabla_{\y\y}^2 \psi_{\mu_k} (\x, \cdot)$ and $\nabla_{\y} F (\x, \cdot)$, we get 
	\begin{align*}
		&\big\| \left[\nabla_{\y\y}^2 \psi_{\mu_k}(\x_k,\y_{k})-\nabla_{\y\y}^2 \psi_{\mu_k}(\x_k,\y^*_{\mu_k}(\x_k))\right] \v^*_{\mu_k}(\x_k)
		+\nabla_{\y} F(\x_k,\y^*_{\mu_k}(\x_k))-\nabla_{\y} F(\x_k,\y_{k})\big\|\\
		\leq & \left(L_{\psi_{\y\y2}} \|\v^*_{\mu_k}(\x_k)\| + L_{F_{\y2}}\right)
		\left\| \y_{k}-\y^*_{\mu_k}(\x_k) \right\|.
	\end{align*}
	Taking $\varepsilon=\eta_k \sigma_{\psi_{\mu_k}}$, we have 
	\begin{align*}
		\left\|\v_{k+1}-\v^*_{\mu_k}(\x_k)\right\|^2
		\leq &
		\left(1+\eta_k \sigma_{\psi_{\mu_k}}\right)\left(1-\eta_k \sigma_{\psi_{\mu_k}}\right)^2 \left\|\v_k-\v^*_{\mu_k}(\x_k)\right\|^2\\
		&+\left(1+\frac{1}{\eta_k \sigma_{\psi_{\mu_k}}}\right) \eta_{k}^2\left(L_{\psi_{\y\y2}} \|\v^*_{\mu_k}(\x_k)\| + L_{F_{\y2}}\right)^2
		\left\| \y_{k}-\y^*_{\mu_k}(\x_k) \right\|^2,
	\end{align*}
	which implies the desired result since $\eta_{k}^2\leq \eta_{k}/L_{\psi_{\mu_k}}\leq \eta_{k}/\sigma_{\psi_{\mu_k}}$. 
\end{proof}

\begin{lemma}\label{lemmay*v*}
	Suppose Assumptions in Lemma~\ref{lem13} hold. Then the sequence of $\x_k,\y_k,\v_k$ generated by sl-BAMM satisfies
	\begin{align*}
		\| \y^*_{\mu_k}(\x_k)-\y^*_{\mu_{k+1}}(\x_{k+1}) \|^2
		\leq & \frac{2L_{\psi_{\y1}}^2}{\sigma_{\psi_{\mu_k}}^2}
		\|\x_k-\x_{k+1}\|^2
		+\frac{8 \| \nabla_{\y} F(\x_{k+1},\y^*_{\mu_{k+1}}(\x_{k+1})) \|^2}{\sigma_F^2}
		\left(\frac{\mu_{k}-\mu_{k+1}}{\mu_k}\right)^2,\\
		\left\|\v^*_{\mu_k}(\x_k)-\v^*_{\mu_{k+1}}(\x_{k+1})\right\|^2
		\leq &
		\frac{2\left(L_{\v1} \|\v^*_{\mu_{k}}(\x_{k})\|+L_{\v2}\right)^2}{\sigma_{\psi_{\mu_k}}^4}
		\| \x_k-\x_{k+1} \|^2\\
		&+2\left(C_{\v1} \|\v^*_{\mu_{k+1}}(\x_{k+1})\|
		+C_{\v2}\right)^2  
		\|\nabla_{\y} F(\x_{k+1},\y^*_{\mu_{k+1}}(\x_{k+1}))\|^2
		\left(\frac{\mu_{k}-\mu_{k+1}}{\mu_{k}^2}\right)^2.
	\end{align*}
\end{lemma}
\begin{proof}
	By the triangle inequality, and Lemma \ref{lemma*}(i) and \ref{lem13}(i) with $\mu'=\mu_{k}$ and $\mu=\mu_{k+1}$, 
	\begin{align*}
		\| \y^*_{\mu_k}(\x_k)-\y^*_{\mu_{k+1}}(\x_{k+1}) \|\leq&
		\| \y^*_{\mu_k}(\x_k)-\y^*_{\mu_k}(\x_{k+1}) \|
		+ \| \y^*_{\mu_k}(\x_{k+1})-\y^*_{\mu_{k+1}}(\x_{k+1}) \| \\
		\leq & \frac{L_{\psi_{\y1}}}{\sigma_{\psi_{\mu_k}}}\|\x_k-\x_{k+1}\|
		+\frac{2 \| \nabla_{\y} F(\x_{k+1},\y^*_{\mu_{k+1}}(\x_{k+1})) \| }{\sigma_F}\left(\frac{\mu_{k}-\mu_{k+1}}{\mu_k}\right),
	\end{align*}
	and by Lemma \ref{lemma*}(iii) and \ref{lem13}(ii) with $\mu'=\mu_{k}$ and $\mu=\mu_{k+1}$,
	\begin{align*}
		&\| \v^*_{\mu_k}(\x_k)-\v^*_{\mu_{k+1}}(\x_{k+1}) \|
		\leq
		\| \v^*_{\mu_k}(\x_k)-\v^*_{\mu_k}(\x_{k+1}) \|
		+ \| \v^*_{\mu_k}(\x_{k+1})-\v^*_{\mu_{k+1}}(\x_{k+1}) \| \\
		\leq & \frac{\left(L_{\v1} \|\v^*_{\mu_k}(\x_k)\|  +L_{\v2}\right)}{\sigma_{\psi_{\mu_k}}^2}
		\| \x_k-\x_{k+1} \| 
		+\left(C_{\v1} \|\v^*_{\mu_{k+1}}(\x_{k+1})\|
		+C_{\v2}\right)  
		\|\nabla_{\y} F(\x_{k+1},\y^*_{\mu_{k+1}}(\x_{k+1}))\|
		\left(\frac{\mu_{k}-\mu_{k+1}}{\mu_{k}^2}\right) .
	\end{align*}
	Then the desired result follows from the inequality $\left(\sum_{i=1}^{2} a_i\right)^2\leq 2\sum_{i=1}^{2} a_i^2$.
\end{proof}

By the above lemmas, We can next analyze the error of LL and multiplier variables. 

\begin{lemma}\label{dlemma2a}
	Suppose Assumptions in Lemma~\ref{dlemma1} hold.  If we choose 
	\begin{equation*}
		\beta_k\leq\frac{2}{\sigma_{\psi_{\mu_{k}}}+L_{\psi_{\mu_{k}}}},\quad 
		\eta_k\leq \frac{1}{L_{\psi_{\mu_k}}}.
	\end{equation*}
	Then the sequence of $\x_k,\y_k,\v_k$ generated by  sl-BAMM satisfies
	\begin{equation}\label{estimate2a}
		\begin{aligned}
			\| \y_{k+1}-\y^*_{\mu_{k+1}}(\x_{k+1}) \|^2
			- \| \y_k-\y^*_{\mu_k}(\x_k) \|^2
			\leq &
			-\frac{1}{2} \beta_k \sigma_{\psi_{\mu_{k}}}
			\| \y_k-\y^*_{\mu_k}(\x_k) \|^2
			+\frac{6L_{\psi_{\y1}}^2}{\beta_{k}\sigma_{\psi_{\mu_k}}^3}
			\|\x_k-\x_{k+1}\|^2\\
			&+\frac{24 \| \nabla_{\y} F(\x_{k+1},\y^*_{\mu_{k+1}}(\x_{k+1})) \|^2}{\sigma_F^2 \beta_{k} \sigma_{\psi_{\mu_{k}}} }
			\left(\frac{\mu_{k}-\mu_{k+1}}{\mu_k}\right)^2,
		\end{aligned}
	\end{equation}
	and 
	\begin{equation}\label{estimate3a}
		\begin{aligned}
			& \| \v_{k+1}-\v^*_{\mu_{k+1}}(\x_{k+1}) \|^2
			- \| \v_k-\v^*_{\mu_{k}}(\x_k) \|^2 \\
			\leq &
			-\frac{1}{2}\eta_k\sigma_{\psi_{\mu_{k}}}
			\| \v_k-\v^*_{\mu_{k}}(\x_k) \|^2
			+\frac{6\left(L_{\v1} \|\v^*_{\mu_{k}}(\x_{k})\|+L_{\v2}\right)^2}{\eta_{k}\sigma_{\psi_{\mu_k}}^5}
			\|\x_k-\x_{k+1}\|^2 \\
			&+\frac{3\eta_{k}}{ \sigma_{\psi_{\mu_k}}}   
			\left(L_{\psi_{\y\y2}} \|\v^*_{\mu_k}(\x_k)\|
			+ L_{F_{\y2}}\right)^2
			\| \y_{k}-\y^*_{\mu_k}(\x_k) \|^2\\
			&+\frac{6\left(C_{\v1} \|\v^*_{\mu_{k+1}}(\x_{k+1})\|
				+C_{\v2}\right)^2  
				\|\nabla_{\y} F(\x_{k+1},\y^*_{\mu_{k+1}}(\x_{k+1}))\|^2 }{\eta_k\sigma_{\psi_{\mu_{k}}}}
			\left(\frac{\mu_{k}-\mu_{k+1}}{\mu_{k}^2}\right)^2,
		\end{aligned}
	\end{equation}
	where both $C_{\v1}$ and $C_{\v2}$ are constants given by 
	\begin{align*}
		C_{\v1}:=2\left(L_{F_{\y\y2}}+L_{f_{\y\y2}}\right)/\sigma_{F}^2,
		\quad
		C_{\v2}:=2\left(2L_{F_{\y2}}+L_{f_{\y2}}\right)/\sigma_{F}^2.
	\end{align*}
	In particular, if LL objective $f(\x,\cdot)$ is $\sigma_{f}$-strongly convex, we take $\mu_{k}=0$ for all $k$ and then the last terms in Estimates \eqref{estimate2a}-\eqref{estimate3a} are both redundant.
\end{lemma}
\begin{proof}
	By Cauchy-Schwarz inequality, it is easy to check that for any $\varepsilon>0$ we have 
	\begin{align*}
		\| \y_{k+1}-\y^*_{\mu_{k+1}}(\x_{k+1}) \|^2
		=&\| \y_{k+1}-\y^*_{\mu_{k}}(\x_k)+\y^*_{\mu_{k}}(\x_k)-\y^*_{\mu_{k+1}}(\x_{k+1}) \|^2\\
		\leq&(1+\varepsilon) \| \y_{k+1}-\y^*_{\mu_{k}}(\x_k) \|^2
		+\left(1+\frac{1}{\varepsilon}\right)
		\| \y^*_{\mu_{k}}(\x_k)-\y^*_{\mu_{k+1}}(\x_{k+1}) \|^2.
	\end{align*}
	Taking $\varepsilon= \frac{1}{2} \beta_k \sigma_{\psi_{\mu_{k}}}$, by Lemmas~\ref{lemmayv} and \ref{lemmay*v*}, we have 
	\begin{align*}
		\| \y_{k+1}-\y^*_{\mu_{k+1}}(\x_{k+1}) \|^2
		\leq &
		\left(1-\frac{1}{2} \beta_k \sigma_{\psi_{\mu_{k}}}\right)\|\y_k-\y^*_{\mu_k}(\x_k)\|^2
		+\frac{2L_{\psi_{\y1}}^2}{\sigma_{\psi_{\mu_k}}^2}
		\left(1+\frac{2}{ \beta_{k} \sigma_{\psi_{\mu_{k}}} }\right)
		\|\x_k-\x_{k+1}\|^2\\
		&+\frac{8 \| \nabla_{\y} F(\x_{k+1},\y^*_{\mu_{k+1}}(\x_{k+1})) \|^2}{\sigma_F^2}
		\left(1+\frac{2}{ \beta_{k} \sigma_{\psi_{\mu_{k}}} }\right)
		\left(\frac{\mu_{k}-\mu_{k+1}}{\mu_k}\right)^2,
	\end{align*}
	which implies the desired result since $\beta_k\sigma_{\psi_{\mu_{k}}}\leq 1$ when $\beta_k\leq 2/(\sigma_{\psi_{\mu_{k}}}+L_{\psi_{\mu_{k}}})$. 
	
	Similarly, for any $\delta>0$, by Cauchy-Schwarz inequality,
	\begin{align*}
		\left\|\v_{k+1}-\v^*_{\mu_{k+1}}(\x_{k+1})\right\|^2
		=&\left\|\v_{k+1}-\v^*_{\mu_{k}}(\x_k)+\v^*_{\mu_{k}}(\x_k)-\v^*_{\mu_{k+1}}(\x_{k+1})\right\|^2\\
		\leq&(1+\delta)\left\|\v_{k+1}-\v^*_{\mu_{k}}(\x_k)\right\|^2
		+\left(1+\frac{1}{\delta}\right)\left\|\v^*_{\mu_{k}}(\x_k)-\v^*_{\mu_{k+1}}(\x_{k+1})\right\|^2.
	\end{align*}
	Taking $\delta= \frac{1}{2}\eta_k\sigma_{\psi_{\mu_{k}}}$, Lemmas~\ref{lemmayv} and \ref{lemmay*v*} imply that
	\begin{align*}
		&\left\|\v_{k+1}-\v^*_{\mu_{k+1}}(\x_{k+1})\right\|^2\\
		\leq &
		\left(1-\frac{1}{2}\eta_k\sigma_{\psi_{\mu_{k}}}\right)
		\left\|\v_k-\v^*_{\mu_{k}}(\x_k)\right\|^2
		+\frac{2\eta_{k}}{ \sigma_{\psi_{\mu_k}}}    \left(1+\frac{1}{2}\eta_k\sigma_{\psi_{\mu_{k}}}\right)
		\left(L_{\psi_{\y\y2}} \|\v^*_{\mu_k}(\x_k)\|
		+ L_{F_{\y2}}\right)^2
		\left\| \y_{k}-\y^*_{\mu_k}(\x_k) \right\|^2\\
		&+\frac{2\left(L_{\v1} \|\v^*_{\mu_{k}}(\x_{k})\|+L_{\v2}\right)^2}{\sigma_{\psi_{\mu_k}}^4}
		\left(1+\frac{2}{\eta_k\sigma_{\psi_{\mu_{k}}}}\right)
		\|\x_k-\x_{k+1}\|^2 \\
		&+ 2\left(C_{\v1} \|\v^*_{\mu_{k+1}}(\x_{k+1})\|
		+C_{\v2}\right)^2  
		\|\nabla_{\y} F(\x_{k+1},\y^*_{\mu_{k+1}}(\x_{k+1}))\|^2
		\left(1+\frac{2}{\eta_k\sigma_{\psi_{\mu_{k}}}}\right)
		\left(\frac{\mu_{k}-\mu_{k+1}}{\mu_{k}^2}\right)^2.
	\end{align*}
	The desired result follows since $\eta_k\sigma_{\psi_{\mu_{k}}}\leq 1$ when $\eta_k\leq 1/L_{\psi_{\mu_k}}$. 
\end{proof}

\subsection{Proof of Lemma \ref{dlemma3}}

\begin{lemma}\label{dlemma3a}
	Suppose Assumptions in Lemma \ref{dlemma2} hold. 
	Let $\{ a_k, b_k, c_k \}_{k=1}^\infty$ be a sequence of decreasing positive constants, then the following descent of Lyapunov function holds: 
	\begin{equation}\label{estimate4a}
		\begin{aligned}
			V_{k+1}-V_k\leq & -\frac{a_{k+1} \alpha_k}{2}\|\nabla \Phi_{\mu_{k}}(\x_k)\|^2-\frac{\widehat{\alpha}_k}{2}\|\x_{k+1}-\x_k\|^2	
			-\frac{\widehat{\beta}_k}{2}\|\y_k-\y^*_{\mu_{k}}(\x_k)\|^2-\frac{\widehat{\eta}_k}{2}\|\v_k-\v^*_{\mu_{k}}(\x_k)\|^2 
			\\
			&+\frac{2 \big\|\nabla_{\y} F(\x_{k+1},\y^*_{\mu_{k+1}}(\x_{k+1}))\big\|^2 a_{k+1}}{\sigma_{F}} \left(\frac{\mu_k-\mu_{k+1}}{\mu_k}\right) 
			\\
			&+\frac{24 \left\| \nabla_{\y} F(\x_{k+1},\y^*_{\mu_{k+1}}(\x_{k+1})) \right\|^2 b_{k+1}}{\sigma_F^2 \beta_k \sigma_{\psi_{\mu_{k}}}}
			\left(\frac{\mu_{k}-\mu_{k+1}}{\mu_k}\right)^2\\
			&+\frac{6\left(C_{\v1} \|\v^*_{\mu_{k+1}}(\x_{k+1})\|
				+C_{\v2}\right)^2  
				\|\nabla_{\y} F(\x_{k+1},\y^*_{\mu_{k+1}}(\x_{k+1}))\|^2 c_{k+1}}{\eta_k\sigma_{\psi_{\mu_{k}}}}
			\left(\frac{\mu_{k}-\mu_{k+1}}{\mu_{k}^2}\right)^2,
		\end{aligned}
	\end{equation}
	where the coefficients $\{ \widehat{\alpha}_k, \widehat{\beta}_k, \widehat{\eta}_k \}$ are given as below.
	\begin{equation}\label{keycoofficients0a}
		\begin{aligned}
			\widehat{\alpha}_k
			:=&\frac{a_{k+1}}{\alpha_k}
			-\frac{ a_{k+1} \left(C_{\Phi1}\|\v^*_{\mu_k}(\x_k)\| + C_{\Phi2}\right) }{\sigma_{\psi_{\mu_{k}}}^2}
			-\frac{12L_{\psi_{\y1}}^2 b_{k+1} }{\beta_{k}\sigma_{\psi_{\mu_k}}^3}
			-\frac{ 12 \left(L_{\v1} \|\v^*_{\mu_{k}}(\x_{k})\|+L_{\v2}\right)^2 c_{k+1} }{ \eta_{k}\sigma_{\psi_{\mu_k}}^5 } ,\\
			\widehat{\beta}_k
			:=&
			b_{k+1}\beta_k \sigma_{\psi_{\mu_{k}}}
			-2 a_{k+1}\alpha_k \left( L_{\psi_{\x\y2}} \| \v^*_{\mu_k}(\x_k) \| + L_{F_{\x2}} \right)^2 
			-6 
			\left(L_{\psi_{\y\y2}} \| \v^*_{\mu_k}(\x_k) \| 
			+ L_{F_{\y2}}\right)^2
			\frac{ c_{k+1}\eta_{k}}{ \sigma_{\psi_{\mu_k}}} , \\
			\widehat{\eta}_k
			:=&
			c_{k+1}\eta_k\sigma_{\psi_{\mu_{k}}}
			-2 a_{k+1} \alpha_k  L_{\psi_{\y1}}^2.
		\end{aligned}
	\end{equation}
	Here $L_{\v1}:=L_{\psi_{\y\y2}} L_{\psi_{\y1}}+L_{\psi_{\y\y1}} \sigma_{\psi_{\mu}} $ and $L_{\v2}:=L_{F_{\y2}} L_{\psi_{\y1}}
	+L_{F_{\y1}} \sigma_{\psi_{\mu}} $. In particular, if LL objective $f(\x,\cdot)$ is $\sigma_{f}$-strongly convex, we take $\mu_{k}=0$ for all $k$ and then the terms involving $\mu_{k}-\mu_{k+1}$ in Estimate \eqref{estimate4a} are all redundant. 
\end{lemma}
\begin{proof}
	Denote $L_{\Phi\mu_{k}}:=C_{\Phi1}\|\v^*_{\mu_k}(\x_k)\| + C_{\Phi2}$, $\nabla_{\y} F_{k+1}:=\nabla_{\y} F(\x_{k+1},\y^*_{\mu_{k+1}}(\x_{k+1}))$, $B_{\v\x}:=6\left(L_{\v1} \|\v^*_{\mu_{k}}(\x_{k})\|+L_{\v2}\right)^2$, $B_{\v\y}:=3\left(L_{\psi_{\y\y2}} \| \v^*_{\mu_k}(\x_k) \| 
	+ L_{F_{\y2}}\right)^2$, and $B_{\v\mu_{k+1}}:=6\left(C_{\v1} \|\v^*_{\mu_{k+1}}(\x_{k+1})\|
	+C_{\v2}\right)^2  
	\|\nabla_{\y} F(\x_{k+1},\y^*_{\mu_{k+1}}(\x_{k+1}))\|^2$. Then,
	by Lemmas \ref{dlemma1a}-\ref{dlemma2a}, we have 
	\begin{align*}
		V_{k+1}-V_k
		=  & a_{k+1} \left[ F(\x_{k+1}, \y_{\mu_{k+1}}^*(\x_{k+1}) ) - F(\x_{k}, \y_{\mu_{k}}^*(\x_{k}) ) \right] + (a_{k+1}-a_{k}) F(\x_{k}, \y_{\mu_{k}}^*(\x_{k}) )\\
		& + b_{k+1} \left[ \|\y_{k+1}-\y^*_{\mu_{k+1}}(\x_{k+1})\|^2 - \|\y_k-\y^*_{\mu_k}(\x_k)\|^2 \right] + (b_{k+1}-b_{k}) \|\y_k-\y^*_{\mu_k}(\x_k)\|^2\\
		& + c_{k+1} \left[ \|\v_{k+1}-\v^*_{\mu_{k+1}}(\x_{k+1})\|^2 - \|\v_k-\v^*_{\mu_k}(\x_k)\|^2 \right] + (c_{k+1}-c_{k}) \|\v_k-\v^*_{\mu_k}(\x_k)\|^2\\
		\leq & -\frac{ a_{k+1} \alpha_k}{2}\|\nabla \Phi_{\mu_{k}}(\x_k)\|^2
		-\frac{a_{k+1}}{2} \left(\frac{1}{\alpha_k}
		-\frac{L_{\Phi\mu_{k}}}{\sigma_{\psi_{\mu_k}}^2}\right)
		\| \x_{k+1}-\x_k \|^2
		+ a_{k+1} \alpha_k L_{\psi_{\y1}}^2
		\| \v_{k}-\v^*_{\mu_{k}}(\x_k) \|^2   \nonumber\\
		&+a_{k+1} \alpha_k \left( L_{\psi_{\x\y2}} 
		\| \v^*_{\mu_k}(\x_k) \| + L_{F_{\x2}} \right)^2 
		\|\y_{k}-\y^*_{\mu_k}(\x_k)\|^2
		+\frac{2 a_{k+1} \big\|\nabla_{\y} F_{k+1}\big\|^2}{\sigma_{F}} \left(\frac{\mu_k-\mu_{k+1}}{\mu_k}\right) \\ 
		&-\frac{b_{k+1}}{2} \beta_k \sigma_{\psi_{\mu_{k}}}
		\|\y_k-\y^*_{\mu_k}(\x_k)\|^2
		+\frac{6 b_{k+1} L_{\psi_{\y1}}^2}{\beta_{k}\sigma_{\psi_{\mu_k}}^3}
		\| \x_k-\x_{k+1} \|^2
		+\frac{24 b_{k+1} \| \nabla_{\y} F_{k+1} \|^2}{\sigma_F^2 \beta_k \sigma_{\psi_{\mu_{k}}}}
		\left(\frac{\mu_{k}-\mu_{k+1}}{\mu_k}\right)^2\\
		&-\frac{c_{k+1}}{2}\eta_k\sigma_{\psi_{\mu_{k}}}
		\| \v_k-\v^*_{\mu_{k}}(\x_k) \|^2
		+\frac{ c_{k+1} B_{\v\y} \eta_{k}}{ \sigma_{\psi_{\mu_k}}}
		\| \y_{k}-\y^*_{\mu_k}(\x_k) \|^2 \\
		& +\frac{ c_{k+1} B_{\v\x}}{\eta_{k}\sigma_{\psi_{\mu_k}}^5}
		\|\x_k-\x_{k+1}\|^2 
		+\frac{ c_{k+1} B_{\v\mu_{k+1}}}{\eta_k\sigma_{\psi_{\mu_{k}}}}
		\left(\frac{\mu_{k}-\mu_{k+1}}{\mu_{k}^2}\right)^2,
	\end{align*}
	which implies the desired result.
\end{proof}

\section{Experiments}
\label{section C}

Our experiments were conducted on a PC with Intel i7-9700K CPU (4.2 GHz), 32GB RAM, and NVIDIA RTX 2060 GPU, and the platform is 64-bit Ubuntu 18.04.5 LTS.

\subsection{Numerical Example}

\begin{table}[htbp]
	\centering
	
	\caption{Values for hyperparameters of numerical examples.}
	\label{tab:hyper list0}
	\begin{minipage}[t]{0.4\textwidth}
		\renewcommand\arraystretch{1.2}
		\setlength{\tabcolsep}{1.2mm}{
			
			\begin{tabular}{ll}
				\hline\hline
				General Setting & Value \\
				\hline
				Outer loop $K$& 4000 \\
				UL learning rate & 0.005 \\
				LL learning rate & 0.1 \\
				$\mu$ for BDA & 0.1 \\
				$\epsilon$ for CG & $e^{-10}$\\
				$M$ for NS & 40\\
				\hline\hline
			\end{tabular}
		}
	\end{minipage}
	\begin{minipage}[t]{0.2\textwidth}
		\renewcommand\arraystretch{1.2}
		\setlength{\tabcolsep}{1.2mm}{
			\begin{tabular}{ll}
				\hline \hline Specific Parameter & Value \\
				\hline
				$\tau$ & 0.025\\
				$\beta$ & 0.1\\
				$\bar{\mu}$ & 0.9\\
				\hline\hline
			\end{tabular}
		}
	\end{minipage}	
\end{table}

For the BLO problem within the text, we follow the hyper-parameter settings in Tab.~\ref{tab:hyper list0}. All the methods  follow the general setting of hyperparameters, and sl-BAMM follow the instruction of some algorithm-specific hyperparameters. Note that we adopt SGD optimizer for updating UL variables $\x$ and LL variables $\y$. In addition, we present the settings of inner loop $T$ for different experiments in Section 4.1 in Tab.~\ref{tab:hyper list00}. 

\begin{table}[htbp]
	\centering
	\caption{The setting of the inner loop $ T$ for different experiments in Section 4. }			\label{tab:hyper list00}
	\renewcommand\arraystretch{1.2}
	
	\setlength{\tabcolsep}{1.2mm}{
		\begin{tabular}{cccc} 
			\hline \hline
			Figure & RHG & CG/NS & BDA \\
			\hline
			Fig.1 & 1 & 1 & 100 \\
			Fig.2 & - & - & - \\
			Fig.3 & $1 / 5 / 10 / 20$ & 50 & - \\
			Fig.4 & 100 & 100 & -\\
			\hline \hline
		\end{tabular}
	}
\end{table}

In Tab.~\ref{tab:hyper list000}, we also investigate the influence of $\tau$, which is used to remove the restrictive boundedness assumption of $\nabla_\y F(\x, \y)$. Empirically, smaller $\tau$ is better as follows. The numerical experiment is conducted for the toy example in Eq.~(13) by using Strategy S3 in Theorem 3.6 with the stopping error $10^{-7}$.

\begin{table}[htbp]
	\centering
	\caption{Reporting the maximum outer iterations and time to reach the stopping error  as $\tau$ varies.}			\label{tab:hyper list000}
	\renewcommand\arraystretch{1.2}
	
	\setlength{\tabcolsep}{1.2mm}{
		\begin{tabular}{ccc}
			\hline \hline
			$\tau$ & Outer Iteration & Time (s) \\
			\hline
			$1 /  2$ & 10000 & 52.732 \\
			$1 / 4$ & 10000 & 52.801 \\
			$1 / 8$ & 762 & 3.932 \\
			$1 / 10$ & 619 & 3.202 \\
			$1 / 100$ & 352 & 1.836 \\
			\hline \hline
		\end{tabular}
	}
\end{table}

\subsection{Data Hyper-Cleaning}

We use the FashionMNIST dataset for evaluation, which contains different categories of clothing, and serves as a direct drop-in replacement for the original MNIST dataset. The dataset is randomly split to three disjoint subsets, which contain 5000, 5000, 10000 examples, respectively. We report the results of accuracy and F1 score on the test datasets. Following the experimental settings in BDA \cite{liu2020generic}, UL and LL subproblems of data hyper-cleaning problem are defined as follows: $F(\mathbf{x}, \mathbf{y})=\sum_{\left(\mathbf{u}_i, \mathbf{v}_i\right) \in \mathcal{D}_{\text {val }}} \ell\left(\mathbf{y}(\mathbf{x}) ; \mathbf{u}_i, \mathbf{v}_i\right)+\lambda\|\mathbf{y}(\mathbf{x})\|^2$ and $f(\mathbf{x}, \mathbf{y})=\sum_{\left(\mathbf{u}_{i}, \mathbf{v}_{i}\right) \in \mathcal{D}_{\mathtt{tr}}}[\sigma(\mathbf{x})]_{i} \ell\left(\mathbf{y} ; \mathbf{u}_{i}, \mathbf{v}_{i}\right)$,	where $\sigma(\mathbf{x})$ denotes the element-wise sigmoid function to constrain the element in the range of $\left[0,1\right]$, $\mathbf{u}_i$ and $\mathbf{v}_i$ denote the data pairs, $\mathcal{D}_{\mathtt{tr}}$ denotes the training dataset, and $\ell$ denotes the cross-entropy loss function.
We adopt Adam optimizer to update variables $\x$ for all the methods for fair comparison. The values of hyper parameters are listed in Tab.~\ref{tab: hyper list1}.

\begin{table}[htbp]
	\centering
	
	\caption{Values for hyperparameters of data hyper-cleaning.}
	\label{tab: hyper list1}
	\begin{minipage}[t]{0.4\textwidth}
		\renewcommand\arraystretch{1.2}
		\setlength{\tabcolsep}{1.2mm}{
			\begin{tabular}{ll}
				\hline\hline
				General Setting & Value \\
				\hline
				Outer loop $K$& 500 \\
				Inner loop $T$ & 100 \\
				LL learning rate & 0.1\\
				UL learning rate & 0.01 \\
				\hline\hline
			\end{tabular}
		}
	\end{minipage}
	\begin{minipage}[t]{0.2\textwidth}
		\renewcommand\arraystretch{1.2}
		\setlength{\tabcolsep}{1.2mm}{
			\begin{tabular}{ll}
				\hline \hline Specific Parameter & Value \\
				\hline
				$\tau$& 0.001\\
				$\beta$& 0.1 \\
				$\bar{\mu}$&0.9\\
				$p$ & 0.01\\
				\hline\hline
			\end{tabular}
		}
	\end{minipage}
\end{table}

\subsection{Few-Shot Classification}

We employ the ConvNet-4 network structures, which are commonly used in few shot classification tasks. ConvNet-4 is a 4-layer convolutional neural network with $k$ filters followed by batch normalization, non-linearity, and max-pooling operation, then the baseline classifier consists of the fully connected layer with softmax function. UL and LL subproblems are defined as follows:  $F\left(\mathbf{x},\left\{\mathbf{y}^j\right\}\right)=\sum_j \ell\left(\mathbf{x}, \mathbf{y}^j ; \mathcal{D}_{\mathrm{val}}^j\right)$ and $f\left(\mathbf{x},\left\{\mathbf{y}^{j}\right\}\right)=\sum_{j} \ell\left(\mathbf{x}, \mathbf{y}^{j}; \mathcal{D}_{\mathtt{tr}}^{j}\right)$, where $ \mathcal{D}_{\mathtt{tr}}^{j}$ and $ \mathcal{D}_{\mathtt{val}}^{j}$ denotes the meta-training and meta-validation dataset of the $j\textrm{-th}$ task. We adopt 
Adam to update UL variables $\x$ and SGD to update LL variables $\y$ for all the methods.
Related hyperparameters are stated in Table~\ref{tab:few shot list}.

\begin{table}[htbp]
	\centering
	
	\caption{Values for hyperparameters of few-shot classification.}
	\label{tab:few shot list}
	\begin{minipage}[t]{0.4\textwidth}
		\renewcommand\arraystretch{1.2}
		\setlength{\tabcolsep}{1.2mm}{
			\begin{tabular}{ll}
				\hline\hline
				Hyperparameter Setting & Value\\
				\hline
				Outer loop $K$& 20000 \\
				Inner loop $T$& 15 \\
				Learning rate & 0.1 \\
				Meta learning rate & 0.001 \\
				Meta batch size & 16\\
				Hidden size&32\\
				\hline\hline
			\end{tabular}
		}
	\end{minipage}
	\begin{minipage}[t]{0.2\textwidth}
		\renewcommand\arraystretch{1.2}
		\setlength{\tabcolsep}{1.2mm}{
			\begin{tabular}{ll}
				\hline \hline Specific Parameter & Value \\
				\hline
				$\tau$& 0.0001\\
				$\beta$& 0.1 \\
				$\bar{\mu}$&0.7\\
				$p$ & 0.001\\
				\hline\hline
			\end{tabular}
		}
	\end{minipage}
\end{table}

\end{document}